\documentclass[11pt]{article} 
\usepackage{amsfonts,amsmath,amsxtra}
\usepackage{latexsym}
\usepackage{amssymb}

\def\hybrid{\topmargin 0pt      \oddsidemargin 0pt
        \headheight 0pt \headsep 0pt
        \textwidth 16.5cm
        \textheight 23cm
        \marginparwidth 0.0in
        \parskip 5pt plus 1pt   \jot = 1.5ex}
\catcode`\@=11
\def\marginnote#1{}
\newcount\hour
\newcount\minute
\newtoks\amorpm
\hour=\time\divide\hour by60 \minute=\time{\multiply\hour by60
\global\advance\minute by-\hour}
\edef\standardtime{{\ifnum\hour<12 \global\amorpm={am}%
        \else\global\amorpm={pm}\advance\hour by-12 \fi
        \ifnum\hour=0 \hour=12 \fi
      \number\hour:\ifnum\minute<10 0\fi\number\minute\the\amorpm}}
\edef\militarytime{\number\hour:\ifnum\minute<10 0\fi\number\minute}

\def\draftlabel#1{{\@bsphack\if@filesw {\let\thepage\relax
   \xdef\@gtempa{\write\@auxout{\string
      \newlabel{#1}{{\@currentlabel}{\thepage}}}}}\@gtempa
   \if@nobreak \ifvmode\nobreak\fi\fi\fi\@esphack}
        \gdef\@eqnlabel{#1}}
\def\@eqnlabel{}
\def\@vacuum{}
\def\draftmarginnote#1{\marginpar{\raggedright\scriptsize\tt#1}}

\def\draft{\oddsidemargin -0.1truein
        \def\@oddfoot{\sl preliminary draft \hfil
        \rm\thepage\hfil\sl\today\quad\militarytime}
        \let\@evenfoot\@oddfoot \overfullrule 3pt
        \let\label=\draftlabel
        \let\marginnote=\draftmarginnote
\def\@eqnnum{{\rm (\theequation)}
\rlap{\kern\marginparsep\tt\@eqnlabel}%
\global\let\@eqnlabel\@vacuum}  }

\newfont{\Bbbb}{msbm7 scaled 1\@ptsize00}
\newcommand{\zs}{\raise-1pt\hbox{$\mbox{\Bbbb Z}$}}

scaled 1\@ptsize00

\font\sevenmsa=msam6 
\newfam\msafam
\textfont\msafam=\sevenmsa
\def\hexnumber@#1{\ifnum#1<10 \number#1\else
\ifnum#1=10 A\else\ifnum#1=11 B\else\ifnum#1=12 C\else \ifnum#1=13
D\else\ifnum#1=14 E\else\ifnum#1=15 F\fi\fi\fi\fi\fi\fi\fi}
\def\msa@{\hexnumber@\msafam}
\def\llcorner{\delimiter"4\msa@78\msa@78 }
\def\lrcorner{\delimiter"5\msa@79\msa@79 }
\mathchardef\blacktriangleright="3\msa@49
\mathchardef\blacktriangleleft="3\msa@4A \font\tenmsb=msbm10 scaled
1\@ptsize00
\newfam\msbfam
\textfont\msbfam=\tenmsb \scriptfont\msbfam=\tenmsb

\newdimen\Squaresize \Squaresize=14pt
\newdimen\Thickness \Thickness=0.5pt

\def\Square#1{\hbox{\vrule width \Thickness
   \vbox to \Squaresize{\hrule height \Thickness\vss
      \hbox to \Squaresize{\hss#1\hss}
   \vss\hrule height\Thickness}
\unskip\vrule width \Thickness} \kern-\Thickness}

\def\Vsquare#1{\vbox{\Square{$#1$}}\kern-\Thickness}

\def\numberbysection{\@addtoreset{equation}{section}
        \def\theequation{\thesection.\arabic{equation}}}
\numberbysection

\renewcommand{\theequation}{\thesection.\arabic{equation}}
\def\titlepage{\@restonecolfalse\if@twocolumn\@restonecoltrue\onecolumn
     \else \newpage \fi \thispagestyle{empty}\c@page\z@
        \def\thefootnote{\fnsymbol{footnote}} }

\def\endtitlepage{\if@restonecol\twocolumn \else  \fi
        \def\thefootnote{\arabic{footnote}}
        \setcounter{footnote}{0}}  
\relax

\hybrid
\makeatletter
\newdimen\normalarrayskip            
\newdimen\minarrayskip               
\normalarrayskip\baselineskip \minarrayskip\jot
\newif\ifold             \oldtrue            \def\new{\oldfalse}
\def\arraymode{\ifold\relax\else\displaystyle\fi}
\def\eqnumphantom{\phantom{(\theequation)}} 
\def\@arrayskip{\ifold\baselineskip\z@\lineskip\z@
     \else
     \baselineskip\minarrayskip\lineskip1\baselineskip\fi}

\def\@arrayclassz{\ifcase \@lastchclass \@acolampacol \or
\@ampacol \or \or \or \@addamp \or
   \@acolampacol \or \@firstampfalse \@acol \fi
\edef\@preamble{\@preamble
  \ifcase \@chnum
     \hfil$\relax\arraymode\@sharp$\hfil
     \or $\relax\arraymode\@sharp$\hfil
     \or \hfil$\relax\arraymode\@sharp$\fi}}

\def\@array[#1]#2{\setbox\@arstrutbox=\hbox{\vrule
     height\arraystretch \ht\strutbox
     depth\arraystretch \dp\strutbox
width\z@}\@mkpream{#2}\edef\@preamble{\halign \noexpand\@halignto
\bgroup \tabskip\z@ \@arstrut \@preamble \tabskip\z@ \cr}%
\let\@startpbox\@@startpbox \let\@endpbox\@@endpbox
  \if #1t\vtop \else \if#1b\vbox \else \vcenter \fi\fi
  \bgroup \let\par\relax
  \let\@sharp##\let\protect\relax
  \@arrayskip\@preamble}
\def\eqnarray{\stepcounter{equation}%
              \let\@currentlabel=\theequation
              \global\@eqnswtrue
              \global\@eqcnt\z@
              \tabskip\@centering              
              \let\\=\@eqncr
              $$%
            \halign to \displaywidth  \bgroup
             \eqnumphantom \@eqnsel
      \hskip\@centering                               
    $\displaystyle  \tabskip\z@ {##}$%
    &\global\@eqcnt\@ne \hskip 2\arraycolsep
         $ \displaystyle  \arraymode{##}$\hfil
    &\global\@eqcnt\tw@ \hskip 2\arraycolsep
         $\displaystyle\tabskip\z@{##}$\hfil
         \tabskip\@centering
    &{##}\tabskip\z@\cr}
\makeatother
\newcommand{\RR}{{\mathbb{R}}}
\newcommand{\CC}{{\mathbb{C}}}

\def\IC{\mathbb{C}}

\def\IF{\mathbb{F}}

\def\IP{\mathbb{P}}
\def\IQ{\mathbb{Q}}
\def\IR{\mathbb{R}}
\def\IZ{\mathbb{Z}}

\def\CB {\mathcal{B}}
\def\CC {\mathcal{C}}
\def\CD {\mathcal{D}}

\def\CF {\mathcal{F}}

\def\CH {\mathcal{H}}

\def\CL {\mathcal{L}}
\def\CM {\mathcal{M}}

\def\CO {\mathcal{O}}
\def\CP {\mathcal{P}}
\def\CQ {\mathcal{Q}}
\def\CR {\mathcal{R}}

\def\CU {\mathcal{U}}
\def\CV {\mathcal{V}}
\def\CW {\mathcal{W}}

\def\CZ {\mathcal{Z}}

\def\a {{\alpha}}
\def\b {{\beta}}

\def\l {{\lambda}}

\def\la{\lambda}
\def\l{\lambda}

\def\pr {\partial}

\def\wb {\bar{w}}
\def\zb {\bar{z}}

\def\Tr{{\rm Tr}}

\def\frak{\mathfrak}
\def\Fg{{\frak g}}

\newtheorem{te}{Theorem}[section]
\newtheorem{de}{Definition}[section]
\newtheorem{prop}{Proposition}[section]           
\newtheorem{cor}{Corollary}[section]
\newtheorem{lem}{Lemma}[section]
\newtheorem{ex}{Example}[section]
\newtheorem{rem}{Remark}[section]

\newcommand\bqa{\begin{eqnarray}}
\newcommand\eqa{\end{eqnarray}}
\def\be{\begin{eqnarray}\new\begin{array}{cc}}
\def\ee{\end{array}\end{eqnarray}}

\def\beq{\begin{equation}}
\def\eeq{\end{equation}}
\def\bse{\begin{subequations}}                
\def\ese{\end{subequations}}
\def\bp{\begin{pmatrix}}
\def\ep{\end{pmatrix}}

\def\h{\hbar}

\newcommand{\proof}{\noindent {\it Proof}. }
\def\stack#1#2{\raise0.7pt\hbox{$\mathrel{\mathop{#2}\limits^{#1}}$}}
\def\tr{\triangleright}
\def\tl{\triangleleft}
\def\sem{\mathsurround=0pt \raise1pt
\hbox{$\scriptscriptstyle>\!\!$}\:\!\!\tl}
\def\mes{\mathsurround=0pt \tr\!\:\!\raise0.8pt
\hbox{$\scriptscriptstyle\!\!<$}\,}
\def\]{\mathsurround=0pt ]\raise-2pt\hbox{$_\ast$}}

\def\<{\langle}
\def\>{\rangle}

\def\CQ{{\cal Q}}
\def\frak{\mathfrak}

\def\CO{{\cal O}}
\def\CU{{\cal U}}
\def\CZ{{\cal Z}}

\def\CH{\mathcal{H}}

\def\we{\raise-1pt\hbox{$\,\stackrel{\wedge}{,}\,$}}
\def\tr{{\rm tr}\,}
\def\Tr{{\rm Tr}\,}
\def\pr {\partial}

\newcounter{pac}[section]

\newcounter{pacc}[subsection]

0

0

\title{\bf Representation Theory over Tropical Semifield \\
and \\
Langlands Correspondence}

\begin{document}
\author{Anton A. Gerasimov and  Dimitri R. Lebedev}
\date{}

\maketitle

\renewcommand{\abstractname}{}

\begin{abstract}

\noindent {\bf Abstract}. Recently we propose 
a class of
infinite-dimensional integral representations  of classical
$\mathfrak{gl}_{\ell+1}$-Whittaker functions  and local Archimedean
local $L$-factors using two-dimensional topological field theory
framework. The local Archimedean Langlands duality
was identified in this setting with the  mirror symmetry of the underlying
topological field theories. In this note we introduce elementary analogs
of the Whittaker functions and the Archimedean $L$-factors given by
$U_{\ell+1}$-equivariant symplectic volumes of appropriate
K\"{a}hler $U_{\ell+1}$-manifolds.
We demonstrate that the functions thus defined have  a dual description as
matrix elements of representations of monoids $GL_{\ell+1}(\CR)$,
$\CR$ being the tropical semifield.   We also show that
the elementary Whittaker functions
can be  obtained from the  non-Archimedean Whittaker functions
over $\IQ_p$ by taking the formal limit  $p\to 1$.
Hence the elementary special functions constructed in this way might be
considered as functions over  the mysterious field $\IQ_1$.
The existence of two representations
for the elementary Whittaker functions, one as an equivariant volume and
 the other as a matrix element,
should be considered as a manifestation of a hypothetical
elementary analog of the local Langlands duality for number fields.
We would like to note that the elementary local $L$-factors coincide with
$L$-factors introduced previously by Kurokawa.

\end{abstract}
\vspace{5 mm}

\maketitle

\renewcommand{\abstractname}{}

\section{Introduction}

We introduce a simplified version of the local
Archimedean  Langlands correspondence 
as a correspondence between various integral representations
of a class of special functions. Under the local Langlands correspondence
below we will  understood a  duality relation between various representations of
a class of special functions including local $L$-factors and the Whittaker
functions.  This new setup arises 
quite straightforwardly out of an  approach to the
Archimedean Langlands correspondence on the level of special functions
proposed recently \cite{GLO6}, \cite{GLO7},
\cite{GLO8} (see also \cite{G} for a general overview).
 It was argued in op. cit. that the
Langlands duality between two constructions
of the local Archimedean $L$-factors can be considered as an instance
of  a  duality  between infinite-dimensional equivariant
symplectic geometry and finite-dimensional complex geometry
(for generalities on the  Langlands correspondence see e.g.
\cite{ABV}, \cite{B}, \cite{L}).
Thus the proper setting for  Archimedean Langlands duality is the mirror symmetry
of two-dimensional topological field theories. In particular two
dual constructions  of the local Archimedean $L$-factors have a clear
interpretation in terms of mirror dual  two-dimensional topological field theories.

It is known that one of the consequences of Langlands duality
is the existence of two mutually  dual constructions  of  Whittaker
functions (see \cite{Sh}, \cite{CS} for the non-Archimedean case and
\cite{GLO3}, \cite{GLO4}, \cite{GLO5} for the Archimedean case). As  was demonstrated in
\cite{GLO8}, the topological field theory approach can be successfully applied
to description of the duality pattern. 
The Langlands dual constructions
of the $\mathfrak{g}$-Whittaker functions arise while calculating correlation functions
in mirror dual 
topological field theories of type $A$ and type $B$. In type $B$ topological field theory
the correlation function reduces to a finite-dimensional integral
which can be identified with an integral form of a matrix element of an
infinite-dimensional representation of the   Lie algebra $\mathfrak{g}$.
On the other hand the integral
representation arising in type $A$ topological field theory  is given
by an infinite-dimensional integral over the  space of holomorphic maps of
a two-dimensional disk $D$ into a flag space $G^{\vee}/B^{\vee}$
where the Lie group
$G^{\vee}$ is  dual to the Lie group $G$, ${\rm Lie}(G)=\Fg$,
and $B^{\vee}\subset G^{\vee}$ is a Borel subgroup.
In the proposed interpretation
the  infinite-dimensional integral representation
of the Whittaker function is identified with
 the arithmetic side
while the finite-dimensional  integral representation
is identified  (via  relations between matrix elements  of an infinite-dimensional
representation of the dual group)
with the representation theory side of the Langlands correspondence.

In this note we define a simplified version of the Archimedean
Langlands  correspondence for a class of special functions
obtained by replacing two-dimensional topological
field theories by zero-dimensional ones in the constructions of 
\cite{GLO6}, \cite{GLO7},\cite{GLO8}.
 As a result, this simplified version of the Langlands
correspondence is formulated purely
in terms of  finite-dimensional geometry.
 We use the adjective ``elementary'' for constructions arising in this
simplified setting. Note that the two-dimensional topological
field theories involved have an $S^1$-equivariance  parameter $\hbar$
 for which the aforementioned  reduction to
zero dimensions naturally arises in the limit $\hbar\to \infty$.
This parameter can also be easily identified in the standard integral
representations for the Whittaker functions. In this note we mostly
avoid topological field theory
considerations by taking the limit $\hbar \to \infty$ in  the explicit
integral representations for the Whittaker functions \cite{KL},
\cite{Giv}. All the standard properties of the Whittaker functions
have their elementary counterparts. For instance,
the property of $\mathfrak{gl}_{\ell+1}$-Whittaker functions
to be a common eigenfunction of $\mathfrak{gl}_{\ell+1}$-Toda chain
quantum Hamiltonians  is replaced by  the  property of the
elementary  $\mathfrak{gl}_{\ell+1}$-Whittaker function
to be an eigenfunction for a
quantum billiard associated with $\mathfrak{gl}_{\ell+1}$.
We  show that the elementary
 $\mathfrak{gl}_{\ell+1}$-Whittaker function is  the  $U_{\ell+1}$-equivariant
symplectic volume  of the flag space $GL_{\ell+1}/B$.
The dual description  is a matrix element
of a representation of the monoid $GL_{\ell+1}(\CR)$
 where $\CR$ is a tropical semifield. Note that the appearance of the
tropical  monoid in the elementary setting is directly related with the
role  played by the monoid of positive elements of $GL_{\ell+1}$
in the Givental type integral representations of the Whittaker functions
\cite{GKLO}, \cite{GLO1}. Thus the elementary version of
mirror symmetry, relating equivariant  symplectic volumes of
flag spaces $GL_{\ell+1}/B$ and matrix elements of representations
of dual Lie monoids over tropical fields, provides an elementary analog
of  Archimedean Langlands duality.
More generally the elementary Langlands correspondence
for  Whittaker functions
relates symplectic finite-dimensional  geometry of the flag spaces $G/B$
and representation theory of tropical monoids associated with dual
reductive groups $G^{\vee}$.
All the results of \cite{GLO2}, \cite{GLO6}, \cite{GLO7}, \cite{GLO8}
have their counterparts in the elementary setting.
In particular we provide  mutually dual descriptions (i.e. in terms of
finite-dimensional symplectic geometry and in terms of tropical
geometry) of the eigenfunction property of the elementary
$\mathfrak{gl}_{\ell+1}$-Whittaker functions with respect to
an elementary version of the Baxter operators \cite{GLO2}. The corresponding
eigenvalues are given by elementary analogs of local $L$-factors
having both a representation as an equivariant symplectic volume and
as an integral over the tropical semifield.

Construction of the elementary Langlands correspondence  reveals the fundamental
role of  geometry over  tropical semifields as a dual description
of the finite-dimensional symplectic geometry.
The notion of a tropical semifield was first proposed by Maslov as a ``dequantization'' of
real numbers in his study of classical asymptotics of quantum
amplitudes (see \cite{MS} for  detailed explanations) and has
 appeared  (under various  names) in various branches of Mathematics and Physics.
Tropical geometry, i.e. geometry over tropical
semifields, has, for example, been  applied to study
 mirror symmetry of Calabi-Yau manifolds  (the
Berkovich geometry approach due to Kontsevich and Soibelman \cite{KS},
asymptotic analysis of Fukaya \cite{F}), and to  counting complex curves in algebraic
manifolds by Viro, Mikhalkin et al (see e.g. \cite{IMS}).
See also  \cite{GKZ}, \cite{Mi}, \cite{EKL} for discussions of
tropical geometry as a limit of real geometry in the case of toric
manifolds.

The  meaning of the elementary Whittaker  functions can be partially
elucidated using the $q$-deformed  Whittaker functions 
\cite{GLO3}, \cite{GLO4}, \cite{GLO5}.
These functions provide an interpolation between  classical
$\mathfrak{gl}_{\ell+1}$-Whittaker  functions and their  non-Archimedean
analogs. The elementary $\mathfrak{g}_{\ell+1}$-Whittaker functions
can be obtained from  the $q$-deformed $\mathfrak{gl}_{\ell+1}$-Whittaker
functions in two different ways.
One can first consider a limit corresponding to  classical
$\mathfrak{gl}_{\ell+1}$-Whittaker function \cite{GLO9} and then obtain, by
further degeneration,  the elementary $\mathfrak{gl}_{\ell+1}$-Whittaker
functions. Equivalently one can first take a limit of the $q$-deformed
$\mathfrak{gl}_{\ell+1}$-Whittaker functions  leading to the non-Archimedean
$\mathfrak{gl}_{\ell+1}$-Whittaker functions over $\IQ_p$. The
elementary $\mathfrak{gl}_{\ell+1}$-Whittaker functions
are then  obtained by taking a formal limit $p\to 1$.
The last limit has a simple explanation based on
the Shintani-Casselman-Shalika formula \cite{Sh}, \cite{CS}
expressing the  non-Archimedean  Whittaker functions in terms of characters of
finite-dimensional irreducible representations of the dual group.
According to  the Kirillov philosophy (see e.g. \cite{K}),
characters of finite-dimensional irreducible representations of
compact Lie groups  admit a limit expressed as an integral over
corresponding coadjoint orbit.
 This limiting representation is precisely the
expression of the elementary Whittaker functions as
equivariant symplectic volume integrals
obtained in the limit $p\to 1$. Taking the limit $p\to 1$ as a way to
produce the elementary analogs implies that  tropical geometry
should be considered as an effective description of  geometry over a
 mysterious field $\IQ_1$. It is known that
amoebas, defined in \cite{GKZ},  for variety over non-archimedean fields
are identical to tropicalizations of the varieties over complex
numbers \cite{Mi}, \cite{EKL}.
Recently in \cite{CC}, among other things,  the  limit of $\IQ_p$ for $p\to 1$ was
also discussed and the relation with  tropical semifields is stressed.
Although there is an obvious similarity with our considerations,  we should
note that in contrast to \cite{CC},
we treat the tropical semifield as a universal target of the
valuation map while $\IQ_1$ being  surjectively mapped (as
multiplicative monoid) 
onto tropical semifield $\CR$ can have a larger kernel. We also should note that elementary $L$-factors
coincide with the $L$-factors introduced by Kurokawa \cite{Ku} (see also
\cite{Ma}).

Let us note that the  relation between quantum billiards and equivariant symplectic
volumes  is a consequence of the general
description of the equivariant cohomology of flag manifolds
in terms of  representation theory of  nil-Hecke algebras
\cite{BGG}, \cite{KK} (see also \cite{CG}, \cite{Gzb} for  reviews).
This points to a hierarchy of generalizations of the
results of \cite{GLO6}, \cite{GLO7}, \cite{GLO8}, and of the present article
associated with multidimensional
analogs of Hecke algebras (for the theory of double Hecke algebras see
\cite{Ch}). As a simple straightforward generalization
one can consider elementary analogs of
spherical functions and their connection with  representation 
theory of graded affine Hecke algebras. Thus
elementary analogs of $GL_{\ell+1}$-spherical
functions have  two dual descriptions. On the one hand they can be expressed
as  $S^1\times U_{\ell+1}$-equivariant symplectic  volumes
of cotangent bundles $T^*\CB$ where $S^1$ acts on the fibres of the
projection $T^*\CB\to \CB$. On the other hand elementary
spherical functions can be realized canonically as spherical functions
on the tropical monoid $G^{\vee}(\CR)$. An elementary analog of
the Calogero-Sutherland integrable systems (i.e. an integrable system
such that spherical functions are common eigenfunctions of the
corresponding quantum Hamiltonians) are systems of quantum
particles with delta-function interactions (also known as Yang
systems, see e.g. \cite{HO}).  The Hecke algebra description, generalizing the
discussion in Section 6, is given in terms of the
 $G\times \IC^*$-equivariant  K-theory $K_G(\CB)$ and the
graded affine  Hecke algebra $\CH(G)$.
Many other examples of these constructions are known (see e.g. \cite{CG}).
 We defer detailed discussions of the case of spherical functions 
and other generalizations for  another occasion.

In this note we introduce an elementary analog of a particular
manifestation of the local Langlands
correspondence expressed as  a relation between various representations of
special functions e.g. the Whittaker functions.
We expect that this can be generalized to an elementary analog of the
full-flagged local  Langlands
correspondence between admissible representations
of local reductive groups $G(K)$
and admissible representations of the
Weil-Deligne group $W_K$ factored through the homomorphism of 
$W_K$ into dual reductive group ${}^LG$. Let us stress that the proposed
elementary version of the local Langlands correspondence
might be useful for understanding  classical Langlands correspondence.
In particular one can expect  Langlands functoriality
to be more accessible in the elementary setting.

The plan of the paper is as follows. In Section 2 we recall
relevant facts about classical $\mathfrak{gl}_{\ell+1}$-Whittaker functions,
local Archimedean $L$-factors and  Baxter  integral operators.
We also briefly discuss  a topological field theory interpretation of
a class of  Whittaker functions and local Archimedean $L$-factors
following \cite{GLO6}, \cite{GLO7}, \cite{GLO8}. In
Section 3 we define elementary analogs of $\mathfrak{gl}_{\ell+1}$-Whittaker functions,
local Archimedean $L$-factors and of Baxter integral operators.
It is argued that  the  quantum billiard associated with the root system of
$\mathfrak{gl}_{\ell+1}$ plays the role of  an elementary analog of the
$\mathfrak{gl}_{\ell+1}$-Toda chain.
In Section 4 we demonstrate that elementary Whittaker functions
can be obtained as a limit of a specialization at $q=0$
of the  $q$-deformed $\mathfrak{gl}_{\ell+1}$-Whittaker functions
introduced in \cite{GLO3}, \cite{GLO4}, \cite{GLO5}.
In Section 5  elementary analogs of classical functions are identified
with equivariant symplectic volumes. Thus
 the $\mathfrak{gl}_{\ell+1}$-Whittaker functions
are identified with $U_{\ell+1}$-equivariant symplectic volumes of the
flag spaces  $GL_{\ell+1}(\IC)/B$ and
elementary  $L$-factors associated with standard representations
$\IC^{\ell+1}$ of $\mathfrak{gl}_{\ell+1}$
are identified with $U_{\ell+1}$-equivariant symplectic volumes of
$\IC^{\ell+1}$.  We also provide a symplectic geometry interpretation
of the eigenfunction  property of an elementary Whittaker function with
respect to an elementary analog of the Baxter operator.  In Section 6
we elucidate some of the previous construction  using
 the Kostant-Kumar  description \cite{KK} of equivariant cohomology of flag
spaces in terms of representation theory of  nil-Hecke
algebras. In Section 7 the dual description of the elementary
$\mathfrak{gl}_{\ell+1}$-Whittaker functions as matrix elements of
infinite-dimensional representations of tropical monoids associated with
$GL_{\ell+1}$ is given (Theorem \ref{Theorem}).
 Finally, in Section 8 we briefly discuss a relation
of our constructions with  geometry over the mysterious limiting field $\IQ_1$.

{\em Acknowledgments}: The research was supported by  Grant
RFBR-09-01-93108-NCNIL-a. The research of AG was  also partly
supported by Science Foundation Ireland grant.  The authors are thankful to
Max-Planck-Institut f\"ur Mathematik in Bonn for hospitality and
excellent working conditions.

\section{Whittaker functions in classical setting}

In \cite{GLO6}, \cite{GLO7}, \cite{GLO8} 
we introduced, using a topological field theory framework,
infinite-dimensional integral representations of the
local Archimedean $L$-factors (given by products of Gamma-functions)
and of a  class of the  Whittaker functions. These integral representations
can be reduced to give representations of  the special
functions as equivariant volumes of the infinite-dimensional spaces of
holomorphic maps of a two-dimensional disk into a K\"{a}hler manifold
such as a flag space  $G/P$ or a complex vector space $V=\IC^{\ell+1}$.
This should be compared with the  classical
finite-dimensional integral representations  of the  same special
functions e.g. the Euler integral representation of
the Gamma-function,  the Mellin-Barnes and the Givental integral
representations of the Whittaker functions \cite{KL}, \cite{GKLO}.
It was argued in \cite{GLO7}, \cite{GLO8}
that both  the finite dimensional and the infinite-dimensional integral
representations arise naturally as correlation functions in two-dimensional
topological field theories and
are related by a mirror symmetry of the underlying
quantum field theories. One can  conjecture that
 analogs of infinite-dimensional integral representations hold for
many other special functions such as spherical functions and the
Whittaker functions
associated with  arbitrary pairs $(\Fg, \mathfrak{p})$ of  Lie
algebras $\Fg$ and their parabolic subalgebras $\mathfrak{p}$.

In this Section we review two particular classes of integral
representations of the $\mathfrak{gl}_{\ell+1}$-Whittaker functions
associated with maximal and minimal parabolic subalgebras.
 We provide a realizations of the
$\mathfrak{gl}_{\ell+1}$-Whittaker functions
associated  with  maximal parabolic subgroups
as equivariant volumes of spaces of holomorphic maps \cite{GLO8},
 recall two integral representations  of the local Archimedean $L$-factors and
its  topological field theory interpretation \cite{GLO6}, \cite{GLO7}.

Let $E_{ij}$, $i,j=1,\ldots \ell+1$ be the standard basis of the Lie
algebra $\mathfrak{gl}_{\ell+1}$.
Let $\CZ(\CU\mathfrak{gl}_{\ell+1})\subset
\CU\mathfrak{gl}_{\ell+1}$ be the  center of the  universal
enveloping algebra $\CU\mathfrak{gl}_{\ell+1}$.
Let  $B_{\pm}\subset GL_{\ell+1}(\IC)$
be  upper-triangular and lower-triangular
Borel subgroups and  $N_{\pm}\subset B_{\pm}$
be  upper-triangular and lower-triangular
unipotent subgroups. We denote by $\mathfrak{b}_{\pm}={\rm Lie}(B_{\pm})$ and
$\mathfrak{n}_{\pm}={\rm Lie}(N_{\pm})$ their Lie algebras.
 Let $\mathfrak{h}\subset \mathfrak{gl}_{\ell+1}$
be the  diagonal Cartan subalgebra and $W=\mathfrak{S}_{\ell+1}$ be the Weyl
group of $GL_{\ell+1}$.  Using the Harish-Chandra
isomorphism between $\CZ(\CU\mathfrak{gl}_{\ell+1})$ and the
$\mathfrak{S}_{\ell+1}$-invariant subalgebra
of the symmetric algebra $S^*\mathfrak{h}$,  we  identify central characters with
homomorphisms $c:\IC[h_1,\cdots, h_{\ell+1}]^{S_{\ell+1}}\to \IC$. 
The central characters are in one to one correspondence 
with $S_{\ell+1}$-orbits of elements of the dual space
$\mathfrak{h}^*\simeq \IC^{\ell+1}$.
 Let $\pi_{\underline{\lambda}}:\CU\mathfrak{gl}_{\ell+1}\to
{\rm End}(\CV_{\underline{\la}})$  be a representation
of the universal enveloping algebra $\CU\mathfrak{gl}_{\ell+1}$
with  a central character associated with the $S_{\ell+1}$-orbit of 
$\imath \hbar^{-1}\underline{\lambda}=(\imath
\hbar^{-1}\lambda_1,\ldots, \imath \hbar^{-1}\lambda_{\ell+1})\in
\imath \IR^{\ell+1}$,  $\hbar\in \IR_+$. We impose an additional
condition that
the action of the Cartan subalgebra $\mathfrak{h}$
can be integrated
to an action of the corresponding Cartan subgroup $H\subset
GL_{\ell+1}(\IC)$ in $\CV_{\underline{\la}}$.
 Let $\CV'_{\underline{\lambda}}$ be the dual module
equipped with the induced action of $\CU\mathfrak{gl}_{\ell+1}^{opp}$
(universal enveloping algebra of $\mathfrak{gl}_{\ell+1}^{opp}$
obtained by inverting the signs of the structure constants of
$\mathfrak{gl}_{\ell+1}$). Denote by $\<\,,\,\>$  the pairing
between $\CV'_{\underline{\lambda}}$ and $\CV_{\underline{\lambda}}$.
We suppose that the action of the Cartan subalgebra $\mathfrak{h}$
 in the representation  $\CV_{\underline{\lambda}}$ can be integrated
to an action of the corresponding Cartan subgroup $H\subset GL_{\ell+1}(\IC)$.

According to Kostant a $\mathfrak{gl}_{\ell+1}$-Whittaker function
is  defined as  a  matrix element
\be\label{Wf}
\Psi_{\underline{\lambda}}(\underline{x})\,=\,e^{\rho(\underline{x})}
\<\psi_L|\,\pi_{\underline{\lambda}}(e^{\sum_{i=1}x_iE_{ii}})\,|\psi_R\>,
\ee
where $\underline{x}=(x_1,\ldots, x_{\ell+1})$, 
$\rho=(\rho_1,\ldots ,\rho_{\ell+1})$ is  a vector in $\mathfrak{h}^*$ with 
components $\rho_k=-1/2(\ell-2k+2),\,k=1,\ldots,\ell+1$ 
(half of the sum of positive roots of
$\mathfrak{gl}_{\ell+1}$).  The  one-dimensional spaces generated by 
vectors $\<\psi_L|\in \CV'_{\underline{\lambda}}$ and
$|\psi_R\>\in \CV_{\underline{\lambda}}$
provide  one-dimensional representations
of $\CU\mathfrak{n}_-$ and $\CU\mathfrak{n}_+$ respectively
\be\label{Wclass}
\<\psi_L|E_{i+1,i}=-\hbar^{-1}\<\psi_L|,\qquad
E_{i,i+1}|\psi_R\>=-\hbar^{-1}|\psi_R\>,\qquad i=1,\ldots ,\ell.
\ee

Standard considerations (see e.g. \cite{STS})
show that  the matrix element \eqref{Wf}  is a common
eigenfunction  of a family of commuting differential operators
coming from the action of  generators of
$\CZ(\CU\mathfrak{gl}_{\ell+1})$  in $\CV_{\underline{\lambda}}$.
These differential operators  can be identified with quantum Hamiltonians of the
$\mathfrak{gl}_{\ell+1}$-Toda chain. The simplest non-trivial quantum
Hamiltonian acts on a Whittaker function as the differential operator
$$
\CH=-\frac{\hbar^2}{2}\sum_{i=1}^{\ell+1}\frac{\pr^2}{\pr
  x_i^2}+\sum_{i=1}^{\ell}e^{x_{i+1}-x_{i}}.
$$
The matrix element representation  \eqref{Wf}
leads to various  integral representations of
the Whittaker function by using explicit realizations of the
universal enveloping algebra representation
$\pi_{\underline{\lambda}}$
via difference/differential operators acting on an appropriate
space of functions. For example one can consider
the principal series
representation ${\rm   Ind}_{B_-}^{GL_{\ell+1}}
\chi_{\underline{\lambda}}$  induced from  a one-dimensional
representation of $B_-$ given by the character 
\be\label{chardef}
\chi_{\underline{\lambda}}(b)=\prod_{j=1}^{\ell+1}\,|b_{jj}|^{
\imath\hbar^{-1}\lambda_k-\rho_k}, \qquad b\in B_-.
\ee
Here
 $\underline{\lambda}=(\lambda_1,\ldots, \lambda_{\ell+1})$ is a vector in $
\IR^{\ell+1}$,  $\hbar$ is in $\IR_+$  and $\rho_k=-1/2(\ell-2k+2),\,\,\,\,
k=1,\ldots,\ell+1$.   
The subspace of analytic vectors  provides
 the corresponding representation of $\CU\mathfrak{gl}_{\ell+1}$.
More generally one can consider  representations of
$\CU\mathfrak{gl}_{\ell+1}$ realized in the subspace of
$B_-$-equivariant functions supported on $B_-$-stable subvarieties 
in  $GL_{\ell+1}$.  They
obviously support an action of $\CU\mathfrak{gl}_{\ell+1}$ with central
character associated to the $S_{\ell+1}$-orbit through $\imath \hbar^{-1}\underline{\lambda}$.

In this article we shall be especially interested in the integral 
representations of $\mathfrak{gl}_{\ell+1}$-Whittaker functions defined in the following two theorems.

\begin{te} The  $\mathfrak{gl}_{\ell+1}$-Whittaker function has the
  following integral representation:
\be\label{levoneW} \Psi_{\underline{\la}}(\underline{x},\hbar)=\int\limits_{\cal
S}\prod_{n=1}^{\ell}
\frac{\prod\limits_{k=1}^{n+1}
\prod\limits_{m=1}^{n}
\Gamma_1\left(\frac{\gamma_{n+1,k}-\gamma_{n,m}}{\imath}|\hbar\right)} {
n!\prod\limits_{s\neq p} \Gamma_1\left(\frac{\gamma_{ns}-
\gamma_{np}}{\imath}|\hbar\right)} e^{\frac{\imath}{\hbar} \sum\limits_{n=1}^{\ell+1}\,
\sum\limits_{j=1}^{n} (\gamma_{nj}-\gamma_{n-1,j})x_n}
\prod_{\stackrel{\scriptstyle n=1}{j\leq n}}^{\ell}\frac{d\gamma_{nj}}{2\pi\hbar}\,,
\ee
where $\Gamma_1(z|\hbar)=\hbar^{\frac{z}{\hbar}}\Gamma(\frac{z}{\hbar})$,
$\underline{\la}=(\la_1,\ldots,\la_{\ell+1}):=(\gamma_{\ell+1,1},\ldots,
\gamma_{\ell+1,\ell+1})\in \IR^{\ell+1}$, $\underline{x}=(x_1,\ldots,x_{\ell+1})$
and the domain of  integration ${\cal S}$  is such that
$\gamma_{kj}+\imath \epsilon_{jk}\in \IR$ and $\epsilon_{jk}$ are fixed
real numbers  such that
 the conditions  $\!\max_{j}\{{\rm Im}\,\gamma_{kj}\}<
\min_m\{{\rm Im}\,\gamma_{k+1,m}\}\!$ for all $ k=1,\ldots,\ell$  hold.
The integral (\ref{levoneW}) converges absolutely.
Recall that we  assume   $\gamma_{nj}=0$ for $j>n$.
\end{te}

This integral representation was
first  introduced in \cite{KL} and rederived
in the framework of representation theory in \cite{GKL}.
Another relevant integral representation
for  the $\mathfrak{gl}_{\ell+1}$-Whittaker function
was introduced by Givental  \cite{Giv}
(see \cite{GKLO} for a representation theoretic derivation).

\begin{te} The $\mathfrak{gl}_{\ell+1}$-Whittaker function has the following integral representation:
\be\label{giv}
 \Psi_{\underline\lambda} (\underline{x}|\hbar)= \int\limits_{\cal{C}}
\exp\left(\frac{\imath}{\hbar}\sum\limits_{k=1}^{\ell+1}
\lambda_k \left(\sum_{i=1}^{k}T_{k,i}-\sum_{i=1}^{k-1}T_{k-1,i}\right)\right)\\
\times
\exp\{\sum_{k=1}^{\ell}-\frac{1}{\hbar}\left(\sum_{i=1}^{k}
e^{T_{ki}-T_{k+1,i}}+
\sum_{i=1}^{k}
e^{T_{k+1,i+1}-T_{k,i}}\right)\}\prod_{k=1}^{\ell}\prod_{i=1}^{k}dT_{k,i}\,,
\ee
where $\underline{\la}=(\lambda_1,\ldots,\lambda_{\ell+1})\in
\IR^{\ell+1}$. Here  we set $x_i=T_{\ell+1,i},\,\,\,\,i=1,\ldots,\ell+1$  
and the domain of integration $\cal{C}$ is a slight deformation of the 
middle-dimensional real subspace of $\IC^{\ell(\ell+1)/2}$
rendering  the integral (\ref{giv})  convergent.
\end{te}

In \cite{GLO8} we introduced generalized $\mathfrak{gl}_{\ell+1}$-Whittaker functions
associated
with  parabolic subalgebras $\mathfrak{p}\subset
\mathfrak{gl}_{\ell+1}$ (in this sense the standard
$\mathfrak{gl}_{\ell+1}$-Whittaker function \eqref{Wf} is associated with the
Borel subalgebra $\mathfrak{b}\subset
\mathfrak{gl}_{\ell+1}$).  For the  $\mathfrak{gl}_{\ell+1}$-Whittaker
function associated with the maximal parabolic subalgebra
$\mathfrak{p}_{1,\ell+1}\subset
\mathfrak{gl}_{\ell+1}$ (i.e. such that $\dim(\mathfrak{gl}_{\ell+1}/
\mathfrak{p}_{1,\ell+1})=\ell$) we coin  the term $(1,\ell+1)$- parabolic
Whittaker function. The following
analogs of the Mellin-Barnes and the Givental integral representations
of a specialization of the $(1,\ell+1)$- parabolic
Whittaker function holds.

\begin{te}
The $(1,\ell+1)$-parabolic
Whittaker function specialized to $x_1=x$ and $x_i=0$, $i\neq 1$
admits the following integral representations:

(1).
\be\label{identone}
\Psi^{(1,\ell+1)}_{\underline{\lambda}}(x,0,\cdots ,0)\,=
\,\frac{1}{2\pi\hbar}\int_{-\infty-\imath\epsilon}^{\infty-\imath\epsilon}
 \,\,\,d \gamma \,\,
e^{\frac{\imath}{\h} x \gamma}\,\prod_{j=1}^{\ell+1}
\Gamma_1\Big(\imath (\gamma-\lambda_j)|\hbar\Big),
\ee
where $\epsilon > 0$.

(2).
\be\label{x=0one}
\Psi^{(1,\ell+1)}_{\underline{\lambda}}(x,0,\cdots ,0)\,=\,\int_{\CC}
\prod_{j=1}^{\ell}dt_j\,\, \exp \left(\frac{\imath}{\h}\left(\sum_{j=1}^{\ell}\lambda_j
 t_j +\lambda_{\ell+1}(x-\sum_{j=1}^{\ell}t_j)\,\right)\right) \,\cdot\\
  \exp\left(-\frac{1}{\hbar}\left(
  \sum_{j=1}^{\ell}e^{-t_j}+
  e^{\sum_{j=1}^{\ell}t_j-x}\right)\right),
  \ee
 where $\underline{\la}=(\lambda_1,\ldots,\lambda_{\ell+1})\in
 \IR^{\ell+1}$ and  ${\cal C} $ is a slight  deformation of the 
real subspace of middle dimension in $\IC^{\ell}$   rendering 
the integral \eqref{x=0one} convergent.
\end{te}

Whilst being a  common
eigenfunction of a family of mutually commuting differential
operators, 
the $\mathfrak{gl}_{\ell+1}$-Whittaker function
is also a common eigenfunction of a one-parameter family of
integral operators  \cite{GLO2}. This family of integral
operators (called the Baxter operators due to their relation with the
Baxter operators in the theory of quantum integrable systems)
has a representation  theoretic origin. In \cite{GLO2}
 the Baxter integral operators were defined as
generating series for the generators of the Archimedean counterpart of
the non-Archimedean spherical Hecke algebra. Explicitly,
a one-dimensional family $\mathcal{Q}^{\mathfrak{gl}_{\ell+1}}(s),\,\,\,s\in \IC$
of the Baxter integral operators acting in
an appropriate space of functions of $\ell+1$ variables
has the following integral kernel:
\be\label{Baxter}
\mathcal{Q}^{\mathfrak{gl}_{\ell+1}}
(\underline{x},\,\underline{y}|\,s)=
\exp\Big\{\,\frac{\imath}{\hbar} s \sum_{i=1}^{\ell+1}(x_i-y_i)\,-
\frac{1}{\hbar}\sum_{k=1}^{\ell}\Big(e^{x_i-y_i}+e^{y_{i+1}-x_{i}}\Big)\,-\,
\frac{1}{\hbar} e^{x_{\ell+1}-y_{\ell+1}} \Big\}.\ee
The Baxter operator  $\mathcal{Q}^{\mathfrak{gl}_{\ell+1}}(s)$
satisfies the following commutativity relations:
\bqa\label{firstpr}
\mathcal{Q}^{\mathfrak{gl}_{\ell+1}}(s)\cdot
 \mathcal{Q}^{\mathfrak{gl}_{\ell+1}}(s')=
\mathcal{Q}^{\mathfrak{gl}_{\ell+1}}(s')\cdot
\mathcal{Q}^{\mathfrak{gl}_{\ell+1}} (s), \eqa
\bqa\label{secondpr}
\mathcal{Q}^{\mathfrak{gl}_{\ell+1}}(s)\cdot
\CH_r^{\mathfrak{gl}_{\ell+1}}=
\CH_r^{\mathfrak{gl}_{\ell+1}}\cdot
\mathcal{Q}^{\mathfrak{gl}_{\ell+1}}(s),\qquad r=1,\ldots \ell+1,
\eqa
and the $\mathfrak{gl}_{\ell+1}$-Whittaker function \eqref{levoneW}, \eqref{giv}
satisfies the following eigenfunction identity:
\bqa\label{eigenprop}
\int_{\RR^{\ell+1}}\,\prod_{i=1}^{\ell+1}\,dy_{i}\,\,
\mathcal{Q}^{\mathfrak{gl}_{\ell+1}}(\underline{x},
\,\underline{y}|\,s)\,
\Psi^{\mathfrak{gl}_{\ell+1}}_{\underline{\lambda}}(\underline{y})\,=\,
L(s,\underline{\lambda},\hbar)\,\,
\Psi^{\mathfrak{gl}_{\ell+1}}_{\underline{\lambda}}(\underline{x}).\eqa
 Here  $\underline{x}=(x_1,\ldots, x_{\ell+1})$,
$\underline{y}=(y_1,\ldots, y_{\ell+1})$,
$\underline{\la}=(\la_1,\ldots, \la_{\ell+1})$ 
and the eigenvalue is  given by the local Archimedean $L$-factor
attached to the principal series representation of $GL_{\ell+1}$
associated with the character \eqref{chardef}
 (see e.g. \cite{B}, \cite{L})
\be\label{lAL}
L(s,\underline{\lambda},\hbar)=\prod_{j=1}^{\ell+1}\,
\hbar^{\frac{ \imath s - \imath \lambda_j}{\hbar}}
\Gamma\Big(\frac{\imath s- \imath \lambda_{j}}{\hbar}\Big).
\ee
The appearance of the local $L$-factors is not accidental and is
related to the fact that  the integral operators with kernels \eqref{Baxter}
are realizations of generating functions of elements of
 the local spherical Archimedean Hecke algebra (see \cite{GLO2} for a
 detailed discussion). 
Note that the local Archimedean $L$-factors \eqref{lAL} also have integral
representations of the Givental type \eqref{giv}
given by the Euler integral representations for the Gamma-functions
\be\label{Lgiv}
L(s,\underline{\lambda},\hbar)=\int_{\IR^{\ell+1}} \prod_{j=1}^{\ell+1}\,d\tau_j\,
e^{\sum_{j=1}^{\ell+1} \left(\frac{\imath}{\hbar}(s-\lambda_j)\tau_j-
\frac{1}{\hbar}e^{\tau_j}\right)},
\ee
where ${\rm Im}(s)<0 ,\,\,j=1,\ldots,\ell+1$.

The expression of the local $L$-factor \eqref{lAL} as a product of
$\Gamma$-functions is on the other hand an analog of the Mellin-Barnes representation
\eqref{levoneW} of the $\mathfrak{gl}_{\ell+1}$-Whittaker function.

Along with the finite-dimensional integral representations, the Whittaker functions and
the local Archimedean $L$-factors have  infinite-dimensional integral
representations. These integral representations  naturally arise from
an interpretation  of the special functions  as equivariant symplectic volumes
of infinite-dimensional space of holomorphic maps
of a two-dimensional disk into  finite-dimensional symplectic
spaces \cite{GLO6}, \cite{GLO7}, \cite{GLO8}.
Given a finite-dimensional symplectic manifold $M$ with a symplectic form $\omega$ and
Hamiltonian action of a compact Lie group $G$ one defines a
$G$-equivariant symplectic volume as the following integral:
\be\label{eqvolL}
Z_M(\lambda)=\int_M\,e^{\omega_G}.
\ee
Here $\omega_G$ is the  $G$-equivariant extension  of the symplectic form
$\omega$ depending on an element $\lambda\in \Fg^*$ of the dual space to the
Lie algebra $\Fg={\rm Lie}(G)$ (see e.g. \cite{Au} for  precise definitions and
details). In \cite{GLO6}, \cite{GLO7}, \cite{GLO8}
this definition was used for infinite-dimensional
spaces of holomorphic maps $\CM(D,X)$ of a two-dimensional
disk $D=\{z\in \IC|\,|z|\leq 1\}$
into  symplectic $U_{\ell+1}$-spaces $X$ such as $\IC^{\ell+1}$, $SL_{\ell+1}/P$,
with $P\subset SL_{\ell+1}$ being a parabolic group. There is a natural
Hamiltonian action of $S^1\times U_{\ell+1}$ on $\CM(D,X)$, where $S^1$
acts by rotations of $D$ and the action of $U_{\ell+1}$ is induced
from the action on $X$. An extension of the definition   \eqref{eqvolL}
to the infinite-dimensional case requires some care and includes a
$\zeta$-function regularization of the infinite-dimensional
integrals.

The equivariant symplectic volumes of the spaces $\CM(D,X)$ of holomorphic maps
can  be reformulated as  particular correlation functions
in  two-dimensional equivariant
topological sigma models on the disk $D$ with  target space $X$. An
 advantage of this reformulation is that one can invoke
mirror symmetry considerations to reformulate the infinite-dimensional
integral representations  in terms of finite-dimensional ones
\cite{GLO6}, \cite{GLO7}, \cite{GLO8}.
For the  local Archimedean $L$-factors and the
$\mathfrak{gl}_{\ell+1}$-Whittaker functions
this leads to the integral representations \eqref{Lgiv}, \eqref{giv}.

Note that $S^1\times U_{\ell+1}$-equivariant
topological field theory depends on
$\lambda\in \mathfrak{u}_{\ell+1}^*$ and $\hbar\in {\rm Lie}(S^1)$.
In the limit $\hbar\to \infty$ one expects that the
$S^1\times U_{\ell+1}$-equivariant symplectic volume
of $\CM(D,X)$ will be reduced to a $U_{\ell+1}$-equivariant volume integral
over $X$.  In the following
Section we take this limit in the integral representations
 \eqref{levoneW}, \eqref{giv}, \eqref{Lgiv} for the $\mathfrak{gl}_{\ell+1}$-Whittaker
functions  and  local $L$-factors directly and demonstrate
that the resulting elementary Whittaker functions and
$L$-factors are given by equivariant symplectic volumes  of
finite-dimensional symplectic spaces.
Note that up
to now the interpretation of (parabolic) Whittaker functions in terms
of equivariant symplectic volumes of infinite-dimensional spaces is
not known in full generality (recent progress has however been presented in \cite{O1}, \cite{O2}).
In cases when the interpretation in terms of topological field theory with a target
space $X$ is already
established \cite{GLO8}
 the obtained limiting expression is compatible with the localization
to the space of constant  maps to $X$ discussed above.
Thus the identification of the $\hbar\to \infty$ limits of the
$\mathfrak{gl}_{\ell+1}$-Whittaker
functions with $U_{\ell+1}$-equivariant symplectic volumes of
the flag spaces $G/B$ should be considered as  additional support
for the approach  of \cite{GLO6}, \cite{GLO7}, \cite{GLO8}.

\section{ Elementary analogs of  Whittaker functions and local L-factors}

In this Section we take the limit $\hbar\to \infty$ of the $\mathfrak{gl}_{\ell+1}$-
Whittaker functions and the local Archimedean $L$-factors and demonstrate
that the resulting expression can be interpreted as an equivariant
symplectic volume of finite-dimensional K\"{a}hler spaces.
Let us start by introducing  elementary analogs of the
$\mathfrak{gl}_{\ell+1}$-Whittaker functions.

\begin{de} The elementary $\mathfrak{gl}_{\ell+1}$-Whittaker function
is defined as the following absolutely convergent integral:
\be\label{levzeroW} \Psi^{(0)}
_{\underline{\la}}(\underline{x})=\imath^{-\frac{\ell(\ell+1)}{2}}\int\limits_{{\cal
S}}\prod_{n=1}^{\ell}
\frac {\prod\limits_{s\neq p} (\gamma_{ns}-\gamma_{np})}
{n!\prod\limits_{k=1}^n\prod\limits_{m=1}^{n+1}
(\gamma_{n k}-\gamma_{n+1,m})}
 e^{\imath\sum\limits_{n=1}^{\ell+1}\,
\sum\limits_{j=1}^{n} (\gamma_{nj}-\gamma_{n-1,j})x_n}
\prod_{\stackrel{\scriptstyle n=1}{j\leq n}}^{\ell}\frac{d\gamma_{nj}}{2\pi\imath}\,,
\ee
where the domain of  integration  $\!{\cal S}\!$ is such that
$\gamma_{jk}+\imath \epsilon_{jk}\in \IR$ and $\epsilon_{jk}$ are
fixed real numbers such that  the
 conditions  $\!\max_{j}\{{\rm Im}\,\gamma_{kj}\}<
\min_m\{{\rm Im}\,\gamma_{k+1,m}\}\!$ for all $ k=1,\ldots,\ell$ hold.
We set $(\gamma_{\ell+1,1},\ldots ,\gamma_{\ell+1,\ell+1})
:=(\lambda_1,\ldots ,\lambda_{\ell+1})\in \IR^{\ell+1}$
and assume   $\gamma_{nj}=0$ for $j>n$. We also use the notations
$\underline{x}=(x_1,\ldots, x_{\ell+1})$ and
$\underline{\lambda}=(\lambda_1,\ldots ,\lambda_{\ell+1})$.
\end{de}

\begin{prop} The elementary  $\mathfrak{gl}_{\ell+1}$-Whittaker function
has the following representation:
\be\label{sumover}
\Psi^{(0)}_{\underline{\lambda}}(\underline{x})\,=\,
\sum_{s\in \mathfrak{S}_{\ell+1}}
 (-1)^{l(s)}\frac{e^{\sum_{k=1}^{\ell+1}\imath\lambda_{s(k)}x_{k}}}
{\prod_{i<j}\imath(\lambda_{i}-\lambda_{j})},
\qquad x_1\geq x_2\geq x_3\geq \ldots \geq x_{\ell+1},
\ee
and
\be\label{sumoverone}
\Psi^{(0)}_{\underline{\lambda}}(\underline{x})\,=\,0,
\ee
otherwise. Here $(s(1), \ldots, s(\ell+1))$ is a permutation
of $(1,\ldots, \ell+1)$ corresponding to an element $s\in
\mathfrak{S}_{\ell+1}$ of the permutation group
$\mathfrak{S}_{\ell+1}$ and $l(s)$ is the length of the permutation (the
number of terms in the minimal decomposition of $s$ into elementary
permutations).
\end{prop}
\noindent{\it Proof}.
1. In the  domain $x_1\geq x_2\ldots\geq x_{\ell+1}$ the
integration domain  ${\cal S}$ can be
deformed so that the integral \eqref{levzeroW}  is given by
a non-trivial sum of residues. To calculate contributions of the
residues let us note that the integrand is symmetric with respect
to $\mathfrak{S}_{tot}=
\mathfrak{S}_2\times \mathfrak{S}_3\times \cdots \times \mathfrak{S}_{\ell+1}$
acting on $\{\gamma_{ij}\}$ via permutations of the second index.
The nontrivial residues corresponding to the poles of the integrand  are at the points
$\gamma_{k,i}=\gamma_{k+1,j}$, $\gamma_{k,i}\neq \gamma_{k,j}$.
Let us first consider the residue contribution at the poles 
$\gamma_{ij}=\gamma_{kj}$
$$
\sum_{n=1}^{\ell+1}(\gamma_{nj}-\gamma_{n-1,j})x_j\longrightarrow
\sum_{n=1}^{\ell+1}\gamma_{\ell+1,j}x_j.
$$
The contribution of the rational function in the integrand  is reduced
after cancellations to
$$
\frac {1}{\prod\limits_{s< p}
  \imath(\gamma_{\ell+1,s}-\gamma_{\ell+1,p})}.
$$
The integrand is symmetric under action of $\mathfrak{S}_{tot}$ and
thus the contributions of other poles
$\gamma_{k,i}=\gamma_{k+1,j}$, $\gamma_{k,i}\neq \gamma_{k,j}$ can be
obtained by averaging  over the permutations. 
Note that the subgroup $\mathfrak{S}_2\times \mathfrak{S}_3\times \cdots
\times \mathfrak{S}_{\ell}\subset \mathfrak{S}_{tot}$  acts trivially on the
residue contribution. Thus the averaging over the subgroup 
cancels the factor $\prod_{n=1}^{\ell}(n!)^{-1}$ and the averaging over
$\mathfrak{S}_{\ell+1}$ gives
$$
\sum_{s\in \mathfrak{S}_{\ell+1}}
\frac{e^{\sum_{k=1}^{\ell+1}\imath\lambda_{s(k)}x_{k}}}
{\prod_{i<j}\imath(\lambda_{s(i)}-\lambda_{s(j)})}=
\sum_{s\in \mathfrak{S}_{\ell+1}}
 (-1)^{l(s)}\frac{e^{\sum_{k=1}^{\ell+1}\imath\lambda_{s(k)}x_{k}}}
{\prod_{i<j}\imath(\lambda_{i}-\lambda_{j})}.
$$
2. Outside  the dominant domain the integral \ref{sumover} is equal
to zero by deformation-of-contour arguments.
$\Box$

It is easy to see that the representation \eqref{sumover} can be
written in the following recursive form
\be\label{sumoverRec}
\Psi^{(0)}_{\gamma_{\ell+1,1},\ldots,\gamma_{\ell+1,\ell+1}}
(x_1,\ldots,x_{\ell+1})\,=\,
\ee
$$
\int_{{\cal{S}}_{\ell,\ell+1}}\prod_{j=1}^{\ell}\frac{d\gamma_{\ell,j}}{2\pi\imath}
\frac{\prod\limits_{s\neq p} (\gamma_{\ell,s}-\gamma_{\ell,p})}
{\imath^{\ell}\ell!\prod\limits_{k=1}^{\ell}\prod\limits_{m=1}^{\ell+1}
(\gamma_{\ell,k}-\gamma_{\ell+1,m})}
\Psi^{(0)}_{\gamma_{\ell,1},\ldots,\gamma_{\ell,\ell}}(x_1,\ldots,x_{\ell}),
$$
where $\!{\cal{S}}_{\ell,\ell+1}\!$ is defined by the
conditions
$\!\max_{j}\{{\rm Im}\,\gamma_{\ell,j}\}<
\min_m\{{\rm Im}\,\gamma_{\ell+1,m}\}\!$.
We call the functions \eqref{levzeroW}  elementary Whittaker functions
 due the following result.

\begin{prop} The elementary $\mathfrak{gl}_{\ell+1}$-Whittaker
  function  \eqref{levzeroW} can  be represented as the limit
\be\label{Wlimit}
\Psi^{(0)}_{\underline{\la}}(\underline{x})=\lim_{\hbar\to
  \infty}(\hbar)^{-\ell(\ell+1)/2}
\Psi_{\underline{\la}}(\hbar\underline{x},\hbar),
\ee
where $\Psi_{\underline{\la}}(\underline{x},\hbar)$
is defined by (\ref{levoneW}).
\end{prop}
\noindent {\it Proof}.
Consider the asymptotic behavior of the
numerator of  (\ref{levoneW}). By definition of the  contour
$\cal S$  in (\ref{levoneW}),  the real parts of  arguments of $\Gamma$-functions
are positive. Then \eqref{levzeroW} is obtained from
\eqref{levoneW} by applying  the following asymptotic form
of the $\Gamma$-function:
$$
\hbar^{-1}\hbar^{\frac{z}{\hbar}}\Gamma\left(\frac{z}{\hbar}\right)\,=\,
\int_{-\infty}^{\infty}dx e^{-z x}\,e^{-\hbar^{-1}e^{-\hbar x}}
\longrightarrow\int_{-\infty}^{\infty}\,dxe^{-zx}\,\Theta(x)= \frac{1}{z}, \qquad   \hbar\to
\infty, \qquad {\rm Re}\, z >0,
$$
where $\Theta(x)$ is the
Heaviside function i.e. $\Theta(x)=1$ for $x\geq 0$ and zero
otherwise. $\Box$

The integral representation \eqref{levzeroW} is an analog
of the Mellin-Barnes integral representation \eqref{levoneW}
and \eqref{sumoverRec} is an analog of  the fundamental recursive
property of the Mellin-Barnes integral representation \cite{KL}.
Taking $\hbar\to \infty$ in \eqref{giv} one obtains
an elementary analog of the Givental integral representation.

\begin{prop} The elementary $\mathfrak{gl}_{\ell+1}$-Whittaker
  function has the following Givental type integral representation:
\be\label{giv0}
\Psi^{(0)}_{\underline{\lambda}}(\underline{x})\,=\,\int\limits_{\IR^{\ell(\ell+1)/2}}
\exp\left(\imath\sum\limits_{k=1}^{\ell+1}
\lambda_k \left(\sum_{i=1}^{k}T_{k,i}-\sum_{i=1}^{k-1}T_{k-1,i}\right)\right)\\
\times
\left(\prod_{k=1}^{\ell}\prod_{i=1}^{k}
\Theta\left(T_{k+1,i}-T_{k,i}\right)
\Theta\left(T_{k,i}-T_{k+1,i+1}\right)\right)\prod_{k=1}^{\ell}\prod_{i=1}^{k}dT_{k,i},
\ee
Here we define $x_i=T_{\ell+1,i},\,\,\,i=1,\ldots,\ell+1$
and we assume $T_{k,i}=0$ for $i>k$. The function $\Theta(x)$ is the
Heaviside function i.e. $\Theta(x)=1$ for $x\geq 0$ and zero
otherwise.
\end{prop}
{\it Proof}.
Let us change the  variables $T_{k,i}\rightarrow\hbar
T_{k,i},\,\,\,k=1,\ldots,\ell+1,\,\,i=1,
\ldots,k$ in (\ref{giv}). Now we take the limit
$\lim_{\hbar\rightarrow\infty}\hbar^{-\ell(\ell+1)/2}\Psi^{}_{\underline{\lambda}}
(\hbar\underline{x}|\hbar)$ of the Givental integral representation
\eqref{giv} using the identity
\be\label{lim}
\lim_{\hbar\to \infty}\,e^{- \hbar^{-1}e^{-\hbar x}}=
\Theta(x).
\ee
This gives us \eqref{giv0}. The second statement
is a direct consequence of \eqref{giv0}. $\Box$

\begin{cor} The elementary $\mathfrak{gl}_{\ell+1}$-Whittaker
  function is given by  
\be\label{giv00}
\Psi^{(0)}_\lambda(\underline{x})\,=\,\int\limits_{\cal{D}}
\exp\{\imath\sum\limits_{k=1}^{\ell+1}
\lambda_k
\left(\sum_{i=1}^{k}T_{k,i}-\sum_{i=1}^{k-1}T_{k-1,i}\right)\}
\prod_{k=1}^{\ell}\prod_{i=1}^{k}dT_{k,i},\\
T_{\ell+1,1}\geq \ldots \geq T_{\ell+1,\ell+1},
\ee
and
\be
\Psi^{(0)}_\lambda(\underline{x})\,=\,0,
\ee
otherwise. Here $\cal{D}$ is a convex  polytope in
$\IR^{\ell(\ell+1)/2}$ defined  by the inequalities
$T_{k,i}\geq T_{k-1,i}\geq T_{k,i+1},\,\,\,\,k=1,
\ldots,\ell+1,\,\,\,i=1,\ldots,\ell+1.$ We assume  $T_{ki}=0$ when
$i>k$.
\end{cor}

\proof Obviously follows from \eqref{giv0}. $\Box$

\begin{de} The elementary $(\ell+1,1)$-Whittaker function specialized
  at $\underline{x}=(x,0,\ldots, 0)$ is defined as the following
  integral:
\be\label{L+1EF}
{}^{(\ell+1,1)}\Psi^{(0)}_{\underline{\la}}(x,0,\ldots,0)\,=\,
\int_{\RR-\imath\epsilon} \,\,\,\frac{d \gamma}{2\pi} \,\,
e^{\imath x \gamma}\,\prod_{j=1}^{\ell+1}
\frac{1}{\imath (\gamma-\lambda_j)},
\ee
where $\underline{\lambda}\,\in\RR^{\ell+1}$
and $\epsilon > 0$.
\end{de}

\begin{prop}.\label{ElmPW} The following limiting expression holds:

(1). $$
{}^{(\ell+1,1)}\Psi^{(0)}_{\underline{\la}}(x,0,\ldots, 0)=\lim_{\hbar\to
  \infty}(\hbar^{-\ell})\,\,
{}^{(\ell+1,1)}\Psi_{\underline{\la}}(\hbar x,0,\ldots ,0|\hbar),
$$
where ${}^{(\ell+1,1)}\Psi_{\underline{\la}}(x,0,\ldots ,0|\hbar)$
is defined by (\ref{identone}).

(2).  The elementary $(\ell+1,1)$-Whittaker function has the following
integral representation:
$$
{}^{(\ell+1,1)}\Psi^{(0)}_{\underline{\la}}(x,0,\ldots,0)
\,=\,\int_{\RR^{\ell}}\,\prod_{j=1}^{\ell}d\tau_j\,\,
  \, e^{\imath\sum_{j=1}^{\ell}\lambda_j
\tau_j+\imath\lambda_{\ell+1}(x-\sum_{j=1}^{\ell}\tau_j)}
\,\,  \Theta(x-\sum_{j=1}^{\ell}\tau_j)\,\,
\prod_{j=1}^{\ell}
\Theta(\tau_j).
$$

\end{prop}

\noindent {\it Proof}. The proof is analogous to the proof of
Proposition \ref{ElmPW}. $\Box$

Let us note that, in contrast with  classical case,
the elementary $(\ell+1,1)$-Whittaker function
can be obtained as a specialization of the
elementary $\mathfrak{gl}_{\ell+1}$-Whittaker function
$$
\Psi^{(0)}_\lambda(x,0,\ldots
,0)={}^{(\ell+1,1)}\Psi^{(0)}_\lambda(x,0,\ldots 0).
$$
Classical
$\mathfrak{gl}_{\ell+1}$-Whittaker function is a common eigenfunction
of a family of mutually commuting differential operators. These
differential operators can be identified with quantum
Hamiltonians of the $\mathfrak{gl}_{\ell+1}$-Toda chain. Similarly
 the elementary  $\mathfrak{gl}_{\ell+1}$-Whittaker functions
are common eigenfunctions of a family of mutually commuting
differential operators defining an elementary analog of
the $\mathfrak{gl}_{\ell+1}$-Toda chains. This  elementary analog of
the $\mathfrak{gl}_{\ell+1}$-Toda chain can be obtained as
the $\hbar\to \infty$ limit of the standard
$\mathfrak{gl}_{\ell+1}$-Toda chain.
We start with the definition of a (well-known) quantum integrable
system which plays the role of the elementary $\mathfrak{gl}_{\ell+1}$-Toda chain
and then explain how this quantum integrable system  arises in the $\hbar\to \infty$
limit from the $\mathfrak{gl}_{\ell+1}$-Toda chain.

Let $(x_1,\ldots ,x_{\ell+1})$ be linear coordinates in
$\IR^{\ell+1}$. The permutation group $\mathfrak{S}_{\ell+1}$ acts in
$\IR^{\ell+1}$ as the group of reflections with respect to the  principal diagonals
$x_i=x_j$. One defines a quantum billiard associated
with the pair $(\IR^{\ell+1},\mathfrak{S}_{\ell+1})$ as a
free quantum particle moving in the closure of the fundamental domain
$\overline{\CD}_{\ell+1}=\{x=(x_1,\ldots,x_{\ell+1})\in \IR^{\ell+1}|x_i\geq x_{i+1}\}$
of $\mathfrak{S}_{\ell+1}$ acting in $\IR^{\ell+1}$. We impose 
Dirichlet boundary conditions on wave functions at the boundary $\pr \overline{\CD}_{\ell+1}$.
The resulting  quantum integrable system  is a special case of the well-known series
of integrable systems on the fundamental domains of  actions of
Weyl groups $W$ on Cartan subalgebras $\mathfrak{h}$
(in our case  the Lie algebra is $A_{\ell}$,
the Cartan subalgebra is $\mathfrak{h}=\IR^{\ell+1}$
and the Weyl group is $W=\mathfrak{S}_{\ell+1}$) (see e.g. \cite{I}).

\begin{prop} The elementary $\mathfrak{gl}_{\ell+1}$-Whittaker function
  \eqref{levzeroW} is a common eigenfunction
of the elementary Toda chain Hamiltonians
\be\label{eigenf1}
P_i(\pr_{x}) \,\Psi_{\underline{\lambda}}(x)=P_i(\lambda)\,
\Psi_{\underline{\lambda}}(x),\qquad P_i(y)\in
\IC[y_1,\ldots ,y_{\ell+1}]^{\mathfrak{S}_{\ell+1}},\qquad
x\in \overline{\CD}_{\ell+1},
\ee
where $\overline{\CD}_{\ell+1}=\{x=(x_1,\ldots,x_{\ell+1})\in
\IR^{\ell+1}|x_i\geq x_{i+1}\}$
 is a compactification of the fundamental domain of
the action of $\mathfrak{S}_{\ell+1}$
in $\IR^{\ell+1}$ and  Dirichlet boundary conditions are imposed
\be\label{BC}
\Psi_{\underline{\lambda}}(x)|_{x_j=x_{j+1}}=0.
\ee
\end{prop}

\proof The representation \eqref{sumover} of the
elementary $\mathfrak{gl}_{\ell+1}$-Whittaker function
 as a sum over the Weyl group $\mathfrak{S}_{\ell+1}$ implies that the
elementary $\mathfrak{gl}_{\ell+1}$-Whittaker function
is a common eigenfunction of the operators $P_i(\pr_{x})$.
Taking into account \eqref{sumoverone}  we infer that the
elementary $\mathfrak{gl}_{\ell+1}$-Whittaker function satisfies the boundary
condition \eqref{BC} and thus solves the eigenvalue problem of the quantum
billiard.  $\Box$

Taking into account \eqref{sumover} we obtain
an expansion of the quantum billiard eigenfunction
\be\label{Dlikef}
\Psi_{\underline{\lambda}}(x)=
\sum_{s\in \mathfrak{S}_{\ell+1}} C^{(0)}(s\cdot \lambda))\,e^{\imath
  \<s\cdot \lambda,x\>},
\ee
where
$$
C^{(0)}(\lambda)=\prod_{\alpha>0}\Gamma^{(0)}(\imath (\lambda,\alpha))=
\prod_{i<j}\Gamma^{(0)}(\imath\lambda_i-\imath\lambda_j),
$$
 the product is over positive
roots of $\mathfrak{gl}_{\ell+1}$ and
 we use the  elementary analog $\Gamma^{(0)}(s)=s^{-1}$ of  classical
Gamma-function. The functions $C^{(0)}(\lambda)$ are elementary
analogs of the Harish-Chandra functions for the Whittaker case.
Indeed, the expansion \eqref{Dlikef}
should  be compared with an  asymptotic expansion
in the region $ x_k >> x_{k+1}$, $k=1,\ldots, \ell$
of  classical class one
$\mathfrak{gl}_{\ell+1}$-Whittaker function (see e.g. \cite{KL})
\be\label{HCexp}
\Psi_\lambda(x)= \sum_{s\in \mathfrak{S}_{\ell+1}} C(s\cdot \lambda)) e^{\imath
 \<s\cdot \lambda,x\>}+O\left({\rm
   max}\left\{e^{x_{k+1}-x_k}\right\}_{k=1}^{\ell}\right),
\ee
where
\be
C(\lambda)=\hbar^{-2\imath \<\lambda,\rho\>/\hbar}\,
\prod_{\alpha>0}\Gamma\left(\frac{(\lambda,\alpha)}{\imath \hbar}\right).
\ee
Thus, the  $\mathfrak{gl}_{\ell+1}$-Whittaker
function \eqref{levoneW} is a common eigenfunction of the
$\mathfrak{gl}_{\ell+1}$-Toda chain and
 the elementary $\mathfrak{gl}_{\ell+1}$-Whittaker function
\eqref{levzeroW} is a common eigenfunction of
 a quantum billiard system  associated with $\mathfrak{gl}_{\ell+1}$.
It is natural expect that the
quantum billiard can be understood as an $\hbar\to \infty$ limit of
$\mathfrak{gl}_{\ell+1}$-Toda chain. Below we demonstrate this relation for
the simplest non-trivial case $\ell=1$.

A ring of quantum Hamiltonians of the $\mathfrak{gl}_2$-Toda chain
is  generated by the two differential operators
$$
\CH_1=-\imath \hbar\frac{\pr}{\pr  x_1}-\imath \hbar\frac{\pr}{\pr  x_2},
$$
$$
\CH_2=-\hbar^2\frac{\pr^2}{\pr x_1^2}-\hbar^2\frac{\pr^2}{\pr x_2^2}+
e^{x_1-x_2}.
$$
Let us make the change of variables $y_i=\hbar^{-1} x_i$ to obtain
$$
\CH_2=-\frac{\pr^2}{\pr y_1^2}-\frac{\pr^2}{\pr y_2^2}+
e^{\hbar(y_2-y_1)}.
$$
Elementary analogs of the quantum Hamiltonians
obtained by taking the limit $\hbar\to \infty$ are given by
\be\label{lzeroHone}
\CH^{(0)}_1=\CH_1=-\imath \frac{\pr}{\pr  y_2}-\imath \frac{\pr}{\pr  y_2},
\ee
\be\label{lzeroHtwo1}
\CH^{(0)}_2=\lim_{\hbar \to \infty} \CH_2=
-\frac{\pr^2}{\pr y_1^2}-\frac{\pr^2}{\pr y_2^2}+
\Theta_{\infty}(y_2-y_1),
\ee
where
$$
\Theta_\infty(y)=0,\quad y>0,\qquad \Theta_\infty(y)=\infty,\quad
y<0.
$$
One can explicitly check  that the elementary Whittaker function
\be\label{Eig2}
\psi_{\lambda_1,\lambda_2}(y_1,y_2)=\frac{e^{\imath(\lambda_1
 y_1+\lambda_2y_2)}-e^{\imath(\lambda_2y_1+\lambda_1y_2)}}{\imath(\lambda_1-\lambda_2)}
\Theta(y_2-y_1),
\ee
is an eigenfunction of \eqref{lzeroHone} and \eqref{lzeroHtwo1}. 
Indeed,  we have
$$
\left(-\frac{\pr^2}{\pr y_1^2}-\frac{\pr^2}{\pr y_2^2}\right)
 \psi_{\lambda_1,\lambda_2}(y_1,y_2)=(\lambda_1^2+\lambda_2^2)
\psi_{\lambda_1,\lambda_2}(y_1,y_2)+e^{\frac{\imath}{2}(\lambda_1+\lambda_2)(y_1+y_2)}
\delta(y_1-y_2).
$$
The function $\Theta_{\infty}(y_2-y_1)\psi_{\lambda_1,\lambda_2}(y_1,y_2)$ is
zero for $y_1\neq y_2$ and is infinite at $y_1=y_2$. The precise character of
the infinity is fixed by the limiting procedure $\hbar\to \infty$ and leads to an
identification
$$
\Theta_{\infty}(y_2-y_1)\,\frac{e^{\imath(\lambda_1
 y_1+\lambda_2y_2)}-e^{\imath(\lambda_2y_1+\lambda_1y_2)}}{\imath(\lambda_1-\lambda_2)}
\Theta(y_2-y_1):= e^{\frac{\imath}{2}(\lambda_1+\lambda_2)(y_1+y_2)}
\delta(y_1-y_2).
$$
Thus, with this definition of the product
$\Theta_{\infty}(y_2-y_1)\psi_{\lambda_1,\lambda_2}(y_1,y_2)$, we obtain that
 the limit $\hbar\to \infty$ of $\mathfrak{gl}_2$-Toda chain is
 indeed given by the quantum $(\IR^2,\mathfrak{S}_2)$-billiard. These
considerations can be straightforwardly generalized to the case of an arbitrary
 rank. We however prefer just to define   elementary
analogs of $\mathfrak{gl}_{\ell+1}$-Toda chain as
quantum $(\IR^{\ell+1},\mathfrak{S}_{\ell+1})$-billiard
to avoid  ill-defined manipulations with $\Theta_{\infty}(x)$.

Finally we define elementary analogs of the local Archimedean
$L$-factors \eqref{lAL} and the Baxter operators \eqref{Baxter}.

\begin{de} Let $V=\IC^{\ell+1}$ be supplied with the standard action
of $U_{\ell+1}$. The elementary  local $L$-factor
associated with $V$ and with an endomorphism
$\Lambda\in {\rm End}(V)$ given by a diagonal matrix
$\Lambda={\rm diag}(\lambda_1,\cdots,\lambda_{\ell+1})$
  is defined as follows:
\be\label{levzeroL}
L_V^{(0)}(s,\underline{\lambda})=\prod_{j=1}^{\ell+1}\frac{1}{s-\lambda_j},
\ee
where $\underline{\lambda}=(\lambda_1, \ldots , \lambda_{\ell+1})$, and $s\in\IC$.
\end{de}

\begin{prop}\label{prop1}
 The following relation between \eqref{levzeroL} and \eqref{lAL}
holds:
$$
L_V^{(0)}(s,\underline{\lambda})=\lim_{\hbar\to
  \infty}(\frac{\imath}{\hbar})^{(\ell+1)}L(s,\underline{\lambda},\hbar).
$$
\end{prop}
\noindent {\it Proof}.
Let ${\rm Re}\, z
>0$, then
$$
\lim_{\hbar\rightarrow\infty}\hbar^{-1}\hbar^{\frac{z}{\hbar}}
\Gamma\left(\frac{z}{\hbar}\right)\,=
\,\int_{-\infty}^{\infty}\,dx \,\,e^{-zx}\,\,\Theta(x)\,=\,\frac{1}{z},
$$
and thus
$$
\lim_{\hbar\to \infty}(\frac{\imath}{\hbar})^{(\ell+1)}L(s,\underline{\lambda},\hbar)=
\imath^{\ell+1}\,\int_{-\infty}^{\infty}\prod_{j=1}^{\ell}dt_j\,\prod_{j=1}^{\ell}\,
e^{\imath t_j(s-\lambda_j)}\,\Theta(t_j)=
\prod_{j=1}^{\ell+1}\frac{1}{s-\lambda_j},
$$
provided    ${\rm Im}(s)<0$. $\Box$

The  eigenvalues of the Baxter integral
operators \eqref{Baxter}  acting on the $\mathfrak{gl}_{\ell+1}$-Whittaker functions
\eqref{levoneW} are given by the local Archimedean $L$-factors.
Similar relations hold for their elementary analogs.
Consider the one-dimensional family ${}^{\mathfrak{gl}_{\ell+1}}\mathcal{Q}^{(0)}(s)$
of integral operators acting in an 
appropriate space of functions of $\ell+1$ variables
and having the integral kernel
\be
{}^{\mathfrak{gl}_{\ell+1}}\mathcal{Q}^{(0)}(\underline{x},\,\underline{y}|\,s)=
e^{\,\imath s \sum_{i=1}^{\ell+1}(x_i-y_i)}\,
\Theta(y_{\ell+1}-x_{\ell+1})\,\,
\prod_{i=1}^{\ell}\Theta(y_i-x_i)\,\Theta(x_i-y_{i+1}),
\ee
where we assume $\underline{x}:=(x_1,\cdots , x_{\ell+1})$ and
$\underline{y}:=(y_1,\ldots , y_{\ell+1})$.

\begin{prop} The following identity holds:
\be
{}^{\mathfrak{gl}_{\ell+1}}\mathcal{Q}^{(0)}(s)=\lim_{\hbar\to
  \infty} \,\,\,
\mathcal{Q}^{\mathfrak{gl}_{\ell+1}}(s).
\ee
\end{prop}

\noindent {\it Proof}. By changing variables
$\underline{x}\to\hbar\underline{x},\,\,\,\underline{y}\to \hbar \underline{y}$
the proof reduces to straightforward application of the identity \eqref{lim}. $\Box$

\begin{cor} The elementary $\mathfrak{gl}_{\ell+1}$-Whittaker function
  has the following eigenfunction property:
\bqa\label{eigenprop0}
\imath^{(\ell+1)}\,\int_{\RR^{\ell+1}}\,\prod_{i=1}^{\ell+1}\,dy_{i}\,\,
{}^{\mathfrak{gl}_{\ell+1}}\mathcal{Q}^{(0)}(\underline{x},
\,\underline{y}|\, s)\,\,\,\,
{}^{\mathfrak{gl}_{\ell+1}}\Psi^{(0)}_{\underline{\lambda}}(\underline{y})\,=\,
L_V^{(0)}(s,\underline{\lambda},\hbar)\,\,
{}^{\mathfrak{gl}_{\ell+1}}\Psi^{(0)}_{\underline{\lambda}}(\underline{x}),\eqa
where $\underline{x}=(x_1,\ldots, x_{\ell+1})$,
$\underline{y}=(y_1,\ldots, y_{\ell+1})$,
$\underline{\la}=(\la_1,\ldots, \la_{\ell+1})$
and the eigenvalue is equal to  the local Archimedean $L$-factor
\be\label{lAL0}
L^{(0)}(s,\underline{\lambda},\hbar)=\prod_{j=1}^{\ell+1}\,\frac{1}{
s- \lambda_{j}}.
\ee
\end{cor}

\section{A limit of $q$-deformed Whittaker functions}

In \cite{GLO3}, \cite{GLO4}, \cite{GLO5}  explicit expressions for
$q$-deformation of the $\mathfrak{gl}_{\ell+1}$-Whittaker
functions were proposed. The $q$-deformed Whittaker functions
are common eigenfunctions  of $q$-deformed Toda chains (also known as
relativistic Toda chains) \cite{R}, \cite{Et} and the standard Whittaker
functions arise in  the limit $q\to 1$. In the previous Section we
defined the elementary $\mathfrak{gl}_{\ell+1}$-Whittaker functions as
a limit of classical
$\mathfrak{gl}_{\ell+1}$-Whittaker functions.
In this Section we demonstrate that the
elementary $\mathfrak{gl}_{\ell+1}$-Whittaker function can
be obtained directly from the $q$-deformed $\mathfrak{gl}_{\ell+1}$-Whittaker function
$\Psi^{q}_{\underline{z}}(\underline{n})$
by specializing at $q=0$ and taking a limit with respect to spectral
variables $\underline{z}=(z_1,\ldots ,z_{\ell+1})$.
Similar relations hold between  the limits $\hbar\to
\infty$ of  classical  local Archimedean $L$-factors/Baxter
operators on the one hand, and the $q=0$ specialization  of   $q$-deformed local
$L$-factors (introduced in \cite{GLO3}, \cite{GLO4}, \cite{GLO5})/$q$-deformed
Baxter operators on the other hand. 

Specialization $q=0$ of the $q$-deformed class one
$\mathfrak{gl}_{\ell+1}$-Whittaker function was already considered in
\cite{GLO3}. Under this specialization the $q$-deformed
$\mathfrak{gl}_{\ell+1}$-Whittaker function
$\Psi^{q}_{\underline{z}}(\underline{n})$
is given by a character of a finite-dimensional irreducible
representation of $\mathfrak{gl}_{\ell+1}$ corresponding to the partition
$n_1\geq\ldots\geq n_{\ell+1}$
\be\label{pW}
\Psi^{q=0}_{z_1,\cdots,
  z_{\ell+1}}(n_1,\cdots,n_{\ell+1})=
\Tr_{V_{n_1,\cdots, n_{\ell+1}}}z_1^{H_1}\cdots
z_{\ell+1}^{H_{\ell+1}},
\ee
and is equal to zero for $(n_1,\cdots ,n_{\ell+1})$ outside
the principal domain $n_1\geq\ldots\geq n_{\ell+1}$.
Using the Weyl character formula (see e.g. \cite{Zh}) we have
\be\label{Weylf}
\Psi^{q=0}_{z_1,\cdots,
  z_{\ell+1}}(n_1,\cdots,n_{\ell+1})=\frac{1}{\prod_{i<j}(z_i-z_j)}\,\sum_{s\in
  \mathfrak{S}_{\ell+1}}\,(-1)^{l(s)}\,\prod_{j=1}^{\ell+1}z_j^{n_{s(j)}+\ell+1-s(j)},
\ee
where $\mathfrak{S}_{\ell+1}$ is the permutation group identified with  the
Weyl group of $\mathfrak{gl}_{\ell+1}$ and for an element $s\in
\mathfrak{S}_{\ell+1}$, $l(s)$ is  the length of the
minimal product decomposition of $s$.  There exists  another
 representation  for the characters of irreducible finite-dimensional representations
of $\mathfrak{gl}_{\ell+1}$ based on  Gelfand-Zetlin bases
in finite-dimensional irreducible representations \cite{GZ} (see also \cite{Zh}).
Let $\CP^{\ell+1}$  be a  set of Gelfand-Zetlin patterns,
that is a set of collections $\underline{p}=\{p_{i,j}\}$, $i=1,\ldots,
\ell+1$, $j=1,\ldots, i$ of
integers satisfying the conditions
 $p_{i+1,j}\geq p_{i,j}\geq p_{i+1,j+1}$.
An irreducible finite-dimensional representation can be realized in
a vector space with the basis $v_{\underline{p}}$ enumerated by the
Gelfand-Zetlin patterns $\underline{p}$ with  fixed $p_{\ell+1,i}$, $i=1,\ldots,
\ell+1$. Action of the Cartan generators on $v_{\underline{p}}$ is
then given by 
\be
z_1^{H_1}\,z_2^{H_2} \cdots
z_{\ell+1}^{H_{\ell+1}}\,v_{\underline{p}}=
z_1^{s_1}\,z_2^{s_2-s_1}\,\cdots
z_{\ell+1}^{s_{\ell+1}-s_{\ell}}\,v_{\underline{p}},\qquad
s_k=\sum_{i=1}^k\,p_{ki}. \ee  Set
 $(n_1,\ldots,n_{\ell+1}):=(p_{\ell+1,1},\ldots,p_{\ell+1,\ell+1})$.
Then for the  character of an irreducible representation of $GL_{\ell+1}$
corresponding to a partition $(n_1\geq \ldots  \geq n_{\ell+1})$ we have
\be\label{qzero}
\Psi^{q=0}_{\underline{z}}(\underline{p}_{\ell+1})
\,=\,\sum_{p_{k,i}\in{\cal P}^{\ell+1}}\,\, \prod_{k=1}^{\ell+1}
z_k^{(\sum_{i=1}^k p_{k,i}-\sum_{i=1}^{k-1} p_{k-1,i})}\, .\ee
This representation can be written in the recursive form
 \be
\Psi^{q=0}_{z_1,\cdots
  z_{\ell+1}}(n_1,\cdots,n_{\ell+1})=
\sum_{p_{\ell,i}\in{\cal P}_{\ell+1,\ell}}
   \,\,\, z_{\ell+1}^{ \sum_{i=1}^{\ell+1}
  p_{\ell+1,i}-\sum_{i=1}^{\ell}  p_{\ell,i}}
\Psi^{q=0}_{z_1,\cdots,
  z_{\ell}}(p_{\ell,1},\ldots ,p_{\ell,\ell}),
\ee
 where the sum runs over the  set
$\CP_{\ell+1,\ell}$ of  $\underline{p}_{\ell}=(p_{\ell,1},\ldots
,p_{\ell,\ell})$ satisfying the   conditions
$p_{\ell+1,i}\geq p_{\ell,i}\geq p_{\ell+1,i+1}$.
One can also write down irreducible characters in the following
recursive integral form
(see e.g. \cite{GLO5}):
\be\label{reccharGZ}
\Psi^{q=0}_{z_1,\cdots,
  z_{\ell+1}}(n_1,\cdots,n_{\ell+1})
\,=\,
\oint_{y_1=\infty}\cdots\oint_{y_{\ell+1}=\infty}\,
\prod_{i=1}^{\ell}\,\frac{\imath\,dy_i}{2\pi y_i} \,
C_{\ell+1,\ell}(z,y^{-1})\times
\ee
$$
\times\Psi^{q=0}_{y_1,\cdots
  y_{\ell}}(n_1,\cdots,n_{\ell})
\,\Delta(y),
$$
$$
\Psi^{q=0}_{z_1,\cdots
  z_{\ell+1}}(n_1+k,\cdots,n_{\ell+1}+k)
\,=\,
\Big(\prod_{j=1}^{\ell+1}x_j^k\Big)
\Psi^{q=0}_{z_1,\cdots
  z_{\ell+1}}(n_1,\cdots,n_{\ell+1}),$$
where
$$
C_{\ell+1,m+1}(x,y)=
\prod_{i=1}^{\ell+1}\prod_{j=1}^{m+1}\,\frac{1}{1-x_i y_j},\qquad
\Delta(y)\,=\,\prod_{i\neq j}\left(1-y_iy_j^{-1}\right).
$$
Let $\chi_{r}(\underline{z})$ be characters of the fundamental representations
$V_{\omega_r}=\bigwedge^r\mathbb{C}^{\ell+1}$ of $\mathfrak{gl}_{\ell+1}$
$$
\chi_{r}(\underline{z})\,=\,
\sum_{I_r}\,z_{i_1}\cdots z_{i_r}\,,\hspace{1.5cm}r=1,\ldots,\ell+1,
$$
and $I_r\,=\,(i_1<i_2<\ldots<i_r)\subseteq\{1,2,\ldots,\ell+1\}$.
The functions $\Psi^{q=0}_{z_1,\cdots z_{\ell+1}}(n_1,\cdots,n_{\ell+1})$, equal
to the characters of irreducible
finite-dimensional  representations of $\mathfrak{gl}_{\ell+1}$
in the principal domain  $n_{1}\geq\ldots\geq n_{\ell+1}$,
 satisfy a system of difference equations expressing the Piery
branching rules
\be\label{Htozero}
\chi_{r}(\underline{z})\,
\Psi^{q=0}_{\underline{z}}(\underline{p}_{\ell+1})
\,=\,\sum_{I_r}\,\,
\Psi^{q=0}_{\underline{z}}(\underline{p}_{\ell+1}+\underline{\nu}_{I_r}),\qquad
r=1,\ldots,\ell+1,
\ee
where $I_r\,=\,(i_1<i_2<\ldots<i_r)\subseteq\{1,2,\ldots,\ell+1\}$ and
$\nu_j=1$ or $\nu_j=0$
 depending on whether or not $j$ is in the set $I_r$. We also omit the terms
 in the right hand side of \eqref{Htozero} corresponding to
non-dominant weights.
The characters are uniquely defined as solutions of the equations
\eqref{Htozero} with the condition
$\Psi^{q=0}_{\underline{z}}(\underline{p}_{\ell+1})=0$ for
$\underline{p}_{\ell+1}$ outside the principal domain $p_{\ell+1,1}\geq\ldots\geq
p_{\ell+1,\ell+1}$.

\begin{prop} Let $z_i=t^{-\lambda_i}$, $\lambda_i\in \IC$, $t\in \IR_+$
and for $x_j\in \IR$ let $n_i(t,x_j)$ be the integer parts  of $x_j/\log t$.
In the limit $t\to 1$ the function \eqref{pW} reduces to
  the elementary $\mathfrak{gl}_{\ell+1}$-Whittaker function \eqref{levzeroW}
 $$
\Psi^{(0)}_{\lambda_1,\ldots, \lambda_{\ell+1}}(x_1,\ldots ,x_{\ell+1})=\lim_{t\to 1}
\,\,(\log t)^{\ell(\ell+1)/2}\,\,
\Psi^{q=0}_{t^{-\lambda_1},\cdots   ,
  t^{-\lambda_{\ell+1}}}(n_1(t,x_1),\cdots,n_{\ell+1}(t,x_{\ell+1})).
$$
\end{prop}

\noindent {\it Proof}. Taking the limit of the recursive relations \eqref{reccharGZ}
we obtain the recursive relations \eqref{sumoverRec}.
$\Box$

Similar relations hold between $q$-deformed local $L$-factors
(introduced in \cite{GLO3}, \cite{GLO4}, \cite{GLO5}) and elementary local $L$-factors
\eqref{lAL0}. Recall that the $q$-deformed local $L$-factor is given by
\be
L^q(u|z_1,\cdots,
z_{\ell+1})=\prod_{j=1}^{\ell+1}\prod_{n=0}^{\infty}\frac{1}{1-u^{-1}z_jq^n},\,\,\,q<1.
\ee
Taking $q=0$ we have
\be\label{q=0LLF}
L^{q=0}(u|z_1,\cdots,
z_{\ell+1})=\prod_{j=1}^{\ell+1}\frac{1}{1-u^{-1}z_j}.
\ee
Now let $z_i=t^{-\lambda_i}$, $u=t^{-s}$. 
Taking the  limit $t\to 1$ we recover the elementary
local $L$-factor
\be\label{LLimq}
L^{(0)}(u|\underline{\lambda})=\lim_{t\to 1}\,(\log
t)^{\ell+1}\,\prod_{j=1}^{\ell+1}\frac{1}{1-t^{s-\lambda_j}}=
\prod_{j=1}^{\ell+1}\frac{1}{s-\lambda_j}.
\ee
Local Archimedean $L$-factors are eigenvalues of the Baxter integral
operators acting on the Whittaker functions \cite{GLO2}. 
Similar relations hold
for their $q$-deformations. Next we give explicit formulas for
$q=0$ specialization of the $q$-deformed Baxter operators.

\begin{prop} The $q=0$ specialization \eqref{pW} of the $q$-deformed
$\mathfrak{gl}_{\ell+1}$-Whittaker function is a common eigenfunction
of the one-parameter family of operators
\be\label{q=0Baxter}
\CQ^{q=0}(u) \cdot f(n_1,\ldots ,n_{\ell+1})\,=\,
\ee
$$ \sum_{m_1,\cdots,
  m_{\ell+1}\in \IZ} \CQ^{q=0}(u|n_1,\ldots n_{\ell+1};m_1,\ldots ,m_{\ell+1})\,\,
f(m_1,\ldots ,m_{\ell+1}),
$$
where the kernel is given by
\be\label{Qker}
\CQ^{q=0}(u|n_1,\ldots n_{\ell+1};m_1,\ldots ,m_{\ell+1})=\ee
$$
u^{\,\sum_{i=1}^{\ell+1}(n_i-m_i)}\,
\Theta(m_{\ell+1}-n_{\ell+1})\,\,
\prod_{i=1}^{\ell}\Theta(m_i-n_i)\,\Theta(n_i-m_{i+1}).
$$
The corresponding  eigenvalues  are  given by  local $L$-factors
\eqref{q=0LLF}
\be\label{q=0B}
\CQ^{q=0}(u) \cdot \Psi^{q=0}_{z_1,\cdots,
  z_{\ell+1}}(n_1,\cdots,n_{\ell+1})\,=\,
\ee
$$
\ L^{q=0}(u|z_1,\ldots,z_{\ell+1} )\,
\Psi^{q=0}_{z_1,\cdots,
  z_{\ell+1}}(n_1,\cdots,n_{\ell+1}).
$$
\end{prop}

\noindent {\it Proof}.  Recall that $\Psi^{q=0}_{z_1,\cdots,
 z_{\ell+1}}(n_1,\cdots,n_{\ell+1})$ can be interpreted as
characters of  a finite-dimensional
 irreducible representation of $\mathfrak{gl}_{\ell+1}$. The local
 $L$-factor \eqref{q=0LLF} can be  considered as a character of
 an infinite-dimensional  representation of $\mathfrak{gl}_{\ell+1}$
$$
L^{q=0}(u|z_1,\cdots,
z_{\ell+1})=\prod_{j=1}^{\ell+1}\frac{1}{1-u^{-1}z_j}=\sum_{(k_1,\ldots,k_{\ell+1})
\in \IZ_+^{\ell+1}}u^{-(k_1+\ldots +k_{\ell+1})}\,z_1^{k_1}\cdots
z_{\ell+1}^{k_{\ell+1}}=
$$
$$
=\Tr_{V_0}(z_1/u)^{H_1}\cdots
(z_{\ell+1}/u)^{H_{\ell+1}},
$$
where $V_0=\IC[\xi_1,\cdots, \xi_{\ell+1}]$ and
$H_j=\xi_j\frac{\pr}{\pr \xi_j}$.
Now \eqref{q=0B} is derived by applying  the standard
Littlewood-Richardson rules for  decomposition of a tensor product
 of  representations (see e.g. \cite{Fu}, \cite{Zh})
 \be
\chi_{r,0,\ldots,0}(\underline{z})\chi_{n_1,\ldots,n_{\ell+1}}(\underline{z})
=\sum_{I_r} \chi_{\underline{n}+\underline{\nu}_{I_r}}(\underline{z}),
\ee
where $I_r\,=\,(i_1<i_2<\ldots<i_r)\subseteq\{1,2,\ldots,\ell+1\}$ and
$\nu_j=1$ or $\nu_j=0$ depending on  whether or not $j$ is in $I_r$.
 We also discard all terms on the right hand side for which
$\underline{n}+\underline{\nu}_{I_r}$ is not in the principal domain.
Now the statement of the Proposition follows from the  decomposition
$$
L^{q=0}(t|z_1,\ldots,
z_{\ell+1})=\sum_{k=0}^{\infty}t^{-k}\sum_{k_1+\ldots+k_{\ell+1}
=k}z_1^{k_1}\cdots z_{\ell+1}^{k_{\ell+1}}=
\sum_{k=0}^{\infty}t^{-k}\chi_{k,0,\ldots}(z_1,\ldots,z_{\ell+1}).
$$
 $\Box$

Let us consider in detail the eigenfunction equation \eqref{q=0B}
for $\ell=0,1$.
 We start with $\ell=0$.
 Irreducible representations of $U_1$ are one-dimensional
$$
e^{\imath \theta}:\,\,\, V_n\longrightarrow e^{\imath n\theta}\,V_n.
$$
Let $H$ be a generator of ${\rm Lie}(U_1)$.  The $q=0$
specialization of the $q$-deformed
$\mathfrak{gl}_1$-Whittaker function is given by the character
$$
\Psi^{q=0}_t(n)=\Tr_{V_n}\,t^H=t^n,\qquad t\in \IC^*,\qquad n\in \IZ.
$$
The local Archimedean $L$-factor has the following trace representation:
$$
L^{q=0}(u|t)=\frac{1}{1-u^{-1}t}=\sum_{n=0}^{\infty}(t/u)^n=
\Tr_{\IC[z]}\,(t/u)^H=\sum_{n=0}^{\infty} u^{-n}\Psi^{q=0}_t(n),
$$
where $H$ acts in $\IC[z]$ as $H=z\pr_z$.
We have the following  decomposition of the product of 
characters of $V_n$ and $\CH$:
$$
\Tr_{V_n}t^H\times
\Tr_{\IC[z]}\,(u^{-1}t)^H=\sum_{m=0}^{\infty}\Tr_{V_{m+n}}t^Hu^{n-H}=
\sum_{m\in \IZ} \Theta(m-n)\,u^{n-m}\,\Tr_{V_m}t^{H}.
$$
This can be rewritten as the action of an  integral
operator on the Whittaker function $\Psi^{q=0}_t(n)=t^n$
$$
\CQ(u)\cdot \Psi^{q=0}_t(n)=L(u|t)\,\Psi^{q=0}_t(n),
$$
where $Q(u)$ acts on functions on $\IZ$ as follows:
$$
\CQ(u)\cdot f(n)=\sum_{m\in \IZ} \Theta(m-n) \,u^{n-m}\,f(m).
$$

Now  consider the case $\ell=1$. Let $V_{n_1,n_2}$ be the
finite-dimensional irreducible representation of $\mathfrak{gl}_2$
corresponding to a partition $(n_1,n_2)$, $n_1\geq n_2$.
 Let $H_i$, $i=1,2$ be  generators of the diagonal Cartan subalgebra.
The $q=0$ specialization of the $q$-deformed $\mathfrak{gl}_2$-Whittaker function is
expressed via  characters of irreducible finite-dimensional
representations as
$$
\Psi^{q=0}_{t_1,t_2}(n_1,n_2)=\Tr_{V_{n_1,n_2}}t_1^{H_1}t_2^{H_2}=
(t_1t_2)^{n_2}\,\left(\sum_{m_1+m_2=n_1-n_2}\,t_1^{m_1}t_2^{m_2}
\right)
=\frac{t_1^{n_1+1}t_2^{n_2}-t_1^{n_2}t_2^{n_1+1}}{t_1-t_2}.
$$
The local $L$-factor has the following representation:
$$
L^{q=0}(u|t_1,t_2)=\frac{1}{(1-u^{-1}t_1)(1-u^{-1}t_2)}=\Tr_{\IC[z_1,z_2]}
(u^{-1}t_1)^{H_1}(u^{-1}t_2)^{H_2}=\sum_{r=0}^{\infty}u^{-r}
\Tr_{V_{r,0}}t_1^{H_1}t_2^{H_2},
$$
where $H_j$, $j=1,2$  act in $\IC[z_1,z_2]$ as $H_j=z_j\pr_{z_j}$.
The following identity (a particular instance of the
Richardson-Littlewood rule) holds:
$$
\chi_{(r_1,r_2)}\cdot \chi_{(k,0)}=\sum_{I_{r,k}\in \IZ^2} \chi_{p_1,p_2},
$$
where $\chi_{n_1,n_2}(t_1,t_2):=\Tr_{V_{n_1,n_2}}t_1^{H_1}t_2^{H_2}$
and the sum goes over the subset $I_{r,k}\in \IZ^2$  for which
$$
p_2+p_1=r_1+r_2+k,\qquad p_1\geq r_1\,\qquad p_2\geq r_2,\qquad
r_1\geq p_2.
$$
This can be rewritten as follows:
\be\label{l=2br}
\chi_{n_1,n_2}(t_1,t_2)\,\chi_{S^r\IC^2}(u^{-1}t_1,u^{-1}t_2)=
\ee
$$
=\sum_{m_1,m_2\in\IZ^2}
\Theta(m_1-n_1)\Theta(n_1-m_2)\Theta(m_2-n_2) u^{n_1+n_2-m_1-m_2}\chi_{m_1,m_2}(z_1,z_2).
$$
Thus we recover a special case of \eqref{q=0B} with
the kernel of the Baxter operator given by
$$
\CQ(u|n_1,n_2;m_1,m_2)=\Theta(m_1-n_1)\Theta(n_1-m_2)
\Theta(m_2-n_2)\,u^{n_1+n_2-m_1-m_2}.
$$
Let us finally note that the elementary Baxter operator \eqref{prop1}
can be obtained as a limit of its $q=0$ counterpart.

\begin{prop} Let $u=t^{\lambda}$ and
let $n_i(t,x_j)$, $m_j(t,y_j)$ be the corresponding  integer parts of $x_j/\log t$,
 $y_j/\log t$ respectively. Then the following identity holds:
\be
\mathcal{Q}^{(0)}(x_1,\cdots x_{\ell+1};y_1,\ldots y_{\ell+1}|\la)=
\ee
$$
=\lim_{t\to 1}
\CQ^{q=0}(u|n_1(t,x_1),\ldots n_{\ell+1}(t,x_{\ell+1});
m_1(t,y_1),\ldots ,m_{\ell+1}(t,y_{\ell+1})).
$$
\end{prop}

\noindent {\it Proof}. The proof is straightforward. $\Box$

\section{Elementary special functions as symplectic volumes}

In this Section we represent the elementary special functions defined in
the previous Section as equivariant  volumes of
finite-dimensional symplectic spaces. As  was already noticed at the
end of Section 2 this representation is  not surprising,
 given the existence of a similar representation of  classical Whittaker
functions and local $L$-factors  as equivariant volumes of
infinite-dimensional symplectic spaces of holomorphic maps of a
two-dimensional disk into finite-dimensional symplectic spaces
\cite{GLO6}, \cite{GLO7}, \cite{GLO8}. In the
limit when the $S^1$-equivariance parameter $\hbar$ corresponding to  disk
rotations goes to infinity, the corresponding equivariant volume of the
space of holomorphic maps of the disk into $X$  tends to the equivariant volume of
$X$ given by an  elementary analog of the corresponding special fucntion.
In this Section we demonstrate directly that $U_{\ell+1}$-equivariant
volumes of flag spaces $\CB_{\ell+1}=GL_{\ell+1}(\IC)/B$ are equal to
 elementary $\mathfrak{gl}_{\ell+1}$-Whittaker  functions \eqref{levoneW}. The same relation also
holds for $(\ell+1,1)$-Whittaker functions and local $L$-factors. Thus,
for elementary analogs, we prove the relation between symplectic
volumes and special functions for a more general case  than was done
in \cite{GLO8}. In this Section we also provide a reformulation of the
eigenfunction property \eqref{eigenprop} of the elementary Whittaker
function in terms of symplectic geometry.

Let $\CB_{\ell+1}=GL_{\ell+1}(\IC)/B$ be a flag
space of $GL_{\ell+1}$. It can be identified with the factor
$U_{\ell+1}/H$ of the unitary group $U_{\ell+1}$ over the Cartan torus
$U_1^{\ell+1}$ of $U_{\ell+1}$. The space $U_{\ell+1}/H$ on the other
hand is obviously a coadjoint orbit $\CO_{u_0}$ of a regular element  $u_0\in
\mathfrak{u}_{\ell+1}^*$. This interpretation provides a family of
 Kirillov-Kostant symplectic structures on the $\CB_{\ell+1}$
 parameterized by positive
cone in the dual to the Cartan subalgebra $\mathfrak{h}\in
\mathfrak{u}_{\ell+1}$ (see e.g. \cite{K}). In the following we identify both
$\mathfrak{u}_{\ell+1}^*$ with $\mathfrak{u}_{\ell+1}$
 as well as $\mathfrak{h}^*$ with $\mathfrak{h}$  via the  Killing
quadratic form on $\mathfrak{u}_{\ell+1}$.

Let $u_0=\imath\, {\rm diag}(x_1,\ldots ,x_{\ell+1})\in \mathfrak{u}_{\ell+1}^*$
 with $x_1>x_2>\ldots >x_{\ell+1}$.
The Kirillov-Kostant symplectic form $\omega$ on an open part
of the  coadjoint orbit $\CO_{u_0}$ can be written explicitly as follows.
Consider the closed two-form on $U_{\ell+1}$
$$
\omega^{(0)}_{u_0}=\delta \Tr u_0\,\delta g\,g^{-1}, \qquad g\in
U_{\ell+1},
$$
where trace is taken in the standard representation
$\mathfrak{u}_{\ell+1}\to {\rm End}(\IC^{\ell+1})$.
This  two-form is a lift of the  Kirillov-Kostant closed non-degenerate
two-form $\omega_{u_0}$ on $\CO_{u_0}$ along a projection
$U_{\ell+1}\to \CO_{u_0}$, such that $g\to u=g^{-1}u_0g$.
The action of the group $U_{\ell+1}$ on $(\CO_{u_0},\omega_{u_0})$
is Hamiltonian  and the
corresponding momentum map is given by
\be\label{mmap}
H(g,u_0):=u(g,u_0)=g^{-1}u_0g.
\ee

\begin{prop} The following representation for the elementary
 $\mathfrak{gl}_{\ell+1}$-Whittaker function \eqref{levzeroW}  holds:
\be\label{eqvol}
\Psi^{(0)}_{\lambda_1,\cdots,\lambda_{\ell+1}}
(x_1,\cdots,x_{\ell+1})=\int_{\CB_{\ell+1}}\,
e^{\omega_{u_0}+
\sum_{j=1}^{\ell+1}\lambda_jH_{jj}},\qquad
u_0=\imath \,{\rm diag}(x_1,\ldots ,x_{\ell+1}),
\ee
where $H_{jj}$ are diagonal components of $H(u_0,g)$.
\end{prop}

\noindent {\it Proof}. The integral \eqref{eqvol} can be
calculated using the Harish-Chandra formula for orbit integrals
\be\label{orbint}
I(u_0,\Lambda)=
\int_{g\in U_{\ell+1}/H}\,\,e^{\Tr(g^{-1}u_0 g \Lambda)} \,\,
{\rm vol}_{U_{\ell+1}/H}(g)=
\Delta^{-1}(x)\,\Delta^{-1}(\imath\lambda)\,\,\det\|e^{\imath x_i\l_j}\|=
\ee
$$=
\Delta^{-1}(x)\,\Delta^{-1}(\imath\lambda)\,\,
\sum_{s\in \mathfrak{S}_{\ell+1}}(-1)^{l(s)}
e^{\sum_j \imath x_{s(j)}\l_j},
$$
where $\Lambda={\rm diag}(\lambda_1, \ldots ,\lambda_{\ell+1})$  with
$\lambda_1> \cdots > \lambda_{\ell+1}$
 and the measure ${\rm vol}_{U_{\ell+1}/H}$ is the canonical
volume form on the factor $U_{\ell+1}/H$ induced by the Killing form on
$\mathfrak{u}_{\ell+1}$. Here we also denote $\Delta(z)=\prod_{1\leq
  i<j\leq \ell+1}(z_i-z_j)$.
 Thus,  taking into account \eqref{mmap}, to
derive \eqref{eqvol} from  \eqref{orbint}
 one should compare
${\rm vol}_{U_{\ell+1}/H}$ with the Liouville measure
$\omega_{u_{0}}^{d}/d!$ where $d=\frac{1}{2}\dim (U_{\ell+1}/H)$.  Both
${\rm vol}_{U_{\ell+1}/H}$ and $\omega_{u_{0}}^{\dim
  (U_{\ell+1}/H)/2}$ are $U_{\ell+1}$-invariant top-dimensional  forms
on $\CB$  and thus it is enough
to compare these measures near the projection of the unit element $e\in
U_{\ell+1}$. Simple calculation  gives
$$
\omega_{u_{0}}^{d}=d!\,\, \Delta(x)\,{\rm vol}_{U_{\ell+1}/H},
$$
Thus we obtain
$$
\int_{\CB}\,
e^{\omega_{u_0}+
\sum_{j=1}^{\ell+1}\lambda_jH_{jj}(u_0)}=\Delta^{-1}(\imath\lambda)\,\,
\sum_{s\in \mathfrak{S}_{\ell+1}}(-1)^{l(s)}
e^{\sum_j \imath x_{s(j)}\l_j}.
$$
Taking into account \eqref{sumover} we obtain \eqref{eqvol}. $\Box$

The identity \eqref{eqvol} can  also be proved explicitly using a
Gelfand-Zetlin type parametrization on an open part of the regular
coadjoint orbit $\CO_{u_0}$ of $U_{\ell+1}$
\cite{AFS}. Recall that an open part of the regular coadjoint
orbit $\CO_{u_0}$
allows a parametrization in  Darboux coordinates  $\{T_{ij},\theta_{ij}\}$, $1\leq
j\leq i<\ell+1$ so that the symplectic form   $\omega_{u_0}$
is given by
$$
\omega_{u_0}=\sum_{i\geq j} \delta T_{ij}\wedge \delta
\theta_{ij}.
$$
Here $\theta_{ij}$ are periodic coordinates $\theta_{ij}\sim
\theta_{ij}+1$ and $T_{ij}\in \IR$,
$1\leq j\leq i<\ell+1$  satisfy the Gelfand-Zetlin conditions
\be\label{GZpolytop}
T_{i+1,j}\geq T_{i,j}\geq T_{i+1,j+1},\qquad 1\leq j\leq
i\leq \ell+1,
\ee
and  we identify $T_{\ell+1,j}:=x_j$, $j=1,\ldots ,\ell+1$. Thus the image of an open
part of $\CO_{u_0}$ under the projection along an $\ell(\ell+1)/2$-dimensional
torus parameterized by $\theta_{ij}\in \IR$, $1\leq j\leq i<\ell+1$
is a convex Gelfand-Zetlin polytope ${\mathcal{P}}_{\ell+1}$
in $\IR^{\ell(\ell+1)/2}$ defined
by the conditions \eqref{GZpolytop}. In  coordinates $(T,\theta)$
the components of the momentum map  of the action of the diagonal Lie subalgebra
$\mathfrak{u}_1^{\ell+1}\subset \mathfrak{u}_{\ell+1}$ are given by
$$
H_{jj}=\sum_{i=1}^{j}\imath T_{ji}-\sum_{i=1}^{j-1}\imath T_{j-1,i}.
$$
Thus the integral,  after integration over $\theta_{ij}$,  can be written in the
following form:
\be\label{GZpoly}
\int_{\CB_{\ell+1}}\,
e^{\omega_{u_0}+
\sum_{j=1}^{\ell+1}\lambda_jH_{jj}}=\int_{{\cal P}_{\ell+1}}
\,e^{\sum_{ij}\imath\lambda_j(T_{ij}-T_{i-1,j})}\,\prod_{1\leq j\leq j\leq
 \ell} d T_{ij},
\ee
where $T_{ij}:=0$ for $i<j$.  This integral
representation coincides with the Givental type integral
representation \eqref{giv00} and thus we again recover the identity
\eqref{eqvol}.

\begin{rem} The appearance of the Gelfand-Zetlin polytope
${\cal P}_{\ell+1}$ in the Givental type integral representation
\eqref{giv00} is related   with a deep duality relation between
the Gelfand-Zetlin and the Givental realizations of representations of
$\CU\mathfrak{gl}_{\ell+1}$ discussed in \cite{GLO2}.
\end{rem}

The relation \eqref{eqvol} can  also be understood as a
limit of the identification of the $q=0$
$\mathfrak{gl}_{\ell+1}$-Whittaker
functions with the characters of irreducible representations of
$\mathfrak{gl}_{\ell+1}$. Indeed,  according to the Kirillov
philosophy,  unitary irreducible representations of a Lie group
$G$ can be obtained by  quantization of coadjoint orbits of $G$
 \cite{K}. In the other direction, a classical limit of the
character of an irreducible representation is given by  a
$G$-equivariant symplectic volume of the corresponding coadjoint
orbit. The precise  relation  in our case is as follows.
The  orbit integral \eqref{orbint} is equal to a
limit of the character \eqref{qzero}  of an irreducible finite-dimensional
  representation of $\mathfrak{gl}_{\ell+1}$
 $$
\int_{\CB_{\ell+1}}\,
e^{\omega_{u_0}+
\sum_{j=1}^{\ell+1}\lambda_jH_{jj}(u_0)}=\Delta(\imath x)\,
\lim_{\epsilon \to
  0}\frac{\chi_{\epsilon^{-1}x_1,\ldots, \epsilon^{-1}x_{\ell+1}}
(e^{\imath \epsilon \lambda_1},\ldots , e^{\imath \epsilon
\lambda_{\ell+1}})}
{\chi_{\epsilon^{-1}x_1,\ldots, \epsilon^{-1}x_{\ell+1}}(1,\ldots ,1)}.
$$
where
$$
\chi_{x_1,\ldots, x_{\ell+1}}
(e^{\imath \lambda_1},\ldots, e^{\imath \lambda_{\ell+1}})=
\frac{1}{\prod_{i<j}(e^{\imath \lambda_i}-e^{\imath \lambda_j})}\,\sum_{s\in
  \mathfrak{S}_{\ell+1}}\,(-1)^{l(s)}\,\prod_{j=1}^{\ell+1}
e^{\imath(x_{s(j)}+\ell+1-s(j))\lambda_j}.
$$
A similar interpretation holds for the elementary $(\ell+1,1)$-Whittaker
functions \eqref{L+1EF}. Note that these functions
can be obtained either by specialization
of the previous formulas or as equivariant symplectic volumes of the
partial flag space $GL_{\ell+1}/P_{\ell,\ell+1}=\IP^{\ell}$.

\begin{prop} The elementary $(\ell+1,1)$-Whittaker function  associated with
the partial flag space $\IP^{\ell}$ has the following integral
representations:
\be\label{orp}
W_{\lambda_1,\ldots ,\lambda_{\ell+1}}(x)=
\int_{\IR_+^{\ell+1}}
 \prod_{j=1}^{\ell+1}dt_j \,e^{\imath t_j\lambda_j}\,\delta(\sum_{j=1}^{\ell+1}t_j-x)=
\int_{\Delta_{\ell}(x)} \prod_{j=1}^{\ell+1}dt_j \,e^{\imath t_j\lambda_j}\, ,
\ee
where $\Delta_{\ell}(x)$ is a  simplex defined by the
equation $\sum_{j=1}^{\ell+1}t_j=x$ in $\IR_+^{\ell+1}$.
\end{prop}

\proof This representation can be derived straightforwardly
from  \eqref{L+1EF} taking into account $\int_0^{\infty}\,dt
e^{-at}=a^{-1}$, $a>0$. $\Box$

The simplex $\Delta_{\ell}$ in \eqref{orp} can be understood as an  image of the
projective space $\IP^{\ell}$ under the $U_1^{\ell}$-momentum
map. In this respect it is an analog of the Gelfand-Zetlin polytope
${\cal P}_{\ell+1}$ in \eqref{GZpoly}.  Indeed the Darboux
coordinates $(T,\theta)$ on an open part $\CO^{(0)}_{u_0}$ of the coadjoint orbit
$\CO_{u_0}$ define the Hamiltonian action of the torus
$T^{\ell(\ell+1)/2}$ acting by rotation on $\theta_{ij}$.  The
momentum map for this action maps $\CO^{(0)}_{u_0}$ onto the
Gelfand-Zetlin polytope ${\cal P}_{\ell+1}$.

Now we provide a similar interpretation of the elementary local
$L$-factors \eqref{levzeroL} as
equivariant symplectic volumes of non-compact symplectic spaces.
Let us  equip the vector space $V=\IC^{\ell+1}$ with the standard
symplectic structure
$$
\omega=\frac{\imath}{2}\sum_{j=1}^{\ell+1}dz^j\wedge d\zb^j,
$$
where $(z^1,\ldots ,z^{\ell+1})$  are complex linear coordinates on $V$.
The standard action of $U_{\ell+1}$ on $V$ is Hamiltonian.
Explicitly the action of the diagonal subgroup $U_1^{\ell+1}\subset U_{\ell+1}$
on $V$ is  generated by  vector fields
$$
v_j=\left(z_j\frac{\pr }{\pr z_j}-\zb_j\frac{\pr }{\pr \zb_j}\right).
$$
The corresponding momenta $H_j$ (i.e. solutions of equations $\iota_{v_j}\omega=dH_j$)
are given by
$$
H_j=\frac{\imath}{2}|z^j|^2,\qquad j=1,\cdots ,(\ell+1).
$$
The $U_1^{\ell+1}$-equivariant volume of $V$ is defined
as the  integral
\be\label{eqvolC}
Z(\lambda_1,\cdots ,\lambda_{\ell+1})=
\frac{1}{(2\pi)^{\ell+1}}\,\int_{\IC^{\ell+1}}
e^{\omega+\sum_{j=1}^{\ell+1}H_j\lambda_j}=
\frac{1}{(2\pi)^{\ell+1}}\,
\int_{\IC^{\ell+1}}\frac{\omega^{\ell+1}}{(\ell+1)!}
e^{\sum_{j=1}^{\ell+1}H_j\lambda_j}.
\ee
\begin{lem} The elementary local $L$-factor
associated with the vector space $V=\IC^{\ell+1}$
is expressed through the equivariant volume \eqref{eqvolC}
as follows:
$$
L_V^{(0)}(s)=Z(\lambda_1-s,\cdots ,\lambda_{\ell+1}-s).
$$
\end{lem}

\noindent{\it Proof}.   Direct calculation of the Gaussian integral
\eqref{eqvol}. $\Box$

Finally let us give an interpretation of the eigenfunction property
\eqref{eigenprop0} of the elementary Whittaker function with respect to the action of
the elementary Baxter operator. We consider the
simplest cases $\ell=0,1$ leaving the general case for another
occasion.  Let $(M,\omega, G,\mu)$ be a quantizable $G$-symplectic space
$M$ with a symplectic form $\omega$  and a fixed momentum map $\mu:
M\to \mathfrak{g}^*$ where $\mathfrak{g}^*$ is dual to
$\mathfrak{g}={\rm Lie}(G)$. Here quantizable means that
one can naturally associate with $(M,\omega,G,\mu)$ a unitary $G$-module $V_M$.
In general the representation $V_M$ associated with the quantizable symplectic
manifold $(M,\omega,G,\mu)$  has a non-trivial decomposition
on irreducible $G$-representation
\be\label{decomp}
V_M=\oplus_{\alpha\in {\rm Rep}(G)}\,V_\alpha\otimes W_\alpha,
\ee
where $W_\alpha=Hom_G(V_\alpha,V_M)=V_\alpha^*\otimes_GV_M$
are multiplicity spaces. Let us recall  a realization of this
decomposition \eqref{decomp} via  symplectic geometry of the underlying
symplectic space $M$ (see e.g. \cite{GLS}). Given
quantizable symplectic spaces $(M_1,G,\omega_1,\mu_1)$,  $(M_2,G,\omega_2,\mu_2)$
with  corresponding unitary $G$-modules $V_{M_1}$ and $V_{M_2}$,
the space of linear $G$-maps $Hom_G(V_{M_1},V_{M_2})$
can be obtained by  quantization of a symplectic space
associated with $(M_1,G,\omega_1,\mu_1)$,  $(M_2,G,\omega_2,\mu_2)$
as follows. Consider the space $(M_1\times
M_2,G,\omega_1-\omega_2,\mu_1-\mu_2)$
where $G$ acts diagonally on $M_1\times M_2$. The vector space of $G$-maps ${\rm
  Hom}_G(V_2,V_1)$ can be obtained by  quantization of the
Hamiltonian reduction of the space $(M_1\times
M_2,G,\omega_1-\omega_2,\mu_1-\mu_2)$ over
 zero values of the $G$-momentum map $\mu^{tot}=\mu_1-\mu_2$.
Let us denote the result of this reduction by $N(M_1,M_2)$.

According  to Kirillov (see e.g. \cite{K})
 an irreducible representation $V_{\lambda}$ of a Lie group $G$
characterized by a weight $\lambda$  shall be  associated with a coadjoint symplectic
orbit $\CO_{\lambda}$ of an element  $\lambda\in \Fg^*$  supplied with
the Kirillov-Kostant  symplectic structure $\omega_\lambda$. Thus for the special case
$M_2=\CO_{\lambda}$ the reduced space $N(M,\CO_{\lambda})$ is given by a factor of
the zero momentum subset
\be\label{defeq1}
\mu_{M}-\mu_{\CO}=0,
\ee
over the diagonal action of $G$. The locus \eqref{defeq1} can be identified with
$M_0=\mu^{-1}_M(\CO_\lambda)$ and thus the
 symplectic counterpart of the multiplicity space
 $Hom_G(V_\lambda,V_{M})$ is obtained by
Hamiltonian reduction  of $M$ over the coadjoint
orbit $\CO_\lambda$
\be\label{sympmult}
N(M,\CO_\lambda)=M//_{\mu_M=\lambda}G:=\mu_M^{-1}(\CO_\lambda)/G.
\ee
The symplectic counterpart of the  decomposition \eqref{decomp}
is given by a realization of $M$ as a stratified
symplectic bundle  (see e.g. \cite{GLS})
with  base $\mathfrak{h}^*$, $\mathfrak{h}={\rm
  Lie}(H)$,  and with fibres over
$\lambda\in \mathfrak{h}^*$ being $H_{\lambda}\times
\CO_\lambda\times \left(M//_{\mu=\lambda}G\right)$.
Here $H_{\lambda}\subset H$ is a subgroup of the Cartan subgroup
$H\subset G$ depending on $\lambda$.  Note that one can have
$H_{\lambda}=\emptyset$.
Let us now apply these considerations to the construction of 
symplectic counterparts to the Baxter operator for low ranks
$\ell=0,1$.

For $\ell=0$ we  consider two $U_1$-modules
$V_n$ and $V_{\IC}=\IC[z]$ such that the generator $H$ of ${\rm Lie}(U_1)$
acts as
$$
H|_{V_n}=n, \qquad H|_{V_{\IC}}=z\pr_z.
$$
The modules $V_n$ and $V_{\IC}$ can be obtained by  quantization
of the symplectic $U_1$-spaces $({\rm pt},U_1,\omega=0,\mu_n=n)$
and $(\IC,U_1,\omega_\IC,\mu=\frac{1}{2}|z|^2)$
where $\omega_{\IC}=\frac{\imath}{2}dz\wedge d\zb$. 
The Hamiltonian action of $U_1$ is given by
$e^{\imath \theta}:\,\,z\longrightarrow e^{\imath \theta}z$.

Recall that the $q=0$ specialization \eqref{q=0Baxter} of the $q$-deformed
Baxter operator is associated with the decomposition of the
product $V_n\otimes \CH$ into irreducible representations
\be\label{decomp0}
V_n\otimes V_{\IC}=\oplus_{m\in \IZ_{\geq 0}}\,V_{n+m}, \qquad
V_{\IC}=\oplus_{m\in \IZ_{\geq 0}} V_m.
\ee
We would like to find a symplectic geometry counterpart of this
decomposition. The  multiplicity spaces are given by
$$
Hom_{U_1}(V_m,V_n\otimes V_\IC)=\IC,\qquad m\geq n,
$$
$$
Hom_{U_1}(V_m,V_n\otimes V_\IC)=0,\qquad m< n.
$$
Using \eqref{sympmult} we obtain for symplectic analogs
$N_{\lambda}=\IC//_{\mu=\lambda}U(1)$  of the  multiplicity
spaces \\ $Hom_{U_1}(V_m,V_n\otimes V_\IC)$,
$$
N_{\lambda\geq 0}={\rm pt},\qquad N_{\lambda<0}=\emptyset.
$$
The symplectic geometry analog  of the decomposition \eqref{decomp0}  is a
bundle over the stratified space $\IR=\IR_{<\lambda}\cup \{\lambda\}\cup
\IR_{>\lambda}$  such that the fibre over $\IR_{<\lambda}$  is  the empty
set, the fibre over $\{\lambda \}$ is a point and the fibre over
 $\IR_{> \lambda}$ is $S^1$. This indeed  provides a
model for $\IC$ via momentum map projection $\mu_{U(1)}\IC\to \IR$.
Classical counterpart of the decomposition \eqref{decomp0} for the
equivariant symplectic volume integral can thus be written in  the
following form:
$$
\int_\IC\,e^{\omega_{\IC}+\tau\mu_{U_1}}\times  e^{-\imath \lambda
\tau}=\int_\IC\,e^{\omega_{\IC}+\tau(\mu_{U_1}-\imath \lambda)}=
\int_{T^*S^1}d\theta d\rho  \,\,\Theta(\rho-\lambda)e^{-\imath \tau \rho}.
$$

Let us now consider the case of $\ell=1$. We have  two $U_2$-modules
$V_{n_1,n_2}$ and $V_{\IC^2}=\IC[z_1,z_2]$ such that diagonal generators
$H_j$, $j=1,2$ of ${\rm Lie}(U_2)$
act as $H_j|_{V_{n_1,n_2}}=n_j$, $H_j|_{V_{\IC^2}}=z_j\pr_{z_j}$.
The Hamiltonian action of $U_2$ on $\IC^2$ is given by
$g:\,\,z_i\longrightarrow \sum_{j=1}^2g_{ij}\,z_j$ and
the corresponding momenta are
$$
(\mu_{\IC^2})_{ij}=\frac{\imath}{2}z_i\zb_j,
$$
with respect to the symplectic structure
$$
\omega_{\IC^2}=\frac{\imath}{2}\left(dz_1\wedge d\zb_1+dz_2\wedge
d\zb_2\right).
$$
The modules $V_{n_1,n_2}$ and $V_{\IC^2}$ can be obtained by  quantization
of the symplectic $U_2$-spaces $(\IP^1,U_2,\omega_{s_1,s_2}=(s_1-s_2)\omega_{FS},
\mu_{\IP^1})$ and $(\IC^2,U_2,\omega_{\IC^2},\mu_{\IC^2})$, where $\omega_{FS}$ is the
Fubini-Studi symplectic form on $\IP^1$,  and the momentum map in 
stereographic coordinates is
$$
\mu_{\IP_1}=\imath \frac{s_1+s_2|z|^2}{(1+|z|^2)}.
$$
The $q=0$ Baxter $\CQ$-operator is associated with the decomposition
of the product $V_{n_1,n_2}\otimes V_{\IC^2}$ into irreducible
$U_2$-representations. The Baxter eigenfunction equation \eqref{l=2br} can be
represented in the following form
\be\label{decomp1}
V_{n_1,n_2} \otimes V_{\IC^2}=\oplus_{m_1+m_2=n_1+n_2+r}
\Theta(m_1-n_1)\Theta(n_1-m_2)\Theta(m_2-n_2) \,\,V_{m_1,m_2}.
\ee
The symplectic geometry counterpart of this decomposition is a
representation of the product of a $U_2$ coadjoint orbit
$\CO_{s_1,s_2}$, $s_1>s_2$  and $\IC^2$ as a bundle over a stratified
space $\mathfrak{u}_1\oplus\mathfrak{u}_1=\IR^2$ with  generic non-empty fibre 
$S^1\times S^1\times \CO_{\lambda_1,\lambda_2}\times
N(s_1,s_2|\lambda_1,\lambda_2)$, where $\lambda_1>\lambda_2$.
The symplectic space $N(s_1,s_2|\lambda_1,\lambda_2)$
 is a symplectic  counterpart  of the multiplicity
space ${\rm Hom}_{U_2}(V_{m_1,m_2},V_{n_1,n_2}\otimes
V_{\IC^2})$. We construct this representation by rewriting the
$U_2$-equivariant symplectic volume integral
\be\label{Int1}
I_{\CO_{s_1,s_2}\times \IC^2}(\xi;\tau_1,\tau_2)=
\int_{\CO_{s_1,s_2}\times \IC^2}\,  \omega_{s_1,s_2}\, d^2z_1\,d^2z_2\,\,
e^{\tau_1\mu_{11}+ \tau_{2}\mu_{22}+ \xi\mu_*},
\ee
where $\mu=\mu_{\IC^2}+\mu_{\IP^1}$ is a momentum map $\CO_{s_1,s_2}\times \IC^2 \to
\mathfrak{u}_2^*$, $\mu_*=(\mu_{\IC^2})_{11}+(\mu_{\IC^2})_{22}$
 and $\omega_{s_1,s_2}$ is the  Kirillov symplectic form  on
$\CO_{s_1,s_2}=\IP^1$. The integral \eqref{Int1}
 can be rewritten as follows:
$$
I_{\CO_{s_1,s_2}\times \IC^2}(\xi;\tau_1,\tau_2)=
\int_{\CO_{s_1,s_2}
\times \IC^2\times \mathfrak{u}_2^*}\, \omega_{s_1,s_2}\,d^2z_1\,d^2z_2\,d^4u\,
\prod_{i,j=1}^2\delta(\mu_{ij}-u_{ij})\,e^{ \tau_1\mu_{11}+\tau_{2}\mu_{22}+ \xi\mu_*}.
$$
A generic  element $u\in \Fg^*$ can be represented as $u=g^{-1}u_0g$, $g\in
U_2/U_1^2$ and  $u_0=\imath \,{\rm
  diag}(\lambda_1,\lambda_2)$ such that $\lambda_1> \lambda_2$ via
a projection $\Fg^*\to \mathfrak{t}^*/W$.
We have the following relation between the integration measures (c.f. the  Weyl
integration formula):
$$
d^4u=\frac{1}{2} \Delta^2(\lambda) \,\,{\rm vol}_{U_2/U_1^2}(g)\wedge \,
d\lambda_1\,\wedge\,d\lambda_2, \qquad
\Delta(\lambda)=\prod_{1\leq i<j\leq 2}(\lambda_i-\lambda_j),
$$
where ${\rm vol}_{U_2/U_1^2}$ is the induced volume form on the
coadjoint orbit $U_2/U_1^2$.
Thus we obtain
$$
I_{\CO_{s_1,s_2}\times \IC^2}(\xi;\tau_1,\tau_2)=
\frac{1}{2}\int_{\IP^1\times \IR^2}\, {\rm vol}_{U_2/U_1^2}(g)\,
 \,d\lambda_1\,d\lambda_2\, \Delta(\lambda)
J_{\CO_s\times \IC^2}(g,\lambda_1,\lambda_2,\tau_1,\tau_2),
$$
where
$$
J_{\CO_{s_1,s_2}\times \IC^2}
(g,\lambda_1,\lambda_2;\xi;\tau_1,\tau_2)=
\Delta(\lambda) \int_{\CO_{s_1,s_2}\times \IC^2}\omega_{s_1,s_2}\,d^2z_1\,d^2z_2\,
\prod_{i,j=1}^2\delta(\mu_{ij}-u_{ij})\,e^{\tau_1\mu_{11}+
\tau_{2}\mu_{22}+ \xi \mu_*}=
$$
$$
=e^{\tau_1 u_{11}(u_0,g)+\tau_{2}u_{22}(u_0,g)}
\Delta(\lambda) \int_{\CO_{s_1,s_2}\times \IC^2}\omega_{s_1,s_2}\,d^2z_1\,d^2z_2\,
\prod_{i,j=1}^2\delta(\mu_{ij}-\lambda_i\delta_{ij})\,
e^{\frac{\imath \xi}{2}(|z_1|^2+|z_2|^2)}\,.
$$
In  the last formula we use $SU(2)$-invariance of the integral.
Thus we obtain
\be\label{l=22}
I_{\CO_{s_1,s_2}\times \IC^2}(\xi;\tau_1,\tau_2)
=\int_{\IR^2}\,d\lambda_1\,d\lambda_2\,J(s_1,s_2;\lambda_1,\lambda_2)\,
F(\lambda_1,\lambda_2;\tau_1,\tau_2),
\ee
where
$$
J(s_1,s_2;\lambda_1,\lambda_2)=\frac{1}{2}
\Delta(\lambda) \int_{\CO_{s_1,s_2}\times \IC^2}
\omega_{s_1,s_2}\,d^2z_1\,d^2z_2\,
\prod_{i,j=1}^2\delta(\mu_{ij}-\lambda_i\delta_{ij})\,
e^{\frac{\imath \xi}{2}(|z_1|^2+|z_2|^2)},
$$
and
\be\label{orbint1}
F(\lambda_1,\lambda_2;\tau_1,\tau_2)=
\int_{\CO_{\lambda_1,\lambda_2}}
 \omega_{\lambda_1,\lambda_2}\,e^{\tau_1u_{11}+-\tau_2u_{22}},
\ee
where we use $\omega_{\lambda_1,\lambda_2}=\Delta(\lambda){\rm
  vol}_{U_2/U_1^2}$.  The integral \eqref{orbint1} is the
$U_2$-equivariant symplectic volume of the coadjoint orbit
$\CO_{\lambda_1,\lambda_2}$ and is  classical counterpart of
the character of the irreducible representation associated with the
coadjoint orbit. On the other hand the function
$J(s_1,s_2;\xi;\lambda_1,\lambda_2)$
shall  be considered to be  classical analog of the
multiplicity function ${\rm dim} \,
{\rm Hom}_{U_2}(V_{m_1,m_2},V_{n_1,n_2}\otimes V_{\IC^2})$.

\begin{prop} One has the following expression for the integral:
\be\label{Int2}
J(s_1,s_2;\xi;\lambda_1,\lambda_2)=\frac{1}{2}\Delta(\lambda)\,
\int_{\CO_{s_1,s_2}\times \IC^2}\,\omega_{s_1,s_2}\, d^2z_1\,d^2z_2\,
\prod_{i,j=1}^2\delta(\mu_{ij}-\lambda_i\delta_{ij})
e^{\frac{\imath \xi}{2}(|z_1|^2+|z_2|^2)}=
\ee
$$
=e^{\imath \xi(s_1+s_2-\lambda_1-\lambda_2)}
\Theta(\lambda_1-s_1)\Theta(s_1-\lambda_2)\Theta(\lambda_2-s_2),
$$
where $\lambda_1>\lambda_2$ and $s_1>s_2$.
\end{prop}

\noindent {\it Proof}. The calculation for general $s_1$, $s_2$ can 
easily be reduced to the case $s_2=0$ and $s_1=s$. In this case the integral can be
calculated using
the representation of the integral over the orbit $\CO=\IP^1$ via the
integral over $\IC^2$
$$
J(s_1,s_2;\xi;\lambda_1,\lambda_2)=\Delta(\lambda)\,
\int_{\IC^2\times \IC^2}\,d^2w_1\,d^2w_2\,\,d^2z_1\,d^2z_2\, \times
$$
$$
\times \prod_{i,j=1}^2\delta(z_i\zb_j+w_i\wb_j-(u_0)_{ij})
\delta(|w_1|^2+|w_2|^2-(s_1-s_2))e^{\frac{\imath \xi}{2}(|z_1|^2+|z_2|^2)}
$$
$$
=\Delta(\lambda)\,
\int_{\IC^2\times \IC^2}\,d^2w_1\,d^2w_2\,\,d^2z_1\,d^2z_2\, \times
\delta(|z_1|^2+|w_1|^2-\lambda_1)
\delta(|z_2|^2+|w_2|^2-\lambda_2)\times
$$
$$
\times \delta^{}(z_1\zb_2+w_1\wb_2)\delta^{}(\zb_1 z_2+\wb_1 w_2)
\delta(|w_1|^2+|w_2|^2-s)e^{\frac{\imath \xi}{2}(|z_1|^2+|z_2|^2)}.
$$
Let us introduce the  variable $\xi=z_1\zb_2+w_1\wb_2$ and integrate
over $\xi$ to  obtain
$$
\Delta(\lambda)\, \int\,d^2w_1\,d^2w_2\,\frac{d^2z_2}{|z_2|^2}\,
\delta\left(\frac{|w_1|^2|w_2|^2}{|z_2|^2}+|w_1|^2-\lambda_1\right)
\delta(|z_2|^2+|w_2|^2-\lambda_2)
\delta(|w_1|^2+|w_2|^2-s) e^{\imath \xi(\lambda_1+\lambda_2-s)}.
$$
Integration over $w_1$ gives
$$
e^{\imath
  \xi(\lambda_1+\lambda_2-s)}\frac{\Delta(\lambda)}{\lambda_2}\,
\int\,d^2w_2\,d^2z_2\,
\delta(|z_2|^2+|w_2|^2-\lambda_2)
\delta\left(\frac{\lambda_1}{\lambda_2}|z_2|^2+|w_2|^2-s\right),
$$
and further integration over $z_2$ provides
$$
e^{\imath \xi(\lambda_1+\lambda_2-s)}
\frac{\Delta(\lambda)}{\lambda_2}\int d^2w_2\, \Theta(\lambda_2-|w_2|^2 )
\delta\left(\frac{\lambda_1}{\lambda_2}(\lambda_2-|w_2|^2)+|w_2|^2-s\right)=
$$
$$
=e^{\imath \xi(\lambda_1+\lambda_2-s)}
\frac{\Delta(\lambda)}{\lambda_2}\int_0^{\infty}dt \Theta(\lambda_2-t)
\delta\left(\lambda_1+\frac{\lambda_2-\lambda_1}{\lambda_2}t-s\right).
$$
Finally after integration over $t$ and taking into account
$\lambda_1>\lambda_2$  we  arrive at 
$$
J(0,s;\lambda_1,\lambda_2)=e^{\frac{\imath \xi}{2}(\lambda_1+\lambda_2-s)}
\Theta(\lambda_2)\Theta(\lambda_1-s)\Theta(s-\lambda_2).
$$
This  finishes the proof of the Proposition. $\Box$

The integral representation \eqref{Int2} is the symplectic integral
representation of the Baxter operator kernel \eqref{Qker}  for
$\ell=1$.

\section{Equivariant symplectic volumes and  nil-Hecke algebras}

In Section 3 we  constructed
an  elementary $\mathfrak{gl}_{\ell+1}$-Whittaker
function and  demonstrated that it provides
a solution of the eigenfunction problem for the  quantum billiard
associated with the Lie algebra $\mathfrak{gl}_{\ell+1}$.
Further, in  Section 5 the elementary $\mathfrak{gl}_{\ell+1}$-Whittaker
functions were identified with $U_{\ell+1}$-equivariant symplectic
volumes of flag spaces $\CB_{\ell+1}=GL_{\ell+1}/B$.
Thus we have managed to express eigenfunctions  of the quantum
billiard  via equivariant symplectic volumes of flag spaces.
This  connection between quantum billiard eigenfunctions and
equivariant symplectic volumes can be considered as a
manifestation of a general relation between
quantum many-body integrable systems,  representation theory of the Hecke algebras  and
(generalized) equivariant cohomology of $G$-spaces (see e.g. \cite{CG}, \cite{Ch}
for  detailed discussions).  Thus
a description of $G$-equivariant cohomology $H_G(\CB,\IC)$
 of the flag space $\CB=G/B$ in terms of
the nil-Hecke algebra $\CH^{nil}(G^{\vee})$ associated with
the dual Lie group $G^{\vee}$  was proposed
in \cite{KK} (see \cite{BGG} for non-equivariant  case).
Below we consider the special case  $G=G^{\vee}=GL_{\ell+1}$
and demonstrate that the results of \cite{KK} are  compatible with
the considerations of our previous Sections.

Let us first recall  some standard facts on
the nil-Hecke algebra associated with the Lie algebra
$\mathfrak{gl}_{\ell+1}$ (see e.g. \cite{CG}).
Let  $\Phi$ be the root system of $\mathfrak{gl}_{\ell+1}$.
We identify the diagonal  Cartan subalgebra
$\mathfrak{h}\subset \mathfrak{gl}_{\ell+1}$
with $\IR^{\ell+1}$ and fix a
basis $\{e_i\}$, $i=1,\ldots \ell+1$ in $\mathfrak{h}$,
orthonormal with respect to the bilinear form
$(\,,\,)$ induced by the Killing form on ${\mathfrak{gl}}_{\ell+1}$. We define coroots as
$\alpha^{\vee}=2\alpha/(\alpha,\alpha)$.
Using the Killing form we identify
$\mathfrak{h}$ and its dual $\mathfrak{h}^*$. The
positive simple roots of $\Phi$ are given then  by
$\alpha_i=e_{i+1}-e_{i}$, $i=1,\ldots ,\ell$. The Weyl group
acts on $\mathfrak{h}$ by reflections
\be
s_{\alpha}:\,\,y\longrightarrow y-\<\alpha,y\>\alpha^{\vee}, \qquad
\alpha\in \Phi,
\ee
and is isomorphic to the permutation group
$\mathfrak{S}_{\ell+1}$ generated by  elementary permutations
$s_i=\sigma_{i,i+1}$  acting on the basis
$\{e_j\}$, $j=1,\ldots, (\ell+1)$ in $\IR^{\ell+1}$ via  permutations of the  indexes.
Let $s_{max}\in \mathfrak{S}_{\ell+1}$ be the element 
with  reduced decomposition of maximal length. It can be written for example as follows
\be\label{maxelm}
s_{max}=s_1s_2s_1s_3s_2s_1\cdots s_{\ell}s_{\ell-1}\cdots s_1.
\ee

\begin{de} The nil-Hecke algebra $\CH^{nil}_{\ell+1}$
associated with the root system $\Phi$ of  $\mathfrak{gl}_{\ell+1}$
is an associative algebra generated by  $R_s$, $s\in
\mathfrak{S}_{\ell+1}$, $D_{i}$, $i=1,\ldots, \ell+1$ and a central
element $c$  with the following relations:
\be\label{d1}
R_s\cdot R_{s'}=\begin{cases} R_{ss'} & \mbox{if}\quad  l(s)+l(s')=l(ss')\\
0 &\mbox{if}\quad  l(s)+l(s')\neq l(ss'),
\end{cases}
\ee
\be\label{d3}
D_{i}\,D_{j}-D_{j}\,D_{i}=0,
\ee
\be\label{d4}
R_i\cdot D_{j}-D_{s_{i}(j)}\cdot R_i\,=\,c\,\delta_{ij}.
\ee
\end{de}
The  generators $R_i:=R_{s_i}$ corresponding to
elementary permutations  $s_i$ satisfy the following relations:
$$
R_iR_j=R_jR_i,\qquad |i-j|\geq 2,
$$
\be\label{nilCox}
R_iR_{i+1}R_i=R_{i+1}R_iR_{i+1},
\ee
$$
R_i^2=0,
$$
and generate the nil-Coxeter subalgebra $W^{nil}\subset \CH^{nil}_{\ell+1}$.  In particular
operators
\be\label{exdec}
R_s=R_{s_{i_1}}R_{s_{i_2}}\cdots R_{s_{i_l}},\qquad
s=s_{i_1}s_{i_2}\cdots s_{i_l},
\ee
are zero unless the product representation of $s$
in terms of elementary reflections has  minimal length.
The corresponding element \eqref{exdec}
is independent of the choice of the product decomposition of $s$.

The  nil-Hecke algebra $\CH^{nil}_{\ell+1}$ is a semidirect
product of $W^{nil}$ and the  algebra $\IC[\mathfrak{h}]$ of
polynomial functions on $\mathfrak{h}$.
The center $\CZ_{\ell+1}$ of  $\CH^{nil}_{\ell+1}$  is isomorphic to the algebra
of $\mathfrak{S}_{\ell+1}$-symmetric
polynomials of the generators $D_{i}$, $i=1,\ldots, (\ell+1)$.
Irreducible  representations of $\CH^{nil}_{\ell+1}$ can be characterized
by their  central characters i.e.  homomorphisms  $\CZ_{\ell+1}\to
\IC$ and a homomorphism $\CZ_{\ell+1}\to
\IC$ is  uniquely defined by
$\mathfrak{S}_{\ell+1}$-orbits of an element $\lambda\in \mathfrak{h}$.

Let $\CW_{\lambda}$ be the linear space generated by the exponential functions
$\Psi_{\lambda,s}(x)=\exp(\imath \<s\cdot \lambda, x\>)$ on $\mathfrak{h}^*$ and let
$\{s\cdot \lambda|s\in \mathfrak{S}_{\ell+1}\}$  be the orbit through a
fixed element $\lambda\in \mathfrak{h}$.
$\CW_{\lambda}$ may be identified with the space of common
eigenfunctions of the ring $\mathfrak{Diff}_c^{\mathfrak{S}_{\ell+1}}$
of $\mathfrak{S}_{\ell+1}$-invariant
differential operators  on $\mathfrak{h}^*$ with constant coefficients
\be\label{setofeq}
\CD_{\chi}\cdot
\Psi_{\lambda,s}(x)=\,\chi(\imath \lambda)\,\Psi_{\lambda,s}(x),
\ee
where $\CD_{\chi}\in \mathfrak{Diff}_c^{\mathfrak{S}_{\ell+1}}$
is  the differential operator  corresponding
to $\chi\in \IC[\mathfrak{h}]^{\mathfrak{S}_{\ell+1}}$ via the
isomorphism
$\mathfrak{Diff}_c^{\mathfrak{S}_{\ell+1}}=\IC[\mathfrak{h}]^{\mathfrak{S}_{\ell+1}}$.
The space of solutions  has cardinality  $|\mathfrak{S}_{\ell+1}|$
and is given by linear combinations of the exponents.
The  following integral/differential  operators  provide a realization of the
irreducible representation $\pi:\CH^{nil}_{\ell+1}\to {\rm End}(\CW_{\lambda})$:
\be\label{conR}
\pi(R_j)\cdot g(x)=\int_0^{(x,\alpha_j)}\,g(x-t\alpha_j)dt,
\ee
\be\label{conD}
\pi(D_j)\cdot g(x)=c\,\frac{\pr g(x)}{\pr x_j},
\ee
where $g\in \CW_{\lambda}$.
Equivalently these operators can be written as follows:
\be\label{rep2}
\pi(R_j)\,e^{\imath \<\lambda,x\>}=
\frac{e^{\imath (\lambda,x)}-e^{\imath (s_j\cdot \lambda,x)}}
{\imath \<\lambda,\alpha_j\>},
\ee
\be\label{rep1}
\pi(D_{j})\, e^{\imath \<\lambda,x\>}=
\imath  c\<\lambda,\a_j\>e^{\imath \<\lambda,x\>}.
\ee
From now on we shall abuse notations  and use the symbols
$R_j$ and $D_j$ for the images $\pi(R_j)$ and  $\pi(D_j)$ of the generators in the
representation $\CW_{\lambda}$.

We  call a vector  $\Psi^{(0)}_{\lambda}\in \CW_{\lambda}$
to be of class one if it satisfies the conditions
\be\label{classone}
R_i\Psi^{(0)}_{\lambda}(x)=0, \qquad i=1,\ldots, (\ell+1).
\ee
This condition is a nil-Hecke analog of the sphericity
condition for affine Hecke algebras \cite{CG} and uniquely defines a
class one vector in  $\CW_{\lambda}$.  Taking into account that
$l(s_is_{max})<l(s_i)+l(s_{max})$ and the relations
\eqref{d1}  for $R_i\cdot R_j$,
 we infer that the class one vector can be written as follows
\be\label{clone1}
\Psi^{(0)}_{\lambda}(x)=R_{s_{max}}\,e^{\imath \<\lambda, x\>}.
\ee

\begin{prop} The class one vector in the representation $\CW_{\lambda}$ of
  $\CH^{nil}_{\ell+1}$ solves the eigenfunction problem of the quantum
  billiard \eqref{eigenf1}, \eqref{BC} associated with the root system of
  $\mathfrak{gl}_{\ell+1}$.
\end{prop}

\proof Let us note that it is possible to write $s_{max}$ as a minimal length
product of the elementary reflections starting with any $s_i$, and thus obtain
\be\label{clone11}
\Psi^{(0)}_{\lambda}(x)=R_{s_i}R_{s_*}\,e^{\imath \<\lambda, x\>}.
\ee
Now using the realization \eqref{conR} for $R_i:=R_{s_i}$ we infer
that
\be\label{bcond1}
\Psi^{(0)}_{\lambda}(x)|_{x_{j}=x_{j+1}}=0, \qquad j=1,\ldots ,\ell.
\ee
Thus taking into account \eqref{setofeq} we prove the statement of
the Proposition. $\Box$

In previous Sections  we show that 
the elementary $\mathfrak{gl}_{\ell+1}$-Whittaker function
$\Psi^{(0)}_{\lambda}(x)$ \eqref{sumover}, \eqref{sumoverone} is a quantum billiard  eigenfunction.
Now we explain how the relation \eqref{eqvol} between  elementary
$\mathfrak{gl}_{\ell+1}$-Whittaker functions
and $U_{\ell+1}$-equivariant symplectic volumes of flag spaces
$\CB_{\ell+1}$  arises via a relation between nil-Hecke algebras
and equivariant cohomology of flag manifolds.

The class one vector \eqref{clone1} solving the set of  equations
\eqref{setofeq} can be expressed through  $U_{\ell+1}$-equivariant symplectic
volume of the flag space $\CB_{\ell+1}$. This follows from a general
connection between equivariant cohomology of flag spaces and nil-Hecke
algebras which we review below following \cite{KK} (see also
\cite{CG}).   Let us start with the non-equivariant case.
Recall that the flag space $\CB_{\ell+1}$ has  a Schubert cell decomposition,
with cells $\CO_s$  enumerated by elements $s\in \mathfrak{S}_{\ell+1}$
of the corresponding Weyl group. The homology classes
$[\overline{\CO}_s]$ provide a  basis in the homology
groups $H_*(\CB_{\ell+1},\IQ)$ and we denote by $\sigma_s\in
H^*(\CB_{\ell+1})$  the dual cohomology classes.
Thus for example $\sigma_{id}=1$ is  dual to the fundamental class of
$\CB_{\ell+1}$ and $\sigma_{s_{max}}$ is dual to the unique zero
homology class. The following orthogonality relation follows from the
standard intersection product relations between  Schubert classes
$$
(\sigma_s\,\sigma_{s'})=\delta_{s,s_{max}s'},\qquad
(\omega):=\int_{\CB_{\ell+1}}\omega.
$$

There exists another description of the cohomology groups of flag
manifolds due to Borel
\be\label{BorelISO}
H^*(\CB_{\ell+1},\IC)=\IC[y_1,\ldots ,y_{\ell+1}]/J,
\ee
where the ideal $J$ is generated by
$f\in \IC[y_1,\ldots ,y_{\ell+1}]^{\mathfrak{S}_{\ell+1}}$, $f(0)=0$.
Here the element $\<\varphi,y\>$, $\varphi \in \mathfrak{h}^*_{\IZ}$ in the
right hand side of \eqref{BorelISO}
corresponds to the image of the first Chern class $c_1(\CL_{\varphi})$
of the line bundle $\CL_{\varphi}$ associated with the weight
$\varphi$. Let us denote by $\overline{f}$ the  element of  $H^*(\CB_{\ell+1},\IC)$
corresponding to $f$ under the isomorphism \eqref{BorelISO}.

The pairing of an element $\overline{f}\in H^*(\CB_{\ell+1},\IC)$ with the fundamental
class of $\CB_{\ell+1}$ can be expressed in the integral form
\be\label{pair1}
(\bar{f})=\frac{1}{(2\pi \imath)^{\ell+1}}\,\,\oint_{\CC_0} dy_1\,\cdots
dy_{\ell+1}\,\frac{f(y)}{\prod_{j=1}^{\ell+1}e_j(y)}, \qquad
\bar{f}\in H^*(\CB_{\ell+1}),
\ee
where $e_j(y)$ are  elementary symmetric functions
\be\label{elsym}
\sum_{j=0}^{\ell+1}(-1)^jz^{\ell+1-j}\,e_j(y)=\prod_{j=1}^{\ell+1}(z-y_j),
\ee
and the integration domain  $\CC_0$ encloses the poles of the denominator.

The representatives of the dual Schubert classes $\sigma_w\in
H^*(\CB_{\ell+1})$ in the polynomial ring \\
$\IQ[y_1,\cdots ,y_{\ell+1}]$ are called Schubert polynomials.  For
example one can take
\be\label{smaxrep}
\sigma_{s_{max}}(y)=\frac{1}{|\mathfrak{S}_{\ell+1}|}\prod_{1\leq
  i<j\leq \ell+1}(y_i-y_j).
\ee

The cohomology $H^*(\CB_{\ell+1})$ of the flag space admits the structure of a module over
nil-Hecke algebra $\CH^{nil}_{\ell+1}$. Precisely, the following operators
\be\label{rep11}
\widetilde{R}_i\cdot P(y)=\frac{P(y)-P(s_i\cdot
  y)}{\imath \<\alpha_i,y\>},
\ee
\be\label{rep12}
\widetilde{D}_i\cdot P(y)=\imath c\,\<\alpha_i,y\>\, P(y),
\ee
define an action of the nil-Hecke algebra
$\CH^{nil}_{\ell+1}$ on the space of polynomials
$P\in \IC[y_1,\cdots, y_{\ell+1}]$.
This action commutes with the multiplication on
$\mathfrak{S}_{\ell+1}$-invariant polynomials and thus
descends to the Borel realization \eqref{BorelISO} of
$H^*(\CB_{\ell+1})$ to provide  an action of $\CH^{nil}_{\ell+1}$ on
$H^*(\CB_{\ell+1})$. The action of the nil-Hecke algebra
provides expressions for  generic Schubert polynomials, thus 
\be\label{Sch}
\sigma_{s}(y)=\widetilde{R}_{s^{-1}s_{max}}\,\sigma_{s_{max}}(y).
\ee
The pairing of the  product of an element $\bar{f}\in H^*(\CB_{\ell+1})$
 and a Schubert  class $\sigma_s$ with the fundamental class
 $[\CB_{\ell+1}]$ is given by
\be\label{pair2}
(f\,\sigma_s)=\widetilde{R}_{s^{-1}s_{max}} \cdot f|_{y=0}.
\ee

Let $\omega_i=\bar{y_i}$ be  generators of $H^2(\CB_{\ell+1})$.
Then, using  the general formula \eqref{pair2}
in the special case  $f(x)=\exp (\imath \sum_{j=1}^{\ell+1}y_jx_j)$,
we obtain the following representation for the symplectic volume of
$\CB_{\ell+1}$:
\be\label{svol}
Z_{\CB_{\ell+1}}(x):=\int_{\CB_{\ell+1}}\,
e^{\sum_{j=1}^{\ell+1}\imath \omega_jx_j}=\widetilde{R}_{s_{max}}\cdot
e^{\sum_{j=1}^{\ell+1}\imath y_jx_j}|_{y=0}.
\ee

\begin{ex} For $\ell=1$ the element of the
  Weyl group with maximal length is  $s_{max}=s_1$, interchanging
  $y_1$ and $y_2$. Thus we have
$$
Z_{\CB_2}(x)=\int_{\IP^1}\,e^{\imath \omega_1x_1+
\imath \omega_2x_2}=\widetilde{R}_{s_1}\cdot
e^{\imath y_1x_1+\imath y_2x_2}|_{y=0}=
\frac{e^{\imath  y_1x_1+\imath  y_2x_2}-
e^{\imath  y_2x_1+\imath  y_1x_2}}{\imath (y_1-y_2)}|_{y=0}
=(x_1-x_2).
$$
\end{ex}

Now let us consider  an equivariant analog of the
relation between the symplectic volume of $\CB_{\ell+1}$ and representation theory
of $\CH^{nil}_{\ell+1}$. Let $T_{\ell+1}\subset U_{\ell+1}$ be the diagonal Cartan
torus. The $T_{\ell+1}$-equivariant
cohomology of the flag space $\CB_{\ell+1}$ has the following
description  generalizing \eqref{BorelISO}:
\be\label{BorelISOeq}
H_{T_{\ell+1}}^*(\CB_{\ell+1},\IC)
=\IC[y_1,\ldots ,y_{\ell+1},\lambda_1,\cdots ,\lambda_{\ell+1}]/J,
\ee
where the ideal $J$ is generated by  functions $f(y)-f(\lambda)$
where $f(y_1,\cdots ,y_{\ell+1})$ is an  arbitrary
$\mathfrak{S}_{\ell+1}$-symmetric polynomial function. 
Linear functions $\<\omega,y\>=\sum_{j=1}^{\ell+1}\omega_jy_j$,
 correspond to $T_{\ell+1}$-equivariant
extensions of the first Chern classes $c_1(\CL_{\omega})$
of the line bundles $\CL_{\omega}$ associated with the integral weights
$\omega$ of the Lie algebra $\mathfrak{gl}_{\ell+1}$. The $U_{\ell+1}$-equivariant cohomology
$H^*_{U_{\ell+1}}(\CB_{\ell+1})$ can be similarly
realized as follows:
$$
H_{U_{\ell+1}}^*(\CB_{\ell+1},\IC)
=\IC[y_1,\ldots ,y_{\ell+1}]\otimes \IC[\lambda_1,\cdots
,\lambda_{\ell+1}]^{\mathfrak{S}_{\ell+1}}/J,
$$
where the ideal $J$ is the same as above.
In the equivariant case the  pairing with the fundamental cycle $[\CB_{\ell+1}]$
is given by
\be\label{pair4}
(\bar{f})_{\lambda}=\frac{1}{(2\pi \imath)^{\ell+1}}\, \oint_{\CC} dy_1\,\cdots
dy_{\ell+1}\,\frac{f(y,\lambda)}{\prod_{j=1}^{\ell+1}(e_j(y)-e_j(\lambda))}, \qquad
\bar{f}\in H_{T_{\ell+1}}^*(\CB_{\ell+1}),
\ee
where $e_j(z)$ are elementary symmetric polynomials \eqref{elsym} and
the integration domain $\CC$ encloses all singularities of the denominator.
The integral reduces to  a sum over the residues and is given by
$$
(\bar{f})_{\lambda}=\frac{1}{\Delta(\lambda)}
\sum_{s\in \mathfrak{S}_{\ell+1}} (-1)^{l(s)}f(s\cdot
y)|_{y_i=\lambda_i}, \qquad \Delta(\lambda):=\prod_{1\leq i,j\leq
  \ell+1} (\lambda_i-\lambda_j).
$$

Define a pairing between  the space $\CF(\mathfrak{h}^*)$
of analytic functions on $\mathfrak{h}^*$ and
polynomial functions $\IC[\mathfrak{h}]$ as  follows:
\be\label{pair3}
((f,g))=\left[f(\pr_y)\cdot g(y)\right]_0, \qquad f\in
\IC[\mathfrak{h}],\quad g\in \CF(\mathfrak{h}^*),
\ee
where $[f]_0$ is the constant term of the  series expansion of $f$
in a neighborhood of $y=0$. The space
$\CW_{\lambda}$ defined by \eqref{setofeq} is dual to the cohomology
groups given by \eqref{BorelISOeq} with respect to the pairing
\eqref{pair3} i.e.
$$
\CW_{\lambda}=\{g\in \IC((\mathfrak{h}^*))|
(f(\pr_y)-f(\lambda))\cdot g(y)=0,\, \mbox{ for
  any }\,f\in J\},
$$
where $J\subset\IC[\mathfrak{h}]$ is the subspace of
$\mathfrak{S}_{\ell+1}$-invariant polynomials.
Note that operators $\widetilde{R}_i$ and $\widetilde{D}_j$ given by
\eqref{rep11}, \eqref{rep12} are
conjugate to the operators $R_i$ and $D_j$ given by \eqref{conR},
\eqref{conD}  with respect to the pairing \eqref{pair3}.
Let us  define the basis $\{\CD_s\}$, $s\in \mathfrak{S}_{\ell+1}$
in $\CW_{\lambda}$ to be dual to the basis of the Schubert polynomials
$\{\sigma_s\}$ in $H^*(\CB_{\ell+1})$.

\begin{prop} One has
\be\label{id2}
(\bar{f})_{\lambda}=\widetilde{R}_{s_{max}}\cdot f(y)|_{y_i=\lambda_i}.
\ee
\end{prop}

\proof First we have
$$
(\bar{f}\,\sigma_{s_{max}})_{\lambda}=\frac{1}{|\mathfrak{S}_{\ell+1}|}\frac{1}{(2\pi
\imath)^{\ell+1}}\,\,
\oint dy_1\cdots dy_{\ell+1}
\frac{f(y)\,\prod_{i<j}(y_i-y_j)}{\prod_{j=1}^{\ell+1}(e_j(y)-e_j(\lambda))}=
\frac{1}{|\mathfrak{S}_{\ell+1}|}\sum_{s\in \mathfrak{S}_{\ell+1}}
f(y)|_{y_{s(i)}=\lambda_i},
$$
where in the last step we replace the integral by the sum over
residues. Now we have
$$
(f)_{\lambda}=((f,\CD_{id}))=((f,R_{s_{max}}\CD_{s_{max}}))=
((\widetilde{R}_{s_{max}}f,\CD_{s_{max}})),
$$
where we use the dual version of \eqref{Sch} for $s={\rm id}$.
Using the duality between $\sigma_{s_{max}}$ and $\CD_{s_{max}}$ we
obtain
$$
(f)_{\lambda}=((\widetilde{R}_{s_{max}}\,f)\,\sigma_{s_{max}})_{\lambda}=
\widetilde{R}_{s_{max}}\, f(y)|_{y_i=\lambda_i}.
$$
Here  at the last step we take into account that  $\widetilde{R}_{s_{max}}\, f(y)$
is a symmetric function.  $\Box$

\begin{prop}
Suppose  that $\omega_i(\lambda)$ are the equivariant
cohomology classes represented by the generator 
$y_i$ in the Borel realization \eqref{BorelISOeq}.
Then we have for equivariant  symplectic volume
\be\label{eqsvol}
Z(x,\lambda):=\int_{\CB_{\ell+1}}\,e^{\sum_{j=1}^{\ell+1}\imath \omega_j(\lambda) x_j}=
(e^{\sum_{j=1}^{\ell+1}\imath y_jx_j})_{\lambda}=\widetilde{R}_{s_{max}}\cdot
e^{\sum_{j=1}^{\ell+1}\imath y_jx_j}|_{y_i=\lambda_i}=
\ee
\be\nonumber
=R_{s_{max}}\cdot
e^{\sum_{j=1}^{\ell+1}\imath \lambda_jx_j}.
\ee
\end{prop}

\proof The result follows from \eqref{id2}
applied to the special case   $f(y)=\exp(\imath
\sum_{j=1}^{\ell+1}y_jx_j)$.  $\Box$

\begin{ex} For  $\ell=1$ we have  $s_{max}=s_1$ interchanging
  $y_1$ and $y_2$. Thus for the equivariant symplectic volume
of $\CB_2=\IP^1$ we obtain
\be
Z_{\CB_2}(x,\lambda)=\int_{\IP^1}\,e^{\imath(\omega_1+\lambda_1H_1)x_1+
\imath(\omega_2+\lambda_2H_2)x_2}=R_{s_1}\cdot
e^{\imath y_1x_1+\imath y_2x_2}|_{y=\lambda}=
\ee
$$
\frac{e^{\imath y_1x_1+\imath y_2x_2}-
e^{\imath y_2x_1+\imath y_1x_2}}{\imath(y_1-y_2)}|_{y=\lambda}
=\frac{e^{\imath\lambda_1x_1+\imath\lambda_2x_2}-
e^{\imath\lambda_2x_1+\imath\lambda_1x_2}}
{\imath(\lambda_1-\lambda_2)}.
$$
\end{ex}

Comparing \eqref{eqsvol} and \eqref{clone1} we obtain the equality of the
class one vector $\Psi^{(0)}_{\lambda}\in \CW_{\lambda}$ and the
equivariant symplectic volume $Z_{\CB_{\ell+1}}$
\be\label{equivONE}
\Psi^{(0)}_{\lambda}(x)=Z_{\CB_{\ell+1}}(x,\lambda).
\ee
Finally note that the expression \eqref{clone1}
for the elementary $\mathfrak{gl}_{\ell+1}$-Whittaker function
$\Psi^{(0)}_{\lambda}(x)$
considered as an eigenfunction of the quantum billiard
provides an integral representation for $\Psi^{(0)}_{\lambda}(x)$. For
example for $\ell=1$ the representation \eqref{clone1} with
$R_{s_1}$ realized by the integral operator \eqref{conR}
gives the following:
\be\label{loneex}
\Psi^{(0)}_{\lambda_1,\lambda_2}(x_1,x_2)=R_{s_1}\cdot
e^{\imath(\lambda_1x_1+\lambda_2x_2)}=
\int_0^{x_2-x_1}e^{\imath((x_1+t)\lambda_1+(x_2-t)\lambda_2)}dt=
\ee
$$
=\int_{x_1}^{x_2}e^{\imath(t\lambda_1+(x_2+x_1-t)\lambda_2)}dt.
$$
This integral representation shall be considered as 
 the special  case ${\mathfrak{g} = \mathfrak{gl}}_1$
  of an elementary version of the integral representation
of a $\mathfrak{g}$-Whittaker function introduced in \cite{GLO1}  for
an arbitrary semisimple Lie algebra. Note that the integral representation
introduced in \cite{GLO1}   shares with the
Givental integral representation the property of positivity.
In the $\ell=1$ case the integral representation \eqref{loneex}  accidentally
coincides with the Givental representation \eqref{giv00}.

\section{Elementary  Whittaker function as a matrix element}

In the previous Sections we prove that 
$U_{\ell+1}$-equivariant symplectic volume of flag
space $\CB_{\ell+1}=GL_{\ell+1}/B$ is  expressed through  the elementary
$\mathfrak{gl}_{\ell+1}$-Whittaker function. 
This identification may be considered as  a limit of the
identification of  classical $\mathfrak{gl}_{\ell+1}$-Whittaker
functions with $S^1\times U_{\ell+1}$-equivariant volumes
of the spaces of holomorphic maps of a
two-dimensional disk $D$ into $\CB_{\ell+1}$
\cite{GLO8}. Recall that the 
$S^1\times U_{\ell+1}$-equivariant volume of the space of holomorphic
maps $D\to \CB_{\ell+1}$  can be  succinctly described as a correlation function
for the equivariant type $A$ topological sigma model \cite{GLO8}, thus
providing an infinite-dimensional integral representation of 
the classical
Whittaker function. It was argued in \cite{GLO8}
that the mirror dual description in terms of the type $B$ equivariant
topological Landau-Ginzburg sigma model leads to a finite-dimensional
integral representation  of
classical $\mathfrak{gl}_{\ell+1}$-Whittaker function. This
finite-dimensional integral representation
has an interpretation as a  matrix element of an infinite-dimensional
representation  of $\CU\mathfrak{gl}_{\ell+1}$ written in an integral form.
This pair of infinite-dimensional and finite-dimensional integral
representations of classical Whittaker functions was advocated in
\cite{GLO8} as a manifestation of the  local Archimedean  Langlands
correspondence.

In this  Section we provide an elementary analog of the mirror
symmetry by constructing a representation of the elementary
$\mathfrak{gl}_{\ell+1}$-Whittaker function \eqref{levzeroW}
as a matrix element of a monoid $GL_{\ell+1}(\CR)$ where $\CR$
is the tropical semifield (for the definition of a tropical semifield
see e.g. \cite{MS}, \cite{IMS}). The correspondence between
realizations of the elementary $\mathfrak{gl}_{\ell+1}$-Whittaker function
as a $U_{\ell+1}$-equivariant volume of $\CB_{\ell+1}$
and as  a matrix element of  a representation of
the monoid $GL_{\ell+1}(\CR)$
can be considered as  an  elementary analog of the local Archimedean
Langlands correspondence.

Let us first discuss some special features of the integral representations \eqref{giv}
of the $\mathfrak{gl}_{\ell+1}$-Whittaker function.
The integral form \eqref{giv} of the Whittaker function arises
from the general expression \eqref{Wf}
using a particular representation 
of $\CU\mathfrak{gl}_{\ell+1}$.
Let $\chi_{\underline{\lambda}}$ be a character of $B_-$ given by
\eqref{chardef}. Then the
representation is realized in  a subspace $\CV_{\underline{\lambda}}\subset
{\rm Ind}_{B_-}^{GL_{\ell+1}}\chi_{\underline{\lambda}}$
of $B_-$-equivariant functions supported at
$N_+^>\times
B_-\subset GL_{\ell+1}$ where $N_+^>$ is  the subset of positive
elements of the maximal unipotent subgroup  $N_+\subset GL_{\ell+1}$.  Recall that
positive elements of $GL_{\ell+1}(\IR)$ are the elements realized  in
the standard matrix representation by positive matrices
i.e.  matrices with all  minors  positive (see e.g. \cite{Lu2}).
Similarly, positive elements of $N_+$ are the elements realized in
the standard matrix representation by
 matrices with all  non-identically zero minors  positive.
The following parametrization of  positive subset $N_+^>$
was introduced in \cite{GKLO} using the
results of \cite{Lu2} (see also \cite{BFZ}).
Let $\epsilon_{i,j}$ be  the  elementary
$(\ell+1)\times(\ell+1)$- matrix with unit  in the $(i,j)$-place and
all other entries zero.
 Consider the  set consisting of the diagonal matrices
$$
 U_k=\sum\limits_{i=1}^ke^{-x_{k,i}}\epsilon_{i,i}\,
 +\,\sum\limits_{i=k+1}^{\ell+1}\epsilon_{i,i}\,,
$$
and of their upper-triangular deformations
\be\label{deformed}
\widetilde{U}_k\,=\,\sum\limits_{i=1}^{k}
e^{-x_{k,i}}\epsilon_{i,i}\,
+\,\sum\limits_{i=k+1}^{\ell+1}\epsilon_{i,i}
+\,\sum\limits_{i=1}^{k-1}e^{-x_{k-1,i}}\epsilon_{i,i+1}\,.
\ee
The factorized parametrization of $N_+^>$  follows from the
fact that the image of any generic unipotent element
  $v\in N^>_+$ in the tautological representation
$\pi_{\ell+1}:\mathfrak{gl}_{\ell+1}\to End(\mathbb{C}_{\ell+1})$
 can be represented in the form
 \be\label{givpar}
  \pi_{\ell+1}(v)=\tilde{U}_{2}U_{2}^{-1}\tilde{U}_{3}U_{3}^{-1}
  \cdots
  \tilde{U}_{\ell}U_{\ell}^{-1}\tilde{U}_{\ell+1},
 \ee
where we assume that $x_{\ell+1,i}=0,\,\,\,\,i=1,\ldots,\ell+1$.
Using this parametrization
the following realization of $\CU\mathfrak{gl}_{\ell+1}$ by
differential operators was constructed in \cite{GKLO}.
\begin{prop}
The following differential operators define a realization of
representation $\pi_{\lambda}$ of $\mathfrak{gl}_{\ell+1}$ in
${\cal V}_{\mu}$ in the space of functions on $N^>_+$:
 \bqa\label{GivRep}
  E_{i,i}&=& \mu_i\,-\,
  \sum_{k=1}^{i-1}\frac{\partial}{\partial x_{\ell+1+k-i,k}}\,
  +\,\sum_{k=i}^{\ell}\frac{\partial}{\partial x_{k,i}},\nonumber \\
  E_{i,i+1}&=& -\sum_{n=1}^ie^{x_{\ell-i+n,\,n}-x_{\ell+1-i+n,\,n}}
  \sum_{k=1}^n\Big\{\frac{\pr}{\pr x_{\ell-i+k,\,k}}\,
  -\,\frac{\pr}{\pr x_{\ell-i+k,\,k-1}}\Big\}\,,
 \label{GGrep}\\
  E_{i+1,i}&=&\sum_{n=1}^{\ell+1-i}e^{x_{n+i,\,i+1}-x_{k+i-1,\,i}}
  \Big[\mu_i-\mu_{i+1}\,+\,\sum_{k=1}^n\Big\{
  \frac{\pr}{\pr x_{i+k-1,\,i}}\,-\,\frac{\pr}{\pr x_{i+k-1,\,i+1}}
  \Big\}\Big]\,,
  \nonumber
 \eqa
where $\mu_k=\imath\hbar^{-1}\lambda_k-\rho_k$. Here  
$E_{i,j}=\pi_{\lambda}(e_{i,j})$, $i,j=1,\ldots,\ell+1$ are the images of the standard
generators $e_{ij}$ of $\mathfrak{gl}_{\ell+1}$ satisfying the
relations
$$
[e_{ij},e_{km}]=\delta_{jk}e_{im}-\delta_{mi}e_{kj}.
$$
\end{prop}
The matrix element \eqref{Wf} written explicitly using this
realization of the representation $\CV_{\mu}$ is given by
the Givental integral formula \eqref{giv}. Note that
the matrix element \eqref{Wf} is
defined for representations such that the action of the
Cartan subalgebra can be integrated to an action of the corresponding
group.
The integration of the action of the whole Lie algebra
$\mathfrak{gl}_{\ell+1}$ is not necessary and actually is not possible
for the representation leading to \eqref{giv}. This is obvious
taking into account that the integration in \eqref{giv} is  over a
subset $N_+^>\subset N_+$ of positive elements of $N_+$.  Thus
the function \eqref{giv}  does not naturally extend to
a function on the group for which
\be\label{equivprop}
\Psi_\lambda(n_-gn_+)=\psi^+(n_+)\psi^-(n_-)\Psi_\lambda(g),\qquad
g\in GL_{\ell+1}(\IR),\quad n_{\pm}\in N_{\pm},
\ee
where the characters of $N_{\pm}$ are the Lie group version of the
characters in \eqref{Wclass}:
\be\label{intchar}
\psi^{\pm}(n_{\pm})=\exp\left(-\frac{1}{\hbar}\sum_{j=1+\frac{1}{2}(1\mp
1)}^{\ell+\frac{1}{2}(1\mp 1)}(n_{\pm})_{j,
    j \pm 1}\right), \qquad n_{\pm}\in N_{\pm}.
\ee
The proper analog of \eqref{equivprop} for \eqref{giv} is defined as
follows. In \cite{Lu2} Lusztig constructed a monoid $G^{>}$ of positive elements
for  an  arbitrary  reductive Lie group $G$
(recall that a monoid has  a multiplication operation and
a unity but an inverse element is not defined).
Consider the monoid of positive elements $GL^{>}_{\ell+1}(\IR)\subset
GL_{\ell+1}(\IR)$. Elements of $GL^{>}_{\ell+1}(\IR)$ can be
represented via the Gauss decomposition as $g=f\,h\,n$ where $n$ is
an upper-triangular positive matrix, $f$   is
a lower-triangular positive matrix and $h$ is a  diagonal matrix with
positive entries \cite{Lu2}.  The  monoid
$GL^{>}_{\ell+1}(\IR)$ naturally acts on the positive subset $N_+^{>}$
of the flag space $\CB_{\ell+1}$.
Indeed, the  monoid property implies that the product of two
positive elements  of $GL_{\ell+1}(\IR)$ is again a positive element.
The Gauss decomposition of the product
of $g\in GL_{\ell+1}^>(\IR)$ and $n\in N_+^>$ has the following form:
\be\label{Gdec}
n\cdot g=\tilde{f}\cdot \tilde{h}\cdot \tilde{n},
\ee
where $\tilde{n}\in N_+^>$, $\tilde{h}\in H^>$  and    $\tilde{f}\in N_-^>$.
We define  an action $GL^{>}_{\ell+1}(\IR)$  on the positive subset $N_+^{>}$
of the flag space $\CB_{\ell+1}$ by taking
$\tilde{n}$ as the result of the action of $g$ on $n\in N_+^>$.

Let us consider a monoid $GL^{>}_{\ell+1}$ over a ring
$\IR[\frak{z}_1,\ldots, \frak{z}_N]$ of nilpotent elements
$\frak{z}_1, \ldots, \frak{z}_N$ such that  non-diagonal
elements of the lower and upper
triangular parts in the  Gauss decomposition  $g=f\,h\,n$
are in nilpotent subalgebras of
$R^{\geq}_{nil}=\IR_{\geq}[\frak{z}_1,\ldots, \frak{z}_N]$ and the
diagonal elements of the Cartan part
are invertible elements of $R^{\geq}_{nil}$. Let ${}^{nil}GL_{\ell+1}^>$ be a submonoid of
$GL_{\ell+1}^>(R_{nil})$
such that all strictly lower-triangular elements $f_{i>j}$ of $f$ are
nilpotent i.e. $f_{ij}^M=0$ for some large $M\in \IZ_+$.

\begin{prop}\label{NilpProp}  The $\mathfrak{gl}_{\ell+1}$-Whittaker function
\eqref{giv} can be lifted to a function on the monoid
${}^{nil}GL^{>}_{\ell+1}$ such that the following functional  equation
holds:
\be\label{equiv}
\Psi_\lambda(fgn)=\psi^+(n)\psi^-(f)\Psi_\lambda(g),\qquad
g\in {}^{nil}GL^{>}_{\ell+1},\quad n\in {}^{nil}N^>_+, \quad
f\in {}^{nil}N^>_-,
\ee
where the functions $\psi^{\pm}$ are  given by \eqref{intchar}.
\end{prop}

\proof  The equivariance property \eqref{equiv} with respect to
${}^{nil}N_+^>$ follows from 
properties of the right Whittaker vector. Equivariance with
respect to ${}^{nil}N_-^>$ may be reduced to  equivariance with respect to
the Lie algebra and thus follows from  properties of the left Whittaker
vector. $\Box$

Now we are ready to provide a matrix element interpretation of  the
elementary $\mathfrak{gl}_{\ell+1}$-Whittaker functions \eqref{giv0}.
Let us start by recalling the definition of  a  tropical semifield
(see e.g. \cite{MS}, \cite{IMS}).
\begin{de} A tropical semifield $\CR$ is a set isomorphic to $\IR$ with the
following operations
\be\label{tropm}
\a\dot{\times}\b=\a+\b,\qquad  \a\dot{+}\b={\rm min}(\a,\b).
\ee
\end{de}

We also introduce the notation $\alpha\dot{/}\beta:=\alpha\dot{+}(-\beta)$.
The tropical semifield $\CR$ can be understood as a degeneration of the
standard semifield structure on the positive subset $\IR_+\subset \IR$ of
 real numbers. Indeed, consider the semifield $\IR^{(\hbar)}_+$ with the
following operations
\be\label{toper}
a\times_\hbar b=a\times b,\qquad  a+_\hbar b=(a^{\hbar}+b^{\hbar})^{1/\hbar}.
\ee
The semifield  $\IR^{(\hbar)}_+$ is  isomorphic to  $\IR_+^{(1)}=\IR_+$
via the map $a\rightarrow a^{\frac{1}{\hbar}}$. Let
$a=e^{-\alpha}$, $b=e^{-\beta}$, $\alpha,\beta\in \IR$. Then in the limit
$\hbar\rightarrow +\infty$ the operations \eqref{toper} are transformed
into the tropical operations \eqref{tropm}:
$$
\a\dot{\times}\b:=-
\lim_{\hbar\rightarrow +\infty}\log \left(e^{-\a}\times_\hbar e^{-\b}\right)=\a+\b,
$$
$$
\a\dot{+}\b = -\lim_{\hbar\rightarrow +\infty}\log
\left(e^{-\a}+_\hbar e^{-\b}\right) ={\rm min}(\a,\b),
$$
and the semifield $\IR_+^{(\hbar)}$ turns
into the  tropical semifield $\CR$.
The set of matrices ${\rm Mat}_{\ell+1}(\CR)$
has  a  monoid structure
arising in the limit $\hbar\to +\infty$ from the monoid structure on the subset of
positive elements $GL_{\ell+1}^>(\IR^{(\hbar)}_+)$.
Note that invertibility automatically holds in the tropical case and
thus in the following we can identify the matrix monoid ${\rm
Mat}_{\ell+1}(\CR)$ with  $GL_{\ell+1}(\CR)$.

Consider a  principal series representation $\pi^{(0)}_{\underline{\lambda}}$ of
$GL_{\ell+1}(\CR)$ in  $\CV_{\underline{\lambda}}={\rm
  Ind}_{B_-(\CR)}^{GL_{\ell+1}(\CR)}\,\chi^{(0)}_{\underline{\lambda}}$
of the monoid $GL_{\ell+1}(\CR)$ induced from a character
$\chi^{(0)}_{\underline{\lambda}}$ of the submonoid  $B_-(\CR)$
of upper triangular matrices
$$
\chi^{(0)}_{\underline{\lambda}}(b)=\prod_{j=1}^{\ell+1}\chi^{(0)}_{\lambda_j}(b_{jj}),
\qquad b\in B_-(\CR).
$$
where $\underline{\lambda}=(\lambda_1,\ldots ,\lambda_{\ell+1})$ and
$\chi^{(0)}_\lambda(\tau)=\exp -\imath \tau \lambda$, $\tau\in \IR$; $\lambda\in \IC$
is a multiplicative character of $\CR$
$$
\chi^{(0)}_\lambda(\tau_1\dot{\times}\tau_2)=\chi^{(0)}_\lambda(\tau_1)
\chi^{(0)}_\lambda(\tau_2).
$$

\begin{ex} The following formulas
$$
\pi^{(0)}_{\underline{\lambda}}\begin{pmatrix} 0 & \alpha  \\ 0 & 0 \end{pmatrix}\cdot
f(\tau)=f(\tau\dot{+}\a),
$$
\be\label{act2}
\pi^{(0)}_{\underline{\lambda}}\begin{pmatrix} \tau_1& 0 \\ 0 & \tau_2
\end{pmatrix}\cdot
f(\tau)=f(\tau\dot{\times}\tau_2\dot{/}\tau_1)\,
e^{-\imath(\lambda_1\tau_1+\lambda_2\tau_2)},
\ee
$$
\pi^{(0)}_{\underline{\lambda}}\begin{pmatrix} 0& 0 \\ \beta &0 \end{pmatrix}\cdot
f(\tau)=e^{-\imath(\lambda_1-\lambda_2)({\rm min}(0,\beta+\tau))}
f(-{\rm min}(-\tau,\beta)),
$$
define an  action of the tropical
monoid  $GL_2(\CR)$ on the space
$$
\CV_{\underline{\lambda}}={\rm Ind}_{B_-(\CR)}^{GL_2(\CR)}\,
\chi_{\underline{\lambda}}\,=\{
f(g)\in {\rm
  Fun}(GL_2(\CR))|f(b_-g)=\chi^{(0)}_{\underline{\lambda}}(b_-)f(g),\,\,
 b_-\in B_-(\CR)\},
$$
where
$$
\chi^{(0)}_{\underline{\lambda}}\left(\begin{pmatrix} \tau_1& 0 \\ \beta & \tau_2
\end{pmatrix}\right)=e^{-\imath(\tau_1\lambda_1+\lambda_2\tau_2)},
$$
and $\alpha$, $\beta$, $\tau_1$, $\tau_2$ $\in \CR$.
\end{ex}

The explicit formulas in  the example above can be obtained as a limit
of the analogous formulas for the  principal series representation
$\pi_{\underline{\lambda}}$ of
$GL_2^>(\IR)$ acting on $\CV_{\underline{\lambda}}={\rm
 Ind}_{B_-(R)}^{GL_2(R)} \chi^{(0)}_{\underline{\lambda}}$, where
$$
\chi_{\underline{\lambda}}(b)=|b_{11}|^{\imath\lambda_1/\hbar}\,|b_{22}|^{\imath\lambda_2/\hbar}.
$$
Thus we have
$$
\pi_{\underline{\lambda}}
\begin{pmatrix} 1&  e^{-\hbar \alpha}  \\ 0 &1\end{pmatrix}\cdot
f( e^{-\hbar
  \tau})=f(e^{-\hbar (\tau+\a)}),
$$
$$
\pi_{\underline{\lambda}}
\begin{pmatrix} e^{-\hbar \tau_1}& 0 \\ 0 & e^{-\hbar\tau_2}
\end{pmatrix}\cdot
f(e^{-\hbar \tau})=f(e^{-\hbar(\tau_2-\tau_1+\tau)})\,
e^{-\imath(\lambda_1\tau_1+\lambda_2\tau_2)},
$$
$$
\pi_{\underline{\lambda}}\begin{pmatrix} 1& 0 \\ e^{-\hbar \beta} &1 \end{pmatrix}\cdot
f(e^{-\hbar \tau})=
(1+e^{-\hbar(\beta+\tau)})^{\frac{\imath}{\hbar}(\lambda_1-\lambda_2)}
f\left(\frac{e^{-\hbar \tau}}{1+ e^{-\hbar (\tau+\beta)}}\right).
$$
Taking the limit $\hbar\to +\infty$ we obtain
$$
\pi^{(0)}_{\underline{\lambda}}\begin{pmatrix} 0 & \alpha  \\ 0 & 0 \end{pmatrix}\cdot
f(\tau)=f(\tau\dot{+}\a),
$$
$$
\pi^{(0)}_{\underline{\lambda}}\begin{pmatrix} \tau_1& 0 \\ 0 & \tau_2
\end{pmatrix}\cdot
f(\tau)=f(\tau_2-\tau_1+\tau)\, e^{-\imath(\lambda_1\tau_1+\lambda_2\tau_2)},
$$
$$
\pi^{(0)}_{\underline{\lambda}}\begin{pmatrix} 0& 0 \\ \beta &0 \end{pmatrix}\cdot
f(\tau)=e^{-\imath(\lambda_1-\lambda_2)({\rm min}(0,\beta+\tau))}
f(-({\rm min}(-\tau,\beta))),
$$
where in the last case we use the identity
$$
-{\rm min}(-\alpha,\beta)=\alpha-{\rm min}(0,\beta+\alpha).
$$
Thus we recover the formulas \eqref{act2} in the limit $\hbar\to
+\infty$.

Now we construct  elementary analogs of the
Whittaker vectors \eqref{Wclass} entering the construction of the
matrix element \eqref{Wf}. Precisely, we would like to construct
vectors $\psi^{(0)}_+(x)$ and $\psi^{(0)}_-(x)$  in
$\CV_{\underline{\lambda}}={\rm
  Ind}_{B_-(\CR)}^{GL_{\ell+1}(\CR)}\,\chi^{(0)}_{\underline{\lambda}}$
satisfying the functional relations
\be\label{elemW}
\psi^{(0)}_+(x\dot{\times}n_+)=\psi^+_0(n_+)\psi^{(0)}_+(x), \qquad
\psi^{(0)}_-(x\dot{\times}n_-)=\psi^-_0(n_-)\psi^{(0)}_-(x),
\qquad n_{\pm}\in N_{\pm}(\CR),
\ee
where
\be\label{tropchar}
\psi^{\pm}_0(n_{\pm})=\prod_{j=1+\frac{1}{2}(1\mp
  1)}^{\ell+\frac{1}{2}(1\mp 1)}
\psi_0(n_{j,j\pm 1}).
\ee
Here  $\psi_0(x)$ is an  additive character of $\CR$
\be\label{adChar}
\psi_0(x)=\Theta(x),\qquad
\psi_0(x_1\dot{+}x_2)=\psi_0(x_1)\,\psi_0(x_2),\qquad x_i\in \CR.
\ee

\begin{prop} The following expressions  for
Whittaker vectors hold:
\bqa\label{wR0}
\psi^{(0)}_+(\tau)
=\prod_{i=1}^{\ell}\prod_{k=1}^{\ell+1-i}\Theta(\tau_{k+i-1,i}-\tau_{k+i,i+1}),\eqa
\bqa\label{wL0}
\psi^{(0)}_-(\tau)= e^{\imath\sum_{k=1}^{\ell}\sum_{i=1}^k
(\lambda_k-\lambda_{k+1})(\tau_{k,i}-\tau_{\ell+1,i})}\,\prod_{i=1}^{\ell}
\prod_{k=1}^{\ell+1-i}\Theta(\tau_{k+i,k}-\tau_{k+i-1,k}).
\eqa
\end{prop}

\proof Recall that all the constructions of \cite{GKLO} are formulated
in terms of positive elements of $GL_{\ell+1}(\IR)$ and thus
can be reformulated in terms of $GL_{\ell+1}^{>}(\IR)$. This allows 
 the definition of the $\hbar\to +\infty$ limit  (see e.g. the derivation of
\eqref{act2} above) and so the obtaining of explicit expressions for Whittaker
vectors as the $\hbar\to \infty$ limit of the Whittaker vectors
constructed previously
(see eqs. (3.12), (3.14), (3.15) and (3.16) in \cite{GKLO})
\be\psi_+(T)=
  \exp\Big\{-\frac{1}{\hbar}\sum\limits_{i=1}^{\ell}\sum\limits_{n=1}^{\ell+1-i}
  e^{T_{n+i,i+1}-T_{n+i-1,i}}\Big\},
 \ee
  \be\psi_-(T)=
  \exp\Big\{\frac{\imath}{\hbar}\sum\limits_{k=1}^{\ell}\sum\limits_{i=1}^k
  (\lambda_k-\lambda_{k+1})
  T_{k,i}\Big\}\,
  \exp\Big\{-\frac{1}{\hbar}\sum\limits_{i=1}^{\ell}\sum\limits_{k=1}^{\ell+1-i}
  e^{T_{k+i-1,k}-T_{k+i,k}} \Big\}.
\ee
Using the variables $\tau_{i,j}(\tau)=\hbar^{-1}\tau_{ij}$
and taking the limit $\hbar\to \infty$, we obtain
\bqa
\psi^{(0)}_+(\tau)=\lim _{\hbar\to +\infty}\psi_+(T(\tau))=
\prod_{i=1}^{\ell}\prod_{k=1}^{\ell+1-i}\Theta(\tau_{k+i-1,i}-\tau_{k+i,i+1}),\eqa
and \bqa
\psi^{(0)}_-(\tau)= \lim _{\hbar\to +\infty}\psi_-(T(\tau))=
e^{\imath\sum_{k=1}^{\ell}\sum_{i=1}^k
(\lambda_k-\lambda_{k+1})(\tau_{k,i}-\tau_{\ell+1,i})}\,\prod_{i=1}^{\ell}
\prod_{k=1}^{\ell+1-i}\Theta(\tau_{k+i,k}-\tau_{k+i-1,k})\,.
\eqa
Here we use the following relation
$$
\lim_{\hbar \to +\infty}\,e^{-\frac{1}{\hbar}e^{-\hbar\tau}}=\Theta(\tau).
$$
$\Box$

\begin{ex} For $\ell=1$ the Whittaker vectors satisfying  the  equations
$$
\pi_{\underline{\lambda}}^{(0)}\begin{pmatrix} 0& \alpha \\ 0
  &0\end{pmatrix}\cdot\psi_+^{(0)}(\tau)=\Theta(\alpha)\psi^{(0)}_+(\tau),
$$
$$
 \pi_{\underline{\lambda}}^{(0)}\begin{pmatrix} 0& 0 \\ \beta &0\end{pmatrix}\cdot
\psi_-^{(0)}(\tau)=\Theta(\beta) \psi_-^{(0)}(\tau),
$$
are given by
$$
\psi_+^{(0)}(\tau)=\Theta(\tau),\qquad
\psi_-^{(0)}(\tau)=e^{-\imath \tau(\lambda_1-\lambda_2)}\Theta(-\tau).
$$
Indeed we have
$$
\pi_{\underline{\lambda}}^{(0)}\begin{pmatrix} 0 & \alpha  \\ 0 & 0 \end{pmatrix}\cdot
\psi_+^{(0)}(\tau)=\psi_+^{(0)}(\tau\dot{+}\a)=\Theta(\alpha)\psi_+^{(0)}(\tau)
,\qquad \tau\dot{+}\alpha={\rm min}(\tau,\alpha).
$$
For the Whittaker vector $\psi_-^{(0)}$ we have
$$
\pi_{\underline{\lambda}}^{(0)}\begin{pmatrix} 0& 0 \\ \beta &0\end{pmatrix}\cdot
\psi_-^{(0)}(\tau)=e^{-\imath 1(\lambda_1-\lambda_2)({\rm min}(0,\beta+\tau))}
e^{-\imath (\tau-({\rm min}(0,\beta+\tau)))(\lambda_1-\lambda_2)}
\Theta({\rm min}(-\tau,\beta)))=
\Theta(\beta)\psi_-^{(0)}(\tau).
$$
\end{ex}

Let us define an analog of the monoid ${}^{nil}GL_{\ell+1}^>(\IR)$
used in the formulation of Proposition \ref{NilpProp}. We define
elements of a submonoid ${}^{nil}GL_{\ell+1}(\CR)$ of
$GL_{\ell+1}(\CR)$ via the Gauss decomposition  $g=fhn$. Precisely we
impose the additional condition that $f\in N_-(\CR_+)$ where $\CR_+$
is modeled on $\IR_+$ with operations  $\dot{\times}$, $\dot{+}$.

\begin{te}\label{Theorem} {\it (i)}. From \eqref{giv0}, the elementary
$\mathfrak{gl}_{\ell+1}$-Whittaker function 
\be\label{giv01}
\Psi^{(0)}_{\underline{\lambda}}(\underline{x})\,=\,\int\limits_{\IR^{\ell(\ell+1)/2}}
\exp\left(\imath\sum\limits_{k=1}^{\ell+1}
\lambda_k \left(\sum_{i=1}^{k}\tau_{k,i}-\sum_{i=1}^{k-1}\tau_{k-1,i}\right)\right)\\
\times
\left(\prod_{i=1}^{\ell}\prod_{k=1}^{i}
\Theta\left(\tau_{i,k}-\tau_{i+1,k+1}\right)
\prod_{k=i}^{\ell}
\Theta\left(\tau_{i+1,k}-\tau_{i,k}\right)\right)
\prod_{k=1}^{\ell}\prod_{i=1}^{k}d\tau_{k,i},
\ee
with $x_i=\tau_{\ell+1,i},\,\,\,i=1,\ldots,\ell+1$ and $\tau_{k,i}=0$ for $i>k$,
allows  the following matrix element representation:
\be\label{MEtrop}
\Psi^{(0)}_{\underline{\lambda}}(\underline{x})=\<\psi_L,\,
\pi^{(0)}_{\underline{\lambda}}
(g(x))\,\psi_R\>,
\qquad g(x)={\rm diag}(x_1,\cdots, x_{\ell+1}),
\ee
where $\psi_L=\psi^{(0)}_+$ and $\psi_R=\psi^{(0)}_-$ are defined by
\eqref{elemW} and the pairing is given by
\be\label{pairNEW}
\<\psi_1,\psi_2\>=
\int_{\CR^{(\ell+1)\ell/2}}d\mu^{\times}(\tau)\,\overline{\psi}_1(\tau)\,
\psi_2(\tau).
\ee
Here
$$
d\mu^{\times}(\tau)=\prod_{k=1}^{\ell}\prod_{j=1}^kd\tau_{ki},
$$
is a product of multiplicative measures $d\mu^{\times}(\tau)=d\tau$ on $\CR$.

\noindent {\it (ii)}. The function \eqref{giv01}
can be naturally lifted to a function  on the  monoid
${}^{nil}GL_{\ell+1}(\CR)$ satisfying the functional relations
$$
\Psi^{(0)}_{\underline{\lambda}}(fgn)=\psi^+_0(n)\psi^-_0(f)
\Psi^{(0)}_{\underline{\lambda}}(g),\qquad
g\in {}^{nil}GL_{\ell+1}(\CR),\quad n\in N_+(\CR),\quad f\in N_-(\CR_+)
$$
where $\psi^{\pm}_0$ are given by  \eqref{tropchar}.
\end{te}

\proof The first part of the Theorem is a direct consequence of the
explicit formulas \eqref{wR0}, \eqref{wL0} for Whittaker vectors.
The second part of the Theorem follows from
Proposition \ref{NilpProp} by taking the limit $\hbar \to
\infty$. $\Box$

\begin{ex} Let us directly check the second part of Theorem
  \ref{Theorem} for the simplest
 nontrivial case of $\ell=1$.  The only property that is not obvious is
the  relation
$$
\<\psi_L,\pi_{\underline{\lambda}}^{(0)}\left(\begin{pmatrix} 0& 0\\ \beta
    & 0\end{pmatrix}\begin{pmatrix} \tau_1 & 0\\ 0
    &\tau_2\end{pmatrix}\right)
\,\psi_R\>=\Theta(\beta)\,\<\psi_L,\pi\left(\begin{pmatrix} \tau_1& 0\\ 0
    &\tau_2\end{pmatrix}\right)\,\psi_R\>, \qquad \beta>0.
$$
Using the integral representation for the pairing and the explicit
form of the Whittaker vectors we obtain
$$
\int_{\IR}d\tau\,e^{\imath(\lambda_1-\lambda_2)\tau}\Theta(-\tau)
\,\,\pi_{\underline{\lambda}}^{(0)}\left(\begin{pmatrix} 0& 0\\ \beta
    & 0\end{pmatrix}\right)
\,\,\Theta(\tau+\tau_2-\tau_1)
e^{-\imath(\lambda_1\tau_1+\lambda_2\tau_2)}.
$$
Taking into account the explicit form of the action
$$
\pi_{\underline{\lambda}}^{(0)}\begin{pmatrix} 0& 0 \\ \beta &0 \end{pmatrix}\cdot
f(\tau)=e^{-\imath(\lambda_1-\lambda_2)({\rm min}(0,\beta+\tau))}
f(-{\rm min}(-\tau,\beta)),
$$
we obtain the following integral:
$$
I=\int_{\IR}d\tau\,e^{\imath(\lambda_1-\lambda_2)\tau}\Theta(-\tau)
\,\,e^{-\imath(\lambda_1-\lambda_2)({\rm min}(0,\beta+\tau+\tau_2-\tau_1))}
\,\Theta(-{\rm min}(-\tau-\tau_2+\tau_1,\beta))
e^{-\imath(\lambda_1\tau_1+\lambda_2\tau_2)}.
$$
Let us represent this integral as the  sum
$$
I=I_1+I_2,
$$
with
$$
I_1=\int_{\tau<\tau_1-\tau_2-\beta}d\tau\,
e^{\imath(\lambda_1-\lambda_2)\tau}\Theta(-\tau)
\,\,e^{-\imath(\lambda_1-\lambda_2)({\rm min}(0,\beta+\tau+\tau_2-\tau_1))}
\,\Theta(-{\rm min}(-\tau-\tau_2+\tau_1,\beta))
e^{-\imath(\lambda_1\tau_1+\lambda_2\tau_2)},
$$
$$
I_2=\int_{\tau>\tau_1-\tau_2-\beta}d\tau\,
e^{\imath(\lambda_1-\lambda_2)\tau}\Theta(-\tau)
\,\,e^{-\imath(\lambda_1-\lambda_2)({\rm min}(0,\beta+\tau+\tau_2-\tau_1))}
\,\Theta(-{\rm min}(-\tau-\tau_2+\tau_1,\beta))
e^{-\imath(\lambda_1\tau_1+\lambda_2\tau_2)}.
$$
Taking into account the condition $\beta>0$ for the first integral
$I_1$ we have
$$
I_1=\int_{\tau<\tau_1-\tau_2-\beta}d\tau\,e^{\imath(\lambda_1-\lambda_2)\tau}\Theta(-\tau)
\,\,e^{-\imath(\lambda_1-\lambda_2)(\beta+\tau+\tau_2-\tau_1)}
\,\Theta(-\beta)e^{-\imath(\lambda_1\tau_1+\lambda_2\tau_2)}=0.
$$
The second integral $I_2$  is given by
$$
I_2=\int_{\tau>\tau_1-\tau_2-\beta}d\tau\,e^{\imath(\lambda_1-\lambda_2)\tau}\Theta(-\tau)
\,\,\Theta(\tau+\tau_2-\tau_1)e^{-\imath(\lambda_1\tau_1+\lambda_2\tau_2)}=
$$
$$
=\int_{\IR}d\tau\,e^{\imath(\lambda_1-\lambda_2)\tau}\Theta(-\tau)
\,\,\Theta(\tau+\tau_2-\tau_1)e^{-\imath(\lambda_1\tau_1+\lambda_2\tau_2)}
\Theta(\tau-\tau_1+\tau_2+\beta).
$$
Note that if $\tau-\tau_1+\tau_2>0$  then
obviously $\tau-\tau_1+\tau_2+\beta>0$  when $\beta>0$. Thus
we obtain
$$
I_2=\int_{\IR}d\tau\,e^{\imath(\lambda_1-\lambda_2)\tau}\Theta(-\tau)
\,\,\Theta(\tau+\tau_2-\tau_1)e^{-\imath(\lambda_1\tau_1+\lambda_2\tau_2)},
\qquad \beta >0.
$$
\end{ex}

To recapitulate: the  elementary $\mathfrak{gl}_{\ell+1}$-Whittaker
functions shall  be considered
as the $\mathfrak{gl}_{\ell+1}$-Whittaker functions over the tropical
semifield $\CR$.

\section{From $\IQ_p$ to $\IQ_1$}

In Section 4 the elementary analog of the  Whittaker
functions \eqref{levzeroW} defined in Section 3
 were obtained as a limit of the $q$-deformed Whittaker functions 
specialized at $q=0$.
In \cite{GLO3}, \cite{GLO4}, \cite{GLO5} it was  demonstrated that  the $q$-deformed
Whittaker functions  interpolate  the Whittaker functions
over the Archimedean field $\IR$ and the Whittaker functions over the
non-Archimedean fields $\IQ_p$, $p$ is a prime number.
The representation  \eqref{pW} can be understood as
an expression of  the non-Archimedean Whittaker functions
in terms of  characters of irreducible finite-dimensional
representations  of  the dual reductive  Lie algebra   and thus  is an instance of the
Shintani-Casselman-Shalika formula \cite{Sh}, \cite{CS}.
Note that the Shintani-Casselman-Shalika formula
provides  an explicit realization  of the local non-Archimedean Langlands
duality in terms of the  Whittaker functions.
The elementary Whittaker functions \eqref{levzeroW} are obtained by
further degeneration of the $q=0$ Whittaker functions. Equivalently
the elementary Whittaker functions can be obtained as the limit $p\to 1$  of
the non-Archimedean Whittaker functions.
Taking into account Theorem \ref{Theorem},
the  elementary $\mathfrak{gl}_{\ell+1}$-Whittaker function
\eqref{levzeroW} should be formally considered as the
$\mathfrak{gl}_{\ell+1}$-Whittaker function over a mysterious field which we may call $\IQ_1$.
 A similar relation holds between
elementary $L$-factors \eqref{lAL0}  and non-Archimedean local $L$-factors. Thus elementary
$L$-factors  shall be considered as local $L$-factors  corresponding
to $\IQ_1$. Below we will argue that these local $L$-factors
have a natural interpretation over the tropical semifield $\CR$
considered as a domain of the valuation map for $\IQ_1$.
In this Section we briefly discuss this interpretation for
the simplest $L$-functions, leaving detailed considerations for another occasion.

To explain the appearance of the tropical semifield $\CR$
in the $p\to1$ limit of constructions over $\IQ_p$,
we first recall  the  notion of a valuation on a non-Archimedean field.
A  valuation on a local non-Arhcimedean field $K$
is a map $\nu: \,\,K\to \IR$ such that
\be
\nu(x)=0\ \  \Leftrightarrow\ \  x=0,
\ee
\be
\nu(x\cdot y)=\nu(x)+\nu(y),
\ee
\be
\nu(x+y)\geq {\rm min}(\nu(x),\nu(y)).
\ee
These conditions can be succinctly summarized as follows. A
non-Arhcimedean valuation is a partial morphism  $\nu:K\to \CR$ of $K$
considered as a semifield (i.e taking into account only addition,
multiplication and division operations) to the tropical semifield
$\CR$. The term partial here means that the morphism property 
does not necessary hold for sums of elements of equal norms $\nu$. 
In the case of the p-adic field $\IQ_p$ a non-Archimedean
valuation can be defined as follows:
$$
\nu_p(p^na)=n,\qquad (p,a)=1.
$$
Thus the image of the partial morphism $\nu_p:\,\IQ_p\to \CR$ is a discrete subsemifield
$(\IZ,\dot{\times},\dot{+})\in \CR$.  Note that the valuation $\nu_p$
has a large   kernel $\IZ_p^*$  consisting of invertible p-adic integers.
Now we can partially clarify what meaning one can assign to a limit of the field
$\IQ_p$ when $p\to 1$. It is clear that the field $\IQ_1$, considered
as a semifield, should allow a partial epimorphism $\nu_1$ onto $\CR$. What the
kernel of $\nu_1$  is not quite  clear. Fortunately many
constructions over $\IQ_p$  can be reformulated directly in terms
of the semifield $(\IZ,\dot{\times},\dot{+})$
considered as an image domain of  the valuation map $\nu_p$. Hence  these
constructions over $\IQ_p$ allow a  limit $p\to 1$ formulated
in terms of the image of $\nu_1$ identified with  the tropical semifield $\CR$.
Below we illustrate this phenomena for the simple case of the local
$L$-factors.

Let us define a monoid which is  a non-Archimedean analog of the
monoid $GL_{\ell+1}^>(\IR)$ introduced in
\cite{Lu2} and discussed in the previous Section. Define on the
set $\CR_p=p^{\IZ}$ a monoid structure
with multiplication and  addition
\be\label{p-monid}
p^{n_1}\dot{\times}p^{n_2}=p^{n_1+n_2},\qquad
p^{n_1}\dot{+}p^{n_2}=p^{{\rm min}(n_1,n_2)}.
\ee
We consider the subset $GL_{\ell+1}(\CR_p)\subset GL_{\ell+1}(\IQ_p)$
consisting of matrices with all  entries of the form $p^n$, $n\in \IZ$.
The monoid structure on $GL_{\ell+1}(\CR_p)$ is defined using
\eqref{p-monid}. In the limit $p\to 1$ the monoid $GL_{\ell+1}(\CR_p)$
reduces to the tropical monoid $GL_{\ell+1}(\CR)$ considered in the
previous Section.

Local non-Archimedean $L$-factors  admit
 integral representations which can be rewritten as
sums over a set of points of the varieties defined  over $\CR_p$.
Below we illustrate how the elementary  local $L$-factors
arise  as $L$-functions over $\IQ_1$. Consider
an additive  character $\psi_p(x)$ of $\IQ_p$
$$
\psi_p(x+y)=\psi_p(x)\psi_p(y),
$$
given by a step-function (characteristic function of the
subset $\IZ_p\in \IQ_p$)
\be\label{additNA}
\psi_p(x)=\Theta(\nu_p(x)).
\ee
Fix a multiplicative character  $\chi_s(x)$ of $\IQ_p$
\be\label{multNA}
\chi_s(x\cdot y)=\chi_s(x)\,\chi_s(y),\qquad \chi_s(x)=p^{-s \nu(x)}.
\ee
The local non-Archimedean $L$-factor
$$
L_p(s)=\frac{1}{1-p^{-s}},
$$
has the following standard integral representation (see e.g. \cite{W}, \cite{MaP}):
\be\label{NALA}
L_p(s)=\int_{\IQ_p^\times}d\mu^\times(x) \psi_p(x)\chi_s(x)=
\sum_{n=0}^{\infty}p^{-ns}= \frac{1}{1-p^{-s}},
\ee
where $d\mu^\times$ is the multiplicative measure on $\IQ_p^\times$
for which ${\rm vol}(\IZ_p^\times)=1$.
Indeed, integrating over the fibres $\IZ_p^*$ of the projection $\IQ^*_p\to p^{\IZ}$,
one obtains a sum over  $\IZ$. Taking into account
the restriction  imposed by $\psi_p(x)$, the sum reduces to  a sum
over $\IZ_{\geq 0}$ and reproduces the right hand side of \eqref{NALA}.

It is easy to check that the $p\to 1$ limit of the integral
representation \eqref{NALA} recovers the integral representation
of the elementary local $L$-factors (see \eqref{LLimq} for $t=p^{-1}$,
$\lambda_j=0$ and  $\ell=0$). It is instructive to
consider the limit of the integral representation \eqref{NALA}.
The limits of the additive and
multiplicative characters \eqref{additNA}, \eqref{multNA}
are  multiplicative and additive characters of $\CR$
$$
\psi_0(\tau)=\Theta(\tau),\qquad \chi^{(0)}_s(\tau)=e^{-s\tau},
$$
where $\Theta(\tau)$ is the  Heaviside function
$$
\Theta(\tau)=1,\quad \tau\geq 0,\qquad \Theta(\tau)=0,\quad \tau<0.
$$
Thus for the integral representation of the elementary $L$-factor we obtain
\be\label{TropIN}
L_{\infty}(s)\,=\,\int_{\IR}\,d\tau\,\,\psi_0(\tau)\,\chi^{(0)}_s(\tau)=\frac{1}{s}.
\ee
This is precisely the elementary local Archimedean $L$-factor
\eqref{levzeroL} for $\ell=0$ and
$\lambda=0$.  To recapitulate, let us stress that we  start with the integral
representation (left hand side of \eqref{NALA}) formulated in terms of
$\IQ_p$ and  rewrite  it in terms of the image of $\IQ_p$ under the
valuation map $\nu_p$ (the sum representation in the right hand side of \eqref{NALA}). Then
we take the limit $p\to 1$ and obtain the 
 integral representation \eqref{TropIN} for elementary $L$-factor.
This integral can be also understood as an integral over tropical
semifield $\CR$ of the product of additive and multiplicative
characters of $\CR$. Thus the integral \eqref{TropIN} over tropical semifield  can be formally considered as
an analog for $\IQ_1$ of the intermediate step in \eqref{NALA}.

\vskip 1cm

\noindent {\small {\bf A.G.} {\sl Institute for Theoretical and
Experimental Physics, 117259, Moscow,  Russia; \hspace{8 cm}\,
\hphantom{xxx}  \hspace{2 mm} School of Mathematics, Trinity College
Dublin, Dublin 2, Ireland; \hspace{6 cm}\hspace{5 mm}\,
\hphantom{xxx}   \hspace{2 mm} Hamilton Mathematics Institute,
Trinity College Dublin, Dublin 2, Ireland;}\\
\hphantom{xxxx} {\it E-mail address}: {\tt anton@maths.tcd.ie}}\\

\noindent{\small {\bf D.L.} {\sl
 Institute for Theoretical and Experimental Physics,
117259, Moscow, Russia};\\
\hphantom{xxxx} {\it E-mail address}: {\tt lebedev.dm@gmail.com}}\\

\end{document}